\input cyracc.def
\input amssym.def
\input amssym.tex
  \overfullrule 0pt
\magnification=1150
\hsize=29pc
\font\sm=cmr10 at 9pt
 \font\smit= cmti10 at 9pt
\advance\vsize-14pt % -24pt for \makeheadline, + 10pt for \topskip

 \def\ord{{\rm \,ord\,}}
 \def\car{{\rm \,car\,}}
 
 \def\Spec{{\rm \,Spec\,}}
 
 \def\Spf{{\rm \,Spf\,}}
 
\def\loccit{$loc.$ $cit. \,$}
\def\cf{{\it cf. }}
\def\ie{{\it i.e. }}
\def\resp{{\it resp. }} \def\Z{{\bf Z}}
\def\A{{\bf A}}

\def\N{{\bf N}}
\def\Q{{\bf Q}}
\def\R{{\bf R}}
\def\Z{{\bf Z}}
\def\sp{{\rm sp}}

\def\sE{{\cal E}}
\def\sF{{\cal F}}

\def\sI{{\cal I}}

\def\sN{{\cal N}}
\def\sO{{\cal O}}

\def\sT{{\cal T}}
\def\sV{{\cal V}}
\def\limp{\mathop{{\lim\limits_{\displaystyle\leftarrow}}}\limits}
 \def\dem{\noindent {\bf D\'emonstration. }}
\def\sq{$\,{\square}$}

\centerline {\bf   Structure des con\-ne\-xions m\'ero\-morphes for\-melles de plusieurs variables  }
\medskip
\centerline {\bf   et semi-continuit\'e de l'irr\'egularit\'e. }
 \bigskip\bigskip
 \bigskip
  \centerline{ Yves Andr\'e}
  \bigskip\bigskip 

\bigskip 
\bigskip

{\sm \bigskip\bigskip 
\medskip 1. Introduction.
\par 2. Rang de Poincar\'e-Katz, d\'e\-com\-po\-si\-tion de Turrittin-Levelt, polyg\^one de Newton. 
\par 3. D\'e\-com\-po\-si\-tion de Turrittin-Levelt avec param\`etres, points tournants.
\par 4. \'Evitement des points tournants par \'ecla\-te\-ment.
  \par 5. Points de croisement, et stabilisation par \'ecla\-te\-ment.
  \par 6. Semi-continuit\'e du rang de Poincar\'e-Katz.
  \par 7. Semi-continuit\'e de l'irr\'egularit\'e.
 \par 8. Autour de la conjecture de Sabbah.
\smallskip Appendice: sp\'ecialisation du polyg\^one de Newton. }
   
 \bigskip 
\bigskip

\bigskip 
\bigskip
\bigskip {\smit Abstract. $-$ }{\sm  We prove Malgrange's conjecture on the absence of confluence phenomena for integrable meromorphic connections. More precisely, if $ Y\to X$ is a complex-analytic fibration by smooth curves, $Z$ a hypersurface of $Y$ finite over $X$, and $\nabla$ an integrable meromorphic connection on $Y$ with poles along $Z$, then the function which attaches to  {\smit x} $\in X$ the sum of the irregularities of the fiber $\nabla_{(x)}$ at the points of $Z_x$ is lower semicontinuous.

The proof relies upon a study of the formal structure of integrable meromorphic connections in several variables.}

\vfill\eject  

\bigskip\bigskip  {\bf 1. Introduction.}     
  
\bigskip\noindent {$\bf  1.1.$} Le r\'esultat principal de cet article a trait aux probl\`emes de confluence, c'est-\`a-dire \`a la variation de l'irr\'egularit\'e dans une famille analytique d'\'equations diff\'erentielles lin\'eaires \`a singularit\'es m\'eromorphes.
  
 \smallskip Soit  $f: Y\to X$ un morphisme lisse de dimension relative $1$, \`a fibres connexes, entre vari\'et\'es analytiques complexes, et soit $Z\subset Y$ une hypersurface de $Y$ finie sur $X$. Notons $\sO_Y(\ast Z)$ le faisceau des germes de fonctions m\'eromorphes sur $Y$ \`a p\^oles le long de $Z$, et $\,\Omega^1_{Y }(\ast Z)\,$ (\resp $\,\Omega^1_{Y/X}(\ast Z)\,$) le faisceau des $1$-formes diff\'erentielles m\'eromorphes sur $Y$ (\resp relatives) \`a p\^oles le long de $Z$. 
  
  Soit $\sN$ un $\sO_Y(\ast Z)$-module localement libre de rang $\mu$ muni d'une con\-ne\-xion relative
  $$\nabla: \sN \to \sN \otimes_{\sO_Y(\ast Z)}\Omega^1_{Y/X}(\ast Z).$$    
   Pour $x\in X$, la fibre $\sN_{(x)}$ sur la courbe $f^{-1}(x)$ est ipso facto munie d'une con\-ne\-xion m\'eromorphe $\nabla_{(x)}$ dont les singularit\'es sont des p\^oles situ\'es aux points $\,z\in Z_x = Z\cap f^{-1}(x)$.
   
   \smallskip   L'un des invariants fondamentaux dans l'\'etude asymptotique de $\nabla_{(x)}$ au voisinage d'un point $z\in Z_x$ est son {\it irr\'egularit\'e} (au sens de Deligne-Komatsu-Malgrange), not\'ee $\,{\rm ir}_z\,\nabla_{(x)}$. Nous nous int\'eressons \`a la fonction sur $X$
   $$i(\nabla, x) =\sum_{z\in Z_x} \, {\rm ir}_z\, \nabla_{(x)}.$$ Il est facile de voir que cette fonction
est constructible, et P. Deligne a d\'emontr\'e, en utilisant la description de G\'erard-Levelt de l'irr\'egularit\'e, que la fonction voisine 
$$ i(\nabla, x)\, + \, \mu.\,\sharp\,Z_x$$ est semi-continue inf\'erieurement ([${\bf De2}$], 1976).
C'est l\`a un r\'esultat optimal, comme le montrent les ph\'enom\`enes bien connus de ``confluence" o\`u $i(\nabla, x)$ cro\^{\i}t aux points de branchement de $Z/X$. 

 \smallskip Au d\'ebut des ann\'ees $80$, B. Malgrange a conjectur\'e que lorsque $\nabla$ provient d'une con\-ne\-xion absolue {\it int\'egrable}
 $$ \sN \to \sN \otimes_{\sO_Y(\ast Z)}\Omega^1_{Y }(\ast Z), $$  la fonction $i(\nabla, x)$ elle-m\^eme {\it est semi-continue inf\'erieurement sur $X$}: il n'y a {\it pas de ph\'enom\`ene de confluence sous la condition d'int\'egrabilit\'e}.
 
\smallskip Dans cet article, nous d\'emontrons cette conjecture (\cf 7.1.2), ainsi que plusieurs variantes (6.1.1, 6.1.3, 7.1.1).

\bigskip\noindent {$\bf 1.2.$}  Parmi les travaux li\'es \`a cette conjecture (\cf {e. g.} [$\bf Me1$] [$\bf Me2$, 4] [$\bf S2$, I.3]), nous nous inspirons de celui de C. Sabbah, qui rapproche, dans un contexte microlocal, les probl\`emes de confluence du probl\`eme de la structure formelle des con\-ne\-xions m\'ero\-morphes in\-t\'e\-gra\-bles en deux variables.
  
  Apr\`es les travaux pionniers de M. Hukuhara puis H. Turrittin, le point de d\'epart de l'\'etude asymptotique des con\-ne\-xions m\'ero\-morphes d'une variable est le ``th\'eor\`eme de Turrittin-Levelt" qui affirme l'exis\-tence, apr\`es ramification aux p\^oles, d'une d\'e\-com\-po\-si\-tion canonique de la con\-ne\-xion for\-melle associ\'ee en con\-ne\-xions \'el\'ementaires (produits tensoriels d'une con\-ne\-xion de rang un et d'une con\-ne\-xion r\'egu\-li\`ere).

  L'analogue de ce r\'esultat pour des con\-ne\-xions m\'ero\-morphes in\-t\'e\-gra\-bles en plusieurs variables est loin d'aller de soi.  Apr\`es une s\'erie de r\'esultats partiels (voir par exemple [$\bf LvdE$], [$\bf M$, app.]), la premi\`ere \'etude syst\'ematique du probl\`eme est due \`a Sabbah, dans le cas de deux variables [$\bf S1$, $\bf S2$]. Il introduit la notion de bonne structure for\-melle: une con\-ne\-xion m\'ero\-morphe in\-t\'e\-gra\-ble ${\nabla}$ sur une surface analytique, \`a p\^oles le long d'un diviseur \`a croisements normaux $Z$, admet une bonne structure for\-melle le long de $Z$ si localement, apr\`es ramification le long des branches de $Z$,  le formalis\'e de ${\nabla}$ tant aux points de croisement que le long des branches de $Z$ (ind\'ependamment) se d\'ecompose en con\-ne\-xions \'el\'ementaires (avec une condition technique que nous n\'egligeons ici, \cf 5.5 ci-dessous).  
  
  Sabbah conjecture l'existence d'une bonne structure for\-melle apr\`es \'eclate\-ments, et montre, en utilisant les cycles microcarac\-t\'e\-ris\-tiques d'Y. Laurent,  comment cela entra\^{\i}ne la conjecture  de Malgrange. 
 
    \medskip\noindent {$\bf 1.3.$}  Dans cet article, nous irons assez loin en direction de la conjecture de Sabbah pour pouvoir adapter inconditionnellement les contours de cette strat\'egie de d\'emonstration de la conjecture de Malgrange. Cela passe par une analyse d\'etaill\'ee des ``points tournants" du diviseur polaire (\S\S 3,4) et de la structure formelle de $\nabla$ aux croisements de ce diviseur (\S 5). 
    
    Plus pr\'ecis\'ement, \'etant donn\'e une con\-ne\-xion m\'ero\-morphe in\-t\'e\-gra\-ble ${\nabla}$ sur un germe de surface $(Y, Q)$ \`a p\^oles le long d'un diviseur r\'eduit  $Z$, il s'agit d'estimer l'irr\'egularit\'e en $Q$ de la restriction de ${\nabla}$ \`a un germe de courbe donn\'e $(C, Q)$ (non contenu dans $Z$) en termes de l'irr\'egularit\'e de ${\nabla}$ le long des branches de $Z$. Pour cela, nous proc\'edons \`a une s\'eries d'\'eclatements pour nous ramener \`a une situation de croisements normaux et o\`u la structure formelle de $\nabla$ est suffisamment bonne aux croisements. Nous contr\^olons cette structure formelle le long de chaque composante $E_i$ du diviseur exceptionnel au moyen d'un certain diviseur $D_i$ sur le fibr\'e normal de $E_i$ dans l'\'eclat\'e de $Y$, d\'efini de mani\`ere \'el\'ementaire \`a l'aide de la d\'ecomposition de Turrittin-Levelt relative (3.6). Nous concluons par le biais d'une analyse g\'eom\'etrico-combinatoire de ces diviseurs $D_i$ et $E_i$ (pour une explication heuristique de laquelle nous renvoyons \`a 7.2.1). 
    
  \smallskip  Notre approche est purement alg\'ebrico-formelle, dans l'esprit de [$\bf AB$], et n'utilise pas la th\'eorie microlocale.
 
    \medskip\noindent {$\bf 1.4.$} Ce point de vue nous permet notamment de donner une preuve purement alg\'ebrique du th\'eor\`eme bien connu de Deligne selon lequel la sp\'ecialisation \`a toute courbe lisse d'une con\-ne\-xion int\'egrable alg\'ebrique r\'eguli\`ere est encore r\'eguli\`ere (6.1.2). Nous renvoyons \`a [$\bf B$] pour l'analyse des difficult\'es auxquelles on est confront\'e en confinant le probl\`eme dans le contexte r\'egulier (faillite des techniques alg\'ebriques de r\'eseaux logarithmiques); c'est en nous pla\c cant r\'esolument en si\-tua\-tion  irr\'eguli\`ere que nous  r\'eglons le probl\`eme.
    
     \medskip\noindent {$\bf 1.5.$}  La plupart des techniques de cet article valent aussi en carac\-t\'e\-ris\-tique $p>0$ pour des modules \`a con\-ne\-xion de rang $<p$ (c'est par exemple le cas du th\'eor\`eme de Turrittin-Levelt). Ce souci d'\'eviter des restrictions inutiles sur la carac\-t\'e\-ris\-tique, loin de compliquer les preuves, conduit parfois \`a en trouver de plus simples que celles qu'offre la litt\'erature (\cf e. g. 2.5.1, 5.3.2, 5.4.2). Mais il a surtout pour origine l'analogie $p$-adique suivante.
     
       On peut voir la conjecture de Crew  sur la r\'eduction semi-stable des $F$-isocristaux surconvergents \`a une variable (th\'eor\`eme de monodromie $p$-adique) comme un analogue $p$-adique - plus profond - du th\'eor\`eme de Turrittin-Levelt (\cf [$\bf A2$, 3] pour des pr\'ecisions sur cette analogie).  Sa g\'en\'eralisation \`a plusieurs variables  conjectur\'ee par A. Shiho actuellement tr\`es \'etudi\'ee (\cf  [$\bf Ke$]) appara\^{\i}t, de m\^eme,  comme un analogue $p$-adique de la conjecture de Sabbah. On peut donc s'attendre \`a ce que celle-ci soit plus abordable que celle-l\`a.
          
     \medskip\noindent {$\bf 1.6.$}   Cette analogie incite par ailleurs \`a reformuler le probl\`eme de la bonne structure formelle apr\`es \'eclatements en termes rigides-analytiques $x$-adiques. Une telle traduction sugg\`ere alors de tirer profit de la compacit\'e des espaces analytiques de Berkovich associ\'es en vue de la conjecture de Sabbah. L'analyse montre que les points de Berkovich qui posent probl\`eme sont ceux \`a valeurs dans une extension imm\'ediate de la compl\'etion  de la cl\^oture alg\'ebrique de $k[[x]]$.    
   
   Nous donnons n\'eanmoins, de ce point de vue rigide-analytique, un r\'esultat sur la structure formelle des connexions int\'egrables qui g\'en\'eralise en plusieurs variables le th\'eor\`eme de Turrittin-Levelt (8.2.2).
        
     \medskip\noindent {$\bf 1.7.$}  Dans l'appendice, nous nous pla\c cons de nouveau dans la situation  1.1 d'une connexion relative, mais en supposant que le diviseur polaire $Z$ est \'etale sur la base. Nous montrons dans ce cas que le polyg\^one de Newton, plac\'e conventionnellement dans le quatri\`eme quadrant, ne peut que cro\^{\i}tre (au sens large) par sp\'ecialisation.

       \bigskip    \bigskip {\smit Remerciements. $-$ }{\sm  Cet article doit beaucoup aux travaux mentionn\'es de C. Sabbah sur l'analyse asymptotique en deux variables. Je le remercie pour ses commentaires \'eclairants sur ces travaux. Je remercie aussi J. F. Mattei pour ses remarques critiques qui m'ont amen\'e \`a corriger une premi\`ere version de ce texte sur un point important.}

 \bigskip \bigskip  
\vfill\eject

\bigskip {\bf 2. Rang de Poincar\'e-Katz, d\'e\-com\-po\-si\-tion de Turrittin-Levelt, polyg\^one de Newton.}

\medskip On peut trouver davantage de d\'etails sur ces pr\'eliminaires dans [$\bf{AB}$, II].

\bigskip\noindent {\bf 2.1. Rang de Poincar\'e-Katz et th\'eor\`eme de Turrittin-Katz ``abstrait".} Soit ${{F}}$ un corps complet pour une valeur absolue non archim\'edienne $\vert \;\vert$. 

Soit ${{M}}$ un ${{F}}$-espace vectoriel de dimension finie $\mu$, muni d'une norme $\vert\;\vert_{{{M}}}$\footnote{$^{(1)}$}{\sm dans le cas d'une alg\`ebre de matrices, nous prendrons  pour  $\vert\;\vert_{{{M}}}$ la norme-sup.}. Soit $D$ un endomorphisme du sous-groupe topologique sous-jacent \`a ${{{M}}}$. La norme spectrale de $D$ est 
$$\displaystyle{ \vert D \vert_{{{{M}}},\sp}  = \lim_n  \, \vert D^n\vert_{{{M}}}^{1\over n}= \inf_n  \, \vert D^n\vert_{{{M}}}^{1\over n} \;\in \R_+.}$$  Elle ne d\'epend pas du choix de la norme sur ${{M}}$ (ces normes sont toutes \'equiva\-lentes). Prendre garde qu'en d\'epit de son nom, la norme spectrale ne d\'efinit pas une norme sur $End_\Z^{cont} {{M}}$ mais, a priori, seulement sur chacune de ses sous-alg\`ebres commutatives.

 Nous nous int\'eresserons au cas o\`u $({{M}}, \, D= {{\nabla}}(\partial))$ est un {\it module diff\'erentiel} relativement \`a une d\'erivation continue $\partial $ de ${{F}}$, c'est-\`a-dire un ${{F}}$-espace vectoriel de dimension finie muni d'une action additive continue $ {{\nabla}}(\partial)$ de $\partial$ v\'erifiant la r\`egle de Leibniz
 $$\forall f\in \sF, \forall m\in {{M}},\;{{\nabla}}(\partial)(fm)= \partial f. m+ f{{\nabla}}(\partial)(m) .$$
 Dans ce cas, bien que ${{\nabla}}(\partial)$ et l'op\'erateur de multiplication par $f\in {{F}}$ ne commutent pas, la r\`egle de Leibniz entra\^{\i}ne n\'eanmoins que $\vert\;\vert_{{{{M}}},\sp}$ d\'efinit une norme sur la sous-${{F}}$-alg\`ebre de $End_\Z^{cont} {{M}}$ engendr\'ee par ${{\nabla}}(\partial)$.  
 
  Nous supposerons toujours que $$\vert \partial\vert_{{{F}}}=1.$$
  
 \noindent  Il est souvent plus commode de parler en termes de valuation plut\^ot que de valeur absolue: $v(f)= -\log \vert f\vert \,\in \R\cup \{\infty\}$. On note 
 $v_{{{{M}}},\sp}({{\nabla}}(\partial)) = -\log \vert {{\nabla}}(\partial) \vert_{{{{M}}},\sp} $ la {\it valuation spectrale} de ${{\nabla}}(\partial)$.

\proclaim D\'efinition 2.1.1. Le {\rm rang de Poincar\'e-Katz}\footnote{$^{(2)}$}{ \rm \sm on dit aussi ``plus grande pente", voir 2.4.} du module diff\'erentiel  $({{M}}, {{\nabla}}(\partial))$ est 
$$\rho_v({{M}}, {{\nabla}}(\partial)) := \max (0, - v_{{{{M}}},\sp}({{\nabla}}(\partial)) ).$$
On dit que $({{M}}, {{\nabla}}(\partial))$ est {\rm r\'egu\-lier} si son rang de Poincar\'e-Katz est nul. 
 \par
\medskip\noindent Nous \'ecrirons souvent $\rho_v({{\nabla}}(\partial)) $ ou $\rho_v({{M}}) $ pour abr\'eger.  Il ne change pas si l'on remplace $\partial$ par $f.\partial$ avec $v(f)=0$.

\noindent Par ailleurs, on d\'eduit des propri\'et\'es g\'en\'erales des valuations spectrales le comportement du rang de Poincar\'e-Katz par dualit\'e, somme directe et produit tensoriel:
$$ \rho_v({{M}}^\vee) = \rho_v({{M}}),\; \rho_v({{M}}\oplus {{M}}')= \max (\rho_v({{M}}), \rho_v({{M}}')),$$ $$ \,\rho_v({{M}}\otimes {{M}}')\leq \max (\rho_v({{M}}), \rho_v({{M}}')) $$ avec \'egalit\'e si $ \rho_v({{M}})\neq \rho_v({{M}}').$ 

\medskip \proclaim Th\'eor\`eme 2.1.2. {\rm (Turrittin-Katz)}. Supposons que $({{M}}, {{\nabla}}(\partial))$ admette une base cyclique, c'est-\`a-dire une base de ${{M}}$ de la forme $${\bf m}\,=\,(m, {{\nabla}}(\partial)(m), \ldots, {{\nabla}}(\partial)^{\mu -1}(m)).$$  
Pour tout $n\in \bf N$, notons $G_n$ la matrice de ${{\nabla}}(\partial)^n$ dans cette base, de sorte que $G_1$ s'\'ecrit sous la forme $ \pmatrix{0&& 0&   \theta_0\cr
    1 &&0&  \theta_1\cr &\ldots&&\cr 
     &&& \cr 
     0  &&1&      \theta_{\mu-1} } $. 
     \medskip 1) Pour tout $ \sigma \geq 1$, les \'enonc\'es suivants sont \'equivalents:
 \medskip\noindent $i)$ $\vert   {{\nabla}}(\partial)  \vert_{{{{M}}},sp}   \leq \sigma ,$ 
\medskip\noindent $ii)$ pour tout $j=0,\dots, \mu-1$, $\,\vert\theta_j\vert \leq \sigma^{\mu-j},$
\medskip\noindent $iii)$ pour tout $n\in \bf N$, $\,\vert G_n\vert \leq \sigma^{n+\mu-1}$.
     \medskip 2) Supposons qu'il existe $\xi\in {{F}}$ de valeur absolue $\min(1,\min \vert \theta_j\vert^{-1/\mu-j})$. Ainsi la matrice de $\xi.{{\nabla}}(\partial)$ dans la base $${{\bf m}'} =  {\bf m}\pmatrix{1&0&&0\cr 0& \xi &&0\cr &&\ldots&\cr 0&0&& \xi^{\mu-1 }} $$ est de norme $\leq 1$. Si $ \vert   {{\nabla}}(\partial)  \vert_{{{{M}}},sp}>1, $ alors la r\'eduction de cette matrice dans le corps r\'esiduel de ${{F}}$ n'est {\rm pas nilpotente}.
     \par
 
 Voir [$\bf{CD}$, 1.5]\footnote{$^{(3)}$}{\sm et aussi [$\bf{T}$], [$\bf{Ka}$] en \'egale carac\-t\'e\-ris\-tique. En in\'egale carac\-t\'e\-ris\-tique, c'est-\`a-dire dans le cas $p$-adique, le point 1) est d\^u \`a Dwork et Young; dans ce cas, $\vert   {{\nabla}}(\partial)  \vert_{{{{M}}},sp}$ n'est autre que l'inverse du rayon de convergence g\'en\'erique, multipli\'e par $\vert p\vert^{1/p-1}$. }.
 Le point 1) donne le calcul du rang de Poincar\'e-Katz dans une base cyclique:
 
   \proclaim Corollaire 2.1.3.   $i)$  $\displaystyle\;  \, \rho_v( {{\nabla}}(\partial))  =  \max (0,   \max_{j=0,\dots,\mu-1}  \, {{-v(\theta_j)}\over {\mu -j}})  $; en par\-ti\-cu\-lier, si $v({{F}}^\ast)=\Z$, c'est un nombre rationnel de d\'enominateur $\leq \mu$.  
 \medskip\noindent $ii)$ La fonction $\, n\in {\bf N}\mapsto \log \vert  {{\nabla}}(\partial)^n\vert_{{{M}} }  - n \cdot \rho_v(  \nabla(\partial))\,$  est born\'ee; en par\-ti\-cu\-lier, elle ne prend qu'un nombre fini de valeurs, toutes rationnelles, si  $v({{F}}^\ast)=\Z$ et si $  \vert \;\vert_{{{M}} }$ est la norme-sup relative \`a une base quelconque. \par
   
     \proclaim Corollaire 2.1.4.  Si ${{M}}$ est r\'egu\-lier, alors la matrice de $ {{\nabla}} $ dans toute base cyclique est de norme $\leq 1$. \par
 
    \medskip\noindent {\bf Remarques 2.1.5.} 1) Il est bien connu que $({{M}}, {{\nabla}}(\partial))$ admet une base cyclique si $\car {{F}}= 0$ ou bien si $\mu \leq \car {{F}}  $ (par la r\`egle de Leibniz, le r\'esultat ne change pas si l'on remplace $\partial $ par un multiple non nul, ce qui permet de supposer l'existence d'un $x$ tel que $\partial(x)=1$; la preuve de [$\bf{DGS}$, III.4.2] s'applique alors). 
    
  \smallskip  \noindent 2) Il d\'ecoule du point $i)$ du corollaire que $({{M}}, {{\nabla}}(\partial))$ est r\'egu\-lier si et seulement si il admet une base dans laquelle la matrice de $ {{\nabla}}(\partial)$ est de norme $\leq 1$.
   
    \medskip\noindent {\bf Exemples 2.1.6.} 1) Soit ${{K}}$ un corps et prenons pour ${{F}} $ le corps de s\'eries for\-melles ${{K}}((x))  $ (\resp ${{F}} = {{K}}((x_2))((x_1))\,$) muni de la valuation $x$-adique $\ord_x$ (\resp $x_1$-adique), et pour $\partial$ la d\'erivation $ x {d\over dx }$ (\resp ${{\partial \over {\partial x_2}}}$). Alors $\vert \partial \vert_{{{F}}} =1$, et $\vert \partial \vert_{{{F}}, \sp} =1$ ou $0$ selon que $\car K= 0$ ou non.
    
    \smallskip    \noindent 2) Prenons pour ${{K}}$ un corps $p$-adique, pour ${{F}}$ le compl\'et\'e $p$-adique de  $F(x)$ (ou bien le compl\'et\'e  de $\sO_K((x))\otimes_{\sO_K} K$), et $\partial={d\over dx}$. Alors $\vert \partial \vert_{{{F}}} =1$, et $\vert \partial \vert_{{{F}}, \sp} =\vert p\vert^{1/(p-1)}$.   
   
  Comme me l'a fait remarquer F. Baldassarri, il serait plus naturel dans cette situation de consid\'erer, plut\^ot que $\max (0, - v_{{{{M}}},\sp}({{\nabla}}(\partial)) )$, la quantit\'e 
 $$\sup_D \, \max (0, - v_{{{{M}}},\sp}({{\nabla}}(D)) ) \leqno{(2.1.6.1)}$$ o\`u $D$ parcourt les op\'erateurs diff\'erentiels $K$-lin\'eaires de norme $p$-adique $1$ (sur $F$), quantit\'e qui co\"{\i}ncide avec le rayon de convergence g\'en\'erique (valuatif); par exemple, si $\vert \pi\vert = \vert p\vert^{1/(p-1)}$, elle distingue entre le module diff\'erentiel $ (F, \nabla( {d\over dx})(1)$ $= 1)$ (pour lequel elle vaut $1$), et $ (F, \nabla( {d\over dx})(1)= \pi)$ (pour lequel elle est nulle).

   \medskip\noindent {\bf 2.2. Lemme de d\'e\-com\-po\-si\-tion.} 
  Des lemmes de d\'e\-com\-po\-si\-tion des con\-ne\-xions for\-melles ``suivant les valeurs propres" apparaissent sous de nombreuses formes dans la litt\'erature depuis Turrittin. En voici une version tr\`es g\'en\'erale due \`a Levelt et van den Essen  [$\bf{LvdE}$].
 
\proclaim Lemme 2.2.1. Soit ${ R}$ un anneau local noeth\'erien complet, d'id\'eal maximal  $\frak{m}$ et de corps r\'esiduel $k$. Soit $\delta$ une d\'erivation de ${R}$ telle que $\delta({ R})\subset \frak{m}$ et $\delta(\frak{m})\subset \frak{m}^2$. Soit ${\Bbb M}$ un ${ R}$-module libre\footnote{$^{(4)}$}{\rm\sm l'hypoth\`ese de libert\'e n'est m\^eme pas impos\'ee dans [$\bf{LvdE}$], mais elle permet de simplifier sensiblement la preuve.}  de type fini, muni d'une action additive $ \nabla(\delta)$ de $\delta$ v\'erifiant la r\`egle de Leibniz, et notons  ${\bar \nabla(\delta)  }$ l'action $k$-lin\'eaire induite par $ \nabla(\delta)$ sur $\bar {\Bbb M}:= {\Bbb M}\otimes_{ R} k$. 
 Supposons que $\bar {\Bbb M}$ se d\'ecompose sur $k$ en somme directe  $${\bar {\Bbb M} }=   \oplus \, { \bar {\Bbb M}}_j$$ de sous-espaces ca\-rac\-t\'e\-ris\-ti\-ques.  Alors:
  \medskip\noindent  $1)$ cette d\'e\-com\-po\-si\-tion se rel\`eve de mani\`ere unique en une d\'e\-com\-po\-si\-tion 
 $$  {\Bbb M} =   \oplus \,  {\Bbb M}_j$$ en sous-${ R}$-modules stables sous $ \nabla(\delta)$.
   \medskip\noindent $2)$ Soit $\delta'$ est une autre d\'erivation continue de ${ R}$ commutant \`a $\delta$, et soit  $\nabla(\delta')$ une application additive ${\Bbb M}\to  {\Bbb M} $  v\'erifiant la r\`egle de Leibniz vis-\`a-vis de $\delta'$, et  commutant \`a $ \nabla(\delta)$. Alors la d\'e\-com\-po\-si\-tion pr\'ec\'edente est stable sous $ \nabla(\delta')$. 
      \medskip\noindent $3)$ La m\^eme conclusion vaut si, au lieu de supposer que $\delta$ et $\delta'$ (\resp $\nabla (\delta)$ et $\nabla (\delta')$) commutent, on suppose seulement l'existence de deux non-diviseurs de z\'ero $g, g'\in  R$ tels que ${1\over g}\delta$ et ${1\over g'} \delta'$ (\resp ${1\over g} \nabla(\delta)$ et ${1\over g'} \nabla (\delta')$) commutent, et que $\delta'(g)/g\in \frak m$.  \par 
   
   L'assertion 3) ne figure pas dans \loccit, mais se d\'emontre par le m\^eme argument que 2): le point est que puisque $\delta'(g)/g, \, \delta(g')/g' \, \in \frak m$, ${\bar \nabla(\delta)  }$ et ${\bar \nabla(\delta')  }$ commutent.
   
  \medskip Pour comparer cette si\-tua\-tion \`a celle obtenue par changement de base injectif (non n\'ecessairement local) $ R\to  R'$, o\`u ${ R'}$ est un autre anneau local noeth\'erien complet de corps r\'esiduel $k'$, on a:
 
 \proclaim Compl\'ement 2.2.2. Supposons que $\delta$ (\resp $ \nabla(\delta)$) s'\'etende \`a $ R'$ (\resp \`a $\,{\Bbb M}'={\Bbb M}\otimes_{ R} R'$) avec les m\^emes propri\'et\'es. Alors la d\'e\-com\-po\-si\-tion canonique $\,{\Bbb M}'= \, \oplus  \, {\Bbb M}'_i\;$ qui rel\`eve la d\'e\-com\-po\-si\-tion en sous-espaces ca\-rac\-t\'e\-ris\-ti\-ques de  $\,\bar {\Bbb M}' := {\Bbb M}'\otimes_{ R'} k'$ raffine la d\'e\-com\-po\-si\-tion $\;  {\Bbb M}' = \,  \oplus \,  {\Bbb M}_j\otimes_{ R} R'$.
  \par
 
 La d\'emonstration, directe, est laiss\'ee au lecteur. 
  
\bigskip\noindent {\bf 2.3. Th\'eor\`eme de Turrittin-Levelt ``abstrait".} 
 On reprend la si\-tua\-tion de 2.1 sous des hypoth\`eses suppl\'ementaires: 
 
\medskip \noindent - on suppose que ${{F}}$ est complet pour une valuation {\it discr\`ete} $v$, qu'on normalise de sorte que   
 $$v({{F}}^\ast)=\Z.$$
 Elle s'\'etend de mani\`ere unique \`a toute extension finie ${{F}}' $ de ${{F}}$, et il y a lieu de multiplier la valuation \'etendue par l'indice de ramification $e$ de l'extension pour pr\'eserver la condition de normalisation. On note $v'$ la valuation normalis\'ee sur $\sF'$ (de sorte que $v'(({{F}}^\ast)=e\Z $).
  
Notons par ailleurs que la d\'erivation (de norme $1$) $\partial$  s'\'etend de mani\`ere unique \`a toute extension finie {\it mod\'er\'ee} de ${{F}} \,$ - c'est-\`a-dire  r\'esiduellement s\'eparable et d'indice de ramification premier \`a la carac\-t\'eristique r\'esiduelle $p\,$ -, la d\'erivation \'etendue \'etant encore de norme $1$\footnote{$^{(5)}$}{\sm on se ram\`ene au cas totalement ramifi\'e; l'anneau de valuation de ${{F}}'$ est alors engendr\'e, sur celui de ${{F}}$, par une uniformisante $\pi$, qui v\'erifie une \'equation d'Eisenstein $P(\pi)=0$. En utilisant le fait que $v'(P'(\pi))= e-1$  [$\bf{S}$, III,7, prop. 13], il est facile de d\'efinir $\partial(\pi)$ en appliquant  $\partial$ \`a l'\'equation d'Eisenstein, et $v'(\partial(\pi))=v'(\pi)=1$.}.
 
\medskip \noindent -  on suppose que la dimension $\mu$ du module diff\'erentiel ${{M}}$ est strictement inf\'erieure \`a la carac\-t\'eristique r\'esiduelle si cette derni\`ere est non nulle:
 $$\mu < p\;\;\;\;{\rm si}\;\; \;\;p\neq 0.$$  

 \proclaim Th\'eor\`eme 2.3.1. {\rm (D\'ecom\-po\-si\-tion de Turrittin-Levelt).} 
\smallskip 1) Il existe une extension    ${{F}}'/{{F}}$ finie mod\'er\'ee et une d\'e\-com\-po\-si\-tion du module diff\'erentiel ${{M}}_{{{F}}'} := {{M}}\otimes_{{F}}{{{F}}'}$ sur ${{F}}'$:
 $${{M}}_{{{F}}'}  =  \oplus_j\;M_{\bar\phi_j}\,, \leqno{(2.3.1.1)}$$ avec $$ M_{\bar\phi_j} \cong  \;  L_{\phi_j} \otimes  R_j\,, \leqno{(2.3.1.2)}$$  
o\`u 
\medskip\noindent - les $\phi_j$ d\'esignent des \'el\'ements de ${{F}}'$, dont les classes $\bar \phi_j$ modulo l'anneau de valuation $\sV'$ de ${{F}}'$ sont {\rm deux \`a deux distinctes},
   \medskip\noindent -  $ L_{\phi_j} = ({{F}}', {{\nabla}}_j(\partial))$ avec ${{\nabla}}_j(\partial)(1) = \phi_j $,
 \medskip\noindent -  $ R_j$ est un module diff\'erentiel r\'egu\-lier. 
\medskip 2) Cette d\'e\-com\-po\-si\-tion $(2.3.1.1)$ est unique.
\medskip 3)  On peut prendre ${{F}}'/{{F}}$ galoisienne de degr\'e divisant $\mu !$.
 \medskip 4) Si l'on remplace $\partial$ par $f.\partial$ avec $v(f)=0$, la d\'e\-com\-po\-si\-tion ne change pas ($\bar\phi_j$ est seulement remplac\'e par $\overline{f. \phi}_j$).
  \par

\noindent {\bf Remarque 2.3.2.} La condition que les $\bar \phi_j$  soient  deux \`a deux distincts 
est raisonnable compte tenu de ce que $L_{\phi }$ est r\'egu\-lier si $\phi\in  \sV'$.
 
\medskip\noindent {\bf D\'emonstration} (esquisse)\footnote{$^{(6)}$}{\sm voir aussi [$\bf{T}$], [$\bf{L}$] en \'egale carac\-t\'e\-ris\-tique nulle 
.}. 1)  On raisonne par r\'ecurrence, pour l'ordre lexicographique), sur le couple $ (\mu \in {\bf N}, \;\rho_v \in {{1}\over{ \mu !  }}{\bf N}) $, les cas $\mu\leq 1$ et $\rho_v = 0$ \'etant triviaux. 
  Soit $\bf m$ une base cyclique (\cf remarque 2.1.5). Avec les notations de 1.1.3, on a 
$$\rho_v( {{\nabla}}(\partial))  =  \max (0,   \max_{j=0,\dots,\mu-1}  \, {{-v(\theta_j)}\over {\mu -j}}).$$ 
Quitte \`a faire une extension mod\'er\'ee de ${{F}}$, on peut supposer  (puisque $\mu<p$ si $p\neq 0$)  que ${{F}}$ contient un \'el\'ement  $\xi$ de valuation $\rho_v $ \footnote{$^{(7)}$}{\sm si le maximum des ${-v(\theta_j)}\over {\mu -j}$ est $\geq 0$ et atteint pour $j=j_0$, on peut prendre $\xi= \theta_{j_0}^{-1/(\mu-j_0)}$ puisque $\mu<p$ si $p\neq 0$.}.   On peut aussi supposer que ${{F}}$ contient les valeurs propres $\zeta_i$ de la matrice de ${{\nabla}}(\partial)$ dans la base ${\bf m}'$ de $M$ introduite en 1.1.2.  

\noindent Si les $\zeta_i$ ne sont pas dans la m\^eme classe modulo l'id\'eal maximal  de $\sV'$ et si $\rho_v >0$, on peut appliquer le point $1)$ du lemme de d\'e\-com\-po\-si\-tion (avec $ R= \sV'$,  ${\Bbb M}=$ le $\sV'$-module engendr\'e par ${\bf m}'$, et $\delta = \xi.\partial$), ce qui diminue $\mu$.

\noindent Sinon, et si $\rho_v>0$, la classe des $\zeta_i$ n'est pas l'id\'eal maximal de $\sV'$ d'apr\`es 1.1.2, 2), et en tensorisant par $ L_{-\zeta_1/\xi}$, on diminue $\rho_v$. Ceci \'etablit le point 1). 

Pour 2) et 3), on s'appuie sur le lemme suivant:

 \proclaim Lemme 2.3.3. Soit ${{F}}''$ une autre extension de ${{F}}$ et soit $\iota: {{F}}'\to {{F}}''$ un 
${{F}}$-isomorphisme. Soient $\phi'\in {{F}}',\, \phi''\in {{F}}''$ tels que $ \iota(\phi')-\phi''$ ne soit pas dans l'anneau de valuation. Soient $R'$ (\resp $R''$) un module diff\'erentiel r\'egu\-lier sur ${{F}}'$ (\resp ${{F}}''$). Alors 
$$Hom_{{{\nabla}}(\partial), \iota}( L_{\phi'}\otimes R',\,  L_{\phi''}\otimes R'')=0.$$
\par
\noindent En effet, cet espace d'homomorphismes horizontaux s'identifie \`a  $$Hom_{{{\nabla}}(\partial), {{F}}''}( L_{\iota(\phi')-\phi''} ,\, \iota(R')^\vee\otimes R''),$$ qui est nul car  $ L_{\iota(\phi')-\phi''}$ est irr\'egu\-lier et irr\'eductible, tandis que $\iota(R')^\vee\otimes R''$ est r\'egu\-lier. \sq

\medskip \noindent Le point 2) s'en d\'eduit imm\'ediatement. Pour en d\'eduire 3), prenons ${{F}}''={{F}}'$ galoisien sur ${{F}}$, de groupe $G$, et prenons pour $\iota$ un \'el\'ement de $G$. Alors on voit que $\iota $ permute les $M_{\bar\phi_j}$ en permutant les $ {\bar\phi_j}$. Le sous-groupe normal $G'$ de $G$ qui agit trivialement sur les $ {\bar\phi_j}$ fixe donc les projecteurs qui d\'efinissent la d\'e\-com\-po\-si\-tion ${{M}}_{{{F}}'} =  \oplus\;M_{\bar\phi_j}\,, $ qui descend donc sur une extension galoisienne de degr\'e $$[G:G']\mid \mu!.$$  
 Le point  4) est imm\'ediat. \sq
      
    \medskip\noindent {\bf Exemple 2.3.4.}  Prenons pour ${{F}}$ un corps de s\'eries for\-melles $K((x))$ avec la valuation $x$-adique, et $\partial = x{d\over dx}$. Rappelons que les extensions mod\'er\'ees de $K((x))$ sont de la forme ${{F}}' = K'((x^{1/e}))$, o\`u $K'$ est une extension s\'eparable de $K$ et o\`u $e$ est premier \`a $p$ si $p\neq 0$.

La projection ${{F}}'\to {{F}}'/ \sV'$ admet une section $K$-lin\'eaire canonique d'image $ x^{-1/e}K'[ x^{-1/e}]$. On peut donc choisir les $  \phi_j$ dans $ x^{-1/e}K'[ x^{-1/e}]$ pour un entier $e\geq 1$ convenable, qui est un multiple de d\'enominateur $d$ de $\rho_v$ (mais on prendra garde qu'il ne lui est pas \'egal en g\'en\'eral). On peut prendre  $\xi= x^{1/d}$ dans 2.1.2. 

 On v\'erifie ais\'ement que les $\phi_j$ engendrent une extension galoisienne ${{F}}'$ de ${{F}}$ qui est {\it l'extension minimale} sur laquelle a lieu la d\'e\-com\-po\-si\-tion ${{M}}_{{{F}}'} =  \oplus\;M_{\bar\phi_j}\,.$

Si $\rho_v >0$, les  $\zeta_i $ apparaissant dans la preuve de 2.3.1 sont alors exactement les termes de degr\'e $x^{-1/d}$ des $  \phi_j $, multipli\'es par $ x^{1/d}$. 

  En carac\-t\'eristique nulle, il est bien connu que tout module diff\'erentiel sur $F$ r\'egu\-lier est extension it\'er\'ee de modules diff\'erentiels de dimension un, et il en est donc de m\^eme de ${{M}}_{{{F}}'}$.
        
  En carac\-t\'eristique non nulle, il n'en est rien en g\'en\'eral (cela d\'epend de la $p$-courbure, \cf [$\bf{DGS}$, III.2], [${\bf A}$, 3.2.2]).

\medskip\noindent {\bf 2.4. Polyg\^one de Newton.}  Un argument \'el\'ementaire de descente galoisienne montre que si l'on regroupe les $M_{\bar\phi_j}$ suivant la valeur de $\,\max (0, -v({\bar\phi_j}))\in {1\over \mu !}\N,$\footnote{$^{(8)}$}{\sm il s'agit ici de l'extension de la valuation $v$ \`a ${{F}}'$, non normalis\'ee.} la d\'e\-com\-po\-si\-tion (moins fine) cor\-res\-pondante descend \`a ${{F}}$ et fournit la {\it d\'e\-com\-po\-si\-tion suivant les pentes}:
 $$\displaystyle {{M}} =  \bigoplus_\sigma\;M_{(\sigma)},\leqno{(2.4.1.1)}$$ 
 et il est facile de voir que 
  $ \rho_v(M_{(\sigma)})=\sigma. $  Cette d\'e\-com\-po\-si\-tion ne change pas si l'on remplace $\partial$ par $f.\partial$ avec $v(f)=0$.
    
 \proclaim D\'efinition 2.4.1. Le {\rm polyg\^one de Newton} de $({{M}}, {{\nabla}}(\partial))$ est le polyg\^one convexe contenu dans $[0,\mu]\times \R$, non born\'e sup\'erieurement, dont le sommet le plus \`a droite est $(\mu, 0)$, et dont le segment de pente $\sigma$ est de longueur \'egale \`a la dimension de $M_{(\sigma)}$ (pour tout $\sigma\in \R_+$).  
 \par
\noindent  On le note $ NP_v( {{M}} ,{{\nabla}}(\partial))\,$ (ou $NP_v( {{\nabla}}(\partial))$ ou encore $NP_v( {{M}})$ pour abr\'eger).
Il ne change pas si l'on remplace $\partial$ par $f.\partial$ avec $v(f)=0$.
Sa plus grande pente est $\rho_v( {{\nabla}}(\partial)) $.
 
\medskip De ce que $ \rho_v(M_{(\sigma)})=\sigma $, et de 1.1.3 $i)$, on d\'eduit:
 
  \proclaim Corollaire 2.4.2. Choisissons un vecteur cyclique $m$, et \'ecrivons  $$\displaystyle   \sum_{ 0}^{ \mu}\, \theta_i  {{\nabla}}(\partial)^{i}(m)=0,$$ avec $\theta_\mu= -1$.   Alors $NP_v( {{\nabla}}(\partial))$ est l'enveloppe convexe des demi-droites  $$\{x_1=i, \, x_2\geq v(\theta_i) \}, \;\; i=0,\ldots, \mu.$$ En par\-ti\-cu\-lier, ses sommets sont \`a co\-or\-don\-n\'ees enti\`eres.
 \par

\proclaim Definition 2.4.3.  {\rm L'irr\'egularit\'e} de $M$, not\'ee ${\rm ir}_v(M)$, est la hauteur du polyg\^one de Newton, c'est-\`a-dire (avec les normalisations que nous avons choisies) l'oppos\'e de l'ordonn\'ee du sommet le plus bas. \par
C'est donc un entier naturel qui v\'erifie l'in\'egalit\'e 
$${\rm ir}_v(M)\leq \mu.\rho_v(M). \leqno{(2.4.3.1)}$$

\medskip\noindent {\bf 2.5. D\'erivations commutantes.}  Soit $\partial'$ une autre d\'erivation de ${{F}}$ de norme $1$. Supposons ${{M}} $ muni d'une action $K$-lin\'eaire ${{\nabla}}(\partial')$ de $\partial'$ v\'erifiant la r\`egle de Leibniz et commutant \`a ${{\nabla}}(\partial)$.

\proclaim Proposition 2.5.1.  $1)$ La d\'e\-com\-po\-si\-tion de Turrittin-Levelt est stable sous ${{\nabla}}(\partial')$. 
\medskip\noindent $2)$ Supposons que pour tout $f\in {{F}}$, $\vert \partial(f)\vert =\vert f\vert$ si $\vert f\vert \neq 1$. Alors $$ \rho_v({{\nabla}}(\partial'))\,\leq \,\rho_v({{\nabla}}(\partial ))\,,\;\;\; NP_v({{\nabla}}(\partial'))\,\subset \,NP_v({{\nabla}}(\partial )). \leqno{(2.5.1.1)}$$\par

\dem Quitte \`a passer \`a une extension mod\'er\'ee, on peut supposer la d\'e\-com\-po\-si\-tion de Turrittin-Levelt d\'efinie sur ${{F}}$. 
En suivant la d\'emonstration de cette d\'e\-com\-po\-si\-tion, il est clair que $1)$ d\'ecoule du point $3)$ du lemme de d\'e\-com\-po\-si\-tion 2.2.1 (en prenant $g=\xi$ et $g'$ convenable pour que $\partial'$ respecte le r\'eseau engendr\'e par ${\bf m}'$ sur l'anneau de valuation de ${{F}}$). 

Pour le point $2)$, on peut supposer que ${{M}} $, muni de ${{\nabla}}(\partial)$ et ${{\nabla}}(\partial')$, est ind\'ecomposable. D'apr\`es le point $1)$, la d\'e\-com\-po\-si\-tion de Turrittin-Levelt tant de  $({{M}}, {{\nabla}}(\partial))$ que de $({{M}}, {{\nabla}}(\partial'))$ n'a qu'un seul facteur, not\'e $ L_\phi \otimes  R$ et $ L'_{\phi'} \otimes  R'$ respectivement. Il est donc isocline et l'assertion sur $NP$ d\'ecoule de celle sur $\rho$. On a alors $\bigwedge^\mu {{M}} \cong  L_\phi^{\otimes \mu}\otimes \bigwedge^\mu  R = L_{\mu\phi} \otimes \bigwedge^\mu  R$ (et de m\^eme pour $\partial'$). 
Compte tenu de ce que $\mu < p$ si $p\neq 0$, il suit que les in\'egalit\'es a priori
$$\;   \rho_v(\bigwedge^\mu {{M}}, \bigwedge^\mu {{\nabla}}(\partial))\leq \rho_v({{M}}, {{\nabla}}(\partial)), \;\;\;    \rho_v(\bigwedge^\mu {{M}}, \bigwedge^\mu {{\nabla}}(\partial'))\leq \rho_v({{M}}, {{\nabla}}(\partial')) $$ sont des \'egalit\'es,
ce qui nous ram\`ene au cas $\mu =1$. Dans ce cas, prenant une base $m$ et \'ecrivant 
$$ {{\nabla}}(\partial)m=\theta m,\; {{\nabla}}(\partial')m= \theta' m , $$ on a $$\;\;  \partial'(\theta)=  \partial(\theta'),  $$ $$ \rho_v({{\nabla}}(\partial')) = \max (0, - v(\theta')) \,\leq  \max(0, - v(\partial (\theta'))=  \max(0, - v(\partial' (\theta)),$$ d'apr\`es l'hypoth\`ese sur $\partial$, d'o\`u 
 $$ \rho_v({{\nabla}}(\partial'))\leq \max(0, - v( \theta))= \rho_v({{\nabla}}(\partial )) . \;\;\;\;\;\;\square$$

 \proclaim Corollaire 2.5.2.  Sous l'hypoth\`ese de $2)$, si ${{M}} $ est r\'egu\-li\`ere vis-\`a-vis de $\partial$, alors ${{M}} $ est aussi r\'egu\-li\`ere vis-\`a-vis de $\partial'$.
 \par
 
    \smallskip \noindent{\bf Remarque 2.5.3.} L'argument de $2)$ montre qu'on peut toujours prendre $\phi= {\theta/\mu}$. $L_\phi$ et $R$ sont alors munis d'actions de $\partial'$ compatibles avec $\nabla(\partial')$ sur le produit tensoriel. En outre, dans la si\-tua\-tion $F= K((x))$ de 2.3.4, on observe que puisque $\partial'$ est de norme $1$, elle commute \`a la troncation $K'((x^{1/e}))\to x^{-1/e} K'[x^{-1/e}]$. On d\'eduit de l\`a qu'on peut associer au choix canonique $\phi_j\in   x^{-1/e} K'[x^{-1/e}]$ un unique $\phi'_j\in  x^{-1/e} K'[x^{-1/e}]$ tel que 
$ \;\;   \partial'(\phi_j)=  \partial(\phi_j') .$

\bigskip \bigskip 

\bigskip {\bf 3. D\'e\-com\-po\-si\-tion de Turrittin-Levelt avec param\`etres, points tournants.} 

\medskip Dans ce paragraphe, nous analysons d'un point de vue purement alg\'ebrique le ph\'enom\`ene classique des points tournants, c'est-\`a-dire (essentiellement) des points du diviseur polaire o\`u la d\'ecomposition de Turrittin-Levelt ne se sp\'ecialise pas. Outre le th\'eor\`eme 3.4.1, ce paragraphe contient  deux r\'esultats importants pour la suite (\S 6): 3.1.1 et 3.3.1.

\medskip\noindent {\bf 3.1. Rang de Poincar\'e-Katz et sp\'ecialisation.} Pla\c cons-nous dans la si\-tua\-tion o\`u $F=K((x))$ muni de la valuation $x$-adique $v$, et o\`u $\partial = x{{d}\over { dx}}$. 
 Supposons que $K$ soit donn\'e comme corps de fractions d'un anneau $A$ noeth\'erien int\'egralement clos. 

Soit $M$ un module diff\'erentiel sur $$A((x))= A[[x]][{{1}\over{x}}],$$ c'est-\`a-dire un $A((x))$-module projectif de type fini muni 
 d'une action $A$-lin\'eaire continue $ {{\nabla}}(\partial)$ de $\partial$ v\'erifiant la r\`egle de Leibniz.
 
   On note $\rho= \rho(M))$  le rang de Poincar\'e-Katz du module diff\'erentiel $M_F$ (relativement \`a la  valuation $x$-adique). La simplicit\'e de la preuve du r\'esultat suivant illustre l'avantage du point de vue spectral sur le rang de Poincar\'e-Katz.
 
\proclaim Lemme 3.1.1. Soit $P$ un point de $\Spec A$ et soit $(M_{(P)}= M\otimes_{A((x))}\,\kappa_P((x)), \nabla_{(P)})$ le module diff\'erentiel sp\'ecialis\'e sur $\kappa_P((x))$. On a 
$$\rho(M_{(P)})\leq \rho(M).\leqno{(3.1.1.1)}$$
\par 
\dem Munissons $M$ d'une norme $x$-adique quelconque, et le $\kappa_P((x))$-espace vectoriel de la norme $x$-adique quotient. Par d\'efinition 
$$ \rho ( {M}) = \max(0, -v_{M, sp }(\nabla (\partial ))),\;\;\rho ( {M}_{(P)})= \max(0, -v_{M_{(P)}, sp }(\nabla_{(P)}(\partial ))),$$
et il est clair que $v_{M_{(P)}, sp }(\nabla_{(P)}(\partial ))\geq v_{M , sp }(\nabla (\partial ))$    (le choix des normes $x $-adiques sur $M$ et $M_{(P)}$ n'a en fait aucune importance puisqu'on ne s'int\'eresse qu'aux valuations spectrales).
 \sq 
\proclaim Corollaire 3.1.2. Si $M$ est r\'egulier, il en est de m\^eme de $M_{(P)}$.\sq
\par 
 
  \bigskip\noindent {\bf 3.2. Points semi-stables.} 
  On suppose d\'esormais que le rang $\mu$ de $ M$ sur $A((x))$ est strictement inf\'erieur aux carac\-t\'e\-ris\-tiques r\'esiduelles de $A$ si celles-ci sont non nulles.
  
   On dispose de la d\'ecom\-po\-si\-tion de Turrittin-Levelt (\cf 2.3.1, 2.3.4) de $M_{F'}= M\otimes_{A((x))} F'$ sur une extension galoisienne $F'= K'((x^{1/e}))$ convenable de $F$ de degr\'e divisant $\mu !$  
  $${{M}}_{{{F}}'}  =  \oplus_j \; M_{ \bar\phi_j, {{F}}'   }\,, \;\;  M_{ \bar\phi_j, {{F}}'   } \cong  \;  L_{ \phi_j, {{F}}'   } \otimes  R_{ j, {{F}}'   }\,  \leqno (3.2.1.1) $$  
  avec $\phi_j\in x^{-1/e}K'[x^{-1/e}]$, de degr\'e $\geq -\rho  $ en $x$, l'un des $\phi_j$ au moins \'etant exactement de degr\'e $-\rho $.

\bigskip Ce paragraphe est d\'evolu \`a l'\'etude de la question suivante.

  \smallskip\noindent  {\bf Question 3.2.1.}  {\it Les coefficients des $\phi_j$ sont-ils dans la cl\^oture int\'egrale $A'$ de $A$ dans $K'$, et cette d\'ecom\-po\-si\-tion de $M_{K'((x^{1/e}))}$ descend-elle en une d\'ecom\-po\-si\-tion 
sur $A'((x^{1/e}))$ }
   $${{M}}_{{{A}}'}  =  \oplus_j \; M_{  \bar\phi_j}\,, \;\;  M_{ \bar\phi_j} \cong  \;  L_{  \phi_j} \otimes  R_{ j}\;? \leqno (3.2.1.2) $$   
   {\bf Remarque 3.2.2.} Soit $K''$ une $K'$-alg\`ebre et $A''\subset K''$ un sous-anneau tel que $A''\cap K'=A'$. Alors la question 3.2.1 a une r\'eponse positive si et seulement s'il en est ainsi de la question analogue avec $A',K'$ remplac\'es par $A'',K''$. 
    
\noindent   En effet, puisque $M$ est projectif de type fini, il en est de m\^eme de $\,{{\sE}nd}\, M = M^\vee\otimes M$, donc il existe des \'el\'ements $n_1,\ldots, n_s$ de  $\,{{\sE}nd}\, M $ et des \'el\'ements $n^\vee_1,\ldots, n^\vee_s$ du dual tels que pour tout $n\in \,{{\sE}nd}\, M$, $n= \sum \langle n_i^\vee,n\rangle n_i$. Dire qu'une d\'ecom\-po\-si\-tion de $M_{K'((x^{1/e}))}$  descend sur $A'((x^{1/e}))$ revient \`a dire que pour les projecteurs $n$ qui la d\'efinissent, les \'el\'ements $\langle n_i^\vee,n\rangle\in K'((x^{1/e}))$ sont en fait dans $A'((x^{1/e}))$. Mais ceci \'equivaut aussi \`a $\langle n_i^\vee,n\rangle\in A''((x^{1/e}))$.
 
 \noindent  Cette remarque montre plus g\'en\'eralement que la r\'eponse \`a 3.2.1 ne d\'epend pas de $F'$.
 
 \bigskip La r\'eponse \`a 3.2.1 est positive en rang $\mu=1$. En revanche, sans hypoth\`ese suppl\'ementaire, la question a une r\'eponse n\'egative:
 
    \smallskip\noindent {\bf Contre-exemple 3.2.3.}  Supposons que $A= k[[x_2]]$, et que $M$ soit libre de rang deux, et que $\nabla(\partial)$ soit donn\'e dans une base $ (m_1,m_2)$ par
    $$\nabla(\partial)(m_1)= {x_2\over x}m_1,\; \nabla(\partial)(m_2)= - m_1 .$$
   Les pentes sont $0$ et $1$, et la d\'ecom\-po\-si\-tion de Turrittin-Levelt sur $F= k((x_2))((x))$ s'\'ecrit 
    $$M_F \cong  L_{{{x_2}\over {x}}, F}  \oplus  L_{ {0}, F}  .$$ 
    Mais en $x_2=0$, $M$ se sp\'ecialise en un module diff\'erentiel r\'egu\-lier ind\'ecomposable. Donc la d\'ecom\-po\-si\-tion de Turrittin-Levelt ne descend pas de $F $ \`a $k[[x_2]]((x))$.

 \proclaim D\'efinition 3.2.4. Soit $P$ un point de $Z=\Spec A$. Nous dirons que $P$ est {\rm semi-stable} pour  $M$ si la r\'eponse \`a 3.2.1 est positive lorsque l'on remplace $A'$ par le semi-localis\'e $A'_P$.\par
 
 En vertu de la remarque 3.2.2, la r\'eponse \`a 3.2.1 est positive si et seulement si tous les points sont semi-stables. Le point g\'en\'erique est \'evidemment semi-stable.
 
 Cette terminologie est motiv\'ee d'une part par l'analogie avec le probl\`eme de la r\'eduction semi-stable des $F$-isocristaux surconvergents dans le cas $p$-adique (\cf [$\bf Ke$]), d'autre part parce que ``stable" semble un antonyme acceptable de l'adjectif ``tournant" classique en analyse asymptotique, \cf infra 3.4.2.  
    
      \medskip\noindent  {\bf Remarques 3.2.5.} $1)$ Si tous les points de $Z$ sont semi-stables, la d\'ecomposition suivant les pentes (2.4.1.1) de $M_F$ provient d'une d\'ecomposition de $M$ lui-m\^eme.
       \smallskip \noindent $2)$ Si $P$ est semi-stable pour $M$, il l'est aussi pour  ${{\sE}nd}\,M = M^\vee\otimes M$.
       
   \proclaim Lemme 3.2.6. Si $P$ est un point semi-stable pour $M$, on a 
$$NP(M_{(P)})\subset NP(M),\;\; NP({{\sE}nd}\,M_{(P)})\subset NP({{\sE}nd}\,M).\leqno{(3.2.6.1)}$$
    \par 
  \noindent   C'est imm\'ediat \`a partir de la d\'efinition 2.4.1 (et compte tenu de la remarque pr\'ec\'edente). \sq  
    
  \smallskip  Nous montrerons dans l'appendice que (3.2.6.1) vaut m\^eme sans la condition de semi-stabilit\'e.

\bigskip \noindent {\bf 3.3. Forme normale.}  Dans la direction de 3.2.1, on a le r\'esultat g\'en\'eral suivant:

\proclaim Proposition 3.3.1. Pour tout $j$, le coefficient $ \phi_{j,-\rho}$ de degr\'e $-\rho$ en $x$ dans $\phi_j$ est entier, \ie appartient \`a $A'$. 
 \smallskip\noindent Plus g\'en\'eralement, \'etant donn\'e $r>0$, supposons que pour tout $s>r$,  le coefficient $\phi_{j, -s}$ de $x^{-s}$ dans $\phi_j$ soit ind\'ependant de $j$; alors pour tout $s\geq r$, $\phi_{j, -s}$ appartient \`a $A'$.\par 
 
La seconde assertion se d\'eduit de la premi\`ere en tensorisant par $L_{-\phi_{j,-\rho}x^{-\rho}}$ et en it\'erant.

 Compte tenu de ce qu'un anneau noeth\'erien int\'egrale\-ment clos est intersection de ses localis\'es en les id\'eaux premiers de hauteur $1$, on peut d'embl\'ee remplacer $A$ par un tel localis\'e puis par son compl\'et\'e, qui est un anneau de valuation discr\`ete complet, et il s'agit de montrer que les $ \phi_{j,-\rho}$ sont entiers sur $A$.    
 La proposition 3.3.1 est alors cons\'equence imm\'ediate du r\'esultat plus pr\'ecis suivant, qui fournit une sorte de forme normale pour $\nabla(\partial )$.
 
 \proclaim Th\'eor\`eme 3.3.2. Soit $A$ un anneau de valuation discr\`ete complet, de corps de fractions $K$, de corps r\'esiduel $k$. Soit $M$ un module diff\'erentiel de rang $\mu$ sur $ A((x))$ (pour la d\'erivation $\partial= x{d\over {dx}}$). Si $ \car k\neq 0$, on suppose que $\mu <  \car k$. 
 Soit $d$ le d\'enominateur  du rang de Poincar\'e-Katz $\,\rho\,$ de $M$.
 \smallskip Alors il existe une base $\bf n$ de $M_{A((x^{1/d}))}$  dans laquelle la matrice de  $\nabla( \partial )$ est de la forme $$x^{-\rho }G(x),\;\; {\rm avec}\;\,G(x)\in M_\mu(A[[x^{1/d}]]).$$  En outre,  si $\rho >0$, l'ensemble des valeurs propres non nulles de $\,G( 0 )$ co\"{\i}ncide avec l'ensemble des $\phi_{j,-\rho}$.\par
     
\dem D'apr\`es la th\'eorie des polyg\^ones de Newton de s\'eries de Laurent (\cf [${\bf DGS}$, II.3]), l'anneau $A((x))$ est  principal, de sorte que $M$ est automatiquement libre sur $A((x))$. Fixons-en une base ${\bf n}'$. 

Le corps de fractions de $A((x))$ est un sous-corps de $F=K((x))$ que nous notons provisoirement $F_0$. Par ailleurs, puisque $\car k=0$ ou $>\mu$, il existe un vecteur cyclique pour le module diff\'erentiel $M_{k((x))}$ (obtenu par r\'eduction modulo l'id\'eal maximal $\frak m$ de $A$). Relevons-le en un \'el\'ement $m$ de $M$. Alors $m$ est un vecteur cyclique pour le module diff\'erentiel $M_{F_0}$. En effet, il suffit de voir qu'il l'est sur une extension diff\'erentielle convenable de $F_0$, par exemple le corps de fractions du compl\'et\'e $\frak m$-adique $\widehat {A((x))}_\frak m$ de $A((x))$ (ce compl\'et\'e $\frak m$-adique est lui-m\^eme un anneau de valuation discr\`ete complet pour la topologie $\frak m$-adique, de corps r\'esiduel $k((x))$). Or, par le lemme de Nakayama, le plus petit sous-module de $ M_{\widehat {A((x))}_\frak m}$ stable sous $\partial $ et contenant $m$ est $ M_{\widehat {A((x))}_\frak m}$ lui-m\^eme puisqu'il en est ainsi modulo $\frak m$. Posons 
$${\bf m}\,=\,(m, {{\nabla}}(\partial)(m), \ldots, {{\nabla}}(\partial)^{\mu -1}(m)).$$ 
Pour simplifier les notations, rempla\c cons d'embl\'ee $x$ par $x^{1/d}$ et $\rho$ par l'entier $d\rho$. 

Alors la matrice de $ \nabla(\partial)$ dans la base 
  $${{\bf m}'} =  {\bf m}\pmatrix{1&0&&0\cr 0&x^\rho &&0\cr &&\ldots&\cr 0&0&& x^{\rho(\mu-1)  }} $$ de $M_{F_0}$ est de la forme 
  $$x^{-\rho} H(x),\;\; {\rm avec}\;\,H(x)\in M_\mu(K[[x]]\cap F_0).$$ 
 Si $\rho(M)>0$, l'ensemble des valeurs propres non nulles de $H(0)$ co\"{\i}ncide avec l'ensemble des  coefficients $\phi_{j,-\rho}$ des termes de degr\'e $-\rho$ en $x $ dans les $\phi_j$, \cf 2.3.4 (ceci prouve d\'ej\`a 3.3.1). 
 
  Pour chasser les d\'enominateurs de $H$, \ie se ramener \`a une matrice de s\'eries \`a coefficients dans une extension finie enti\`ere de $A$, nous allons employer une technique \'eprouv\'ee dans la th\'eorie des \'equations diff\'erentielles ultram\'etriques ({\it cf. e.g.} [${\bf DGS}$, V.5.1]). 

Notons $J\in  M_\mu(A((x)))$ la matrice de passage de ${\bf n}'$ vers $\bf m'$: 
$${\bf m}'= {\bf n}' J.$$ Le probl\`eme est que son d\'eterminant n'est pas n\'ecessairement inversible dans $A((x))$. D'apr\`es la th\'eorie des polyg\^ones de Newton de s\'eries de Laurent ([${\bf DGS}$, II.3]),
$\det J$ est produit d'un polyn\^ome $\varphi\in A[x]$ unitaire et d'un \'el\'ement inversible de $A((x))$. Soit $A''$ l'extension enti\`ere obtenue par adjonction des racines de $\varphi $ et $K''$ son corps de fractions.
  Alors $J$ se factorise en $$J= J'(J'')^{-1},\;\;{\rm avec}\;\;J'\in GL_\mu(A''((x))),\;\;  J''   \in GL_{\mu }(K'' [x]_0) $$  (l'indice $0$ d\'esignant la localisation en $0$).
  
  En effet, soit $\xi$ l'un des z\'eros de $\varphi $. Pour \'etablir cette factorisation, on se ram\`ene par r\'ecurrence (sur le nombre et la multiplicit\'e de ces z\'eros)  \`a trouver une matrice $J''_1\in   GL_{\mu }(K''[x ]_0)$  telle que $J J''_1$ n'ait pas de p\^ole en $\xi$ et que ${\rm ord}_{\xi}\, {\rm det}\,  J J''_1   < {\rm ord}_{\xi}\, {\rm det}\, J $.
Soit $\lambda_{1 },\dots, \lambda_\mu\in A$ les coefficients d'une relation de d\'ependance lin\'eaire non triviale entre les {\it colonnes} of $J_{\mid x =\xi}$. On peut supposer que l'un d'entre eux, soit $\lambda_i$, vaut $1$, et il est facile de voir que 
 $$J''_1= \pmatrix{I_{i-1 } &&  \lambda_1/(x -\xi) &&0_{i-1,\mu-i } \cr 
&&\vdots && \cr
  &  &
 \lambda_{ i -1}/(x -\xi)&  & \cr 
0&&1/(x-\xi) &&0 \cr
&&  \lambda_{i+1}/(x -\xi)&&\cr 
&&\vdots && \cr
0_{\mu- i, i-1}& & \lambda_\mu/(x -\xi) &&I_{\mu-i } }$$  
$${\rm (d'inverse}\;\; (J''_1)^{-1}=\pmatrix{I_{i-1 } && -\lambda_1 &&0_{i-1,\mu-i } \cr 
&&\vdots && \cr
  &  &
-\lambda_{ i -1}&  & \cr 
0&&x -\xi &&0 \cr
&& -\lambda_{i+1}&&\cr 
&&\vdots && \cr
0_{\mu- i, i-1}& &-\lambda_\mu &&I_{\mu-i } })\;\;\rm{remplit}\; {\rm cet}\;{\rm office}.$$
 
  Gr\^ace \`a la factorisation $J= J'(J'')^{-1}$, on voit que $${\bf n}'' :={\bf n}'J' = {\bf m}'J''$$ engendre un $A''[[x]]$-r\'eseau ${\Bbb M}''$ de $ M_{A''((x))}$ stable sous $ x^\rho\partial = x^{\rho+1} {\partial\over \partial x }$. Alors ${\Bbb M}= {\Bbb M}''\cap M$ est un $A[[x]]$-r\'eseau de $M$, stable sous $ x^\rho\partial$, et toute base $\bf n$ de ${\Bbb M}$ convient pour la premi\`ere assertion de 3.3.2.  
  
La matrice de $x^\rho\partial $ dans la base ${\bf n}''$ est $ (J'')^{-1} H(x )  J'' +  ( J'')^{-1}  x^\rho\partial   J''$. Soit $G(x)$ la matrice de $x^\rho\partial $ dans $\bf n$, vue comme base de ${\Bbb M}''$. Comme la classe de conjugaison de $G(0)$ ne d\'epend pas du choix de la base, c'est aussi celle de $H(0)$.    En par\-ti\-cu\-lier, les valeurs propres de $G(0)$ sont les m\^emes que celles de $H(0)$, et celles non nulles sont donc les $\phi_{j,-\rho}$ si  $\rho(M)>0$.  $\square$

    \bigskip \noindent {\bf 3.4. Points stables, points tournants.}  Revenons \`a la si\-tua\-tion de 3.2.   On note  $NP(M)$  le  polyg\^one de Newton  du module diff\'erentiel $M_F$ (relativement \`a la  valuation $x$-adique).

 \proclaim Th\'eor\`eme 3.4.1. Soient $P$ un point de $Z=\Spec A$, $A_P$ le localis\'e de $A$ en $P$,   et ${\kappa_P}$ son corps r\'esiduel. 
  \medskip\noindent  1) Les cinq conditions suivantes sont \'equivalentes:
 \medskip $i)$ quitte \`a remplacer $\Spec A$ par un voisinage de Zariski de $P$, les $\phi_j$ sont dans $ x^{-1/e}A'[x^{-1/e}]$ et v\'erifient la condition 
 \medskip\noindent   {\rm $ (\ast)   \;\;\;\;  \phi_j  $ (\resp $\phi_i -\phi_j$) est inversible dans $A'((x^{1/e}))$ si non nul (\resp si $i\neq j$),}
\medskip\noindent et la d\'ecom\-po\-si\-tion de Turrittin-Levelt de $M_{K'((x^{1/e}))}$ descend sur $A'((x^{1/e}))$, 
 \medskip $ii)$ idem en rempla\c cant $A'$ par le semi-localis\'e $A'_P$,
 \medskip $iii)$ les coefficients des termes de plus bas degr\'e des $\phi_j$ (\resp $\phi_i- \phi_j$)  sont des unit\'es dans $A'_P$,
\medskip $iv)$ on a \'egalit\'e de polyg\^ones de Newton
 \medskip\noindent  {\rm $(\ast\ast) \;\;\;\;\; NP(M_{(P)})= NP(M),\;\;\; NP({{\sE}nd}\,M_{(P)})= NP({{\sE}nd}\,M)  $},
 \medskip $v)$ pour tout germe (formel) de courbe $\frak C$ coupant $x=0$ transversalement en $P$,   on a \'egalit\'e de polyg\^ones de Newton
 \medskip\noindent  {\rm ${\;}\;\;\;\;\;\;\;\;\; NP(M_{\frak C})= NP(M),\;\;\; NP({{\sE}nd}\,M_{\frak C})= NP({{\sE}nd}\,M) .$}
 \medskip\noindent 2) Si $A$ est un anneau r\'egu\-lier, ces conditions 
 entra\^{\i}nent qu'on peut choisir $A'$ {\rm \'etale} sur $A$ au-dessus d'un voisinage de $P$.
 \par
 
 \proclaim D\'efinition 3.4.2. Nous dirons que $P$ de $Z$ est un {\rm point stable} pour $M$ s'il v\'erifie  les conditions \'equivalentes de 1) ci-dessus. Un point qui n'est pas stable est dit {\rm tournant}. \par 
 
Il est clair qu'un point stable est semi-stable au sens de 3.2.4. La terminologie ``point tournant" est conforme \`a la tradition en analyse asymptotique (\cf [${\bf W}$]). Le th\'eor\`eme 3.4.1 est d'ailleurs proche de la caract\'erisation des points tournants donn\'ee dans [$\bf{Sc}$], et raffine des r\'esultats de [$\bf{BaV}$, 5.7].  
   
   \proclaim Corollaire 3.4.3. Les points tournants forment un ferm\'e, vide ou purement de codimension un, de $\Spec A$. \par 
   
\dem   Cela d\'ecoule du crit\`ere $ii)$ pour les points stables et du Haupt\-idealsatz de Krull.\sq

 \medskip\noindent {\bf Exemple 3.4.4.} Pour $Z={\bf A}^1=\Spec k[x_2]$ et $\phi= x_2/x_1$, l'unique point tournant pour $L_\phi$ est $0$, qui est semi-stable.

   \bigskip\noindent {\bf 3.5. Preuve de 3.4.1.} 1) Les implications $i)\Rightarrow ii)\Rightarrow iii)$ et $v)\Rightarrow iv)$ sont imm\'ediates, et $iii)\Rightarrow v)$ est ais\'ee: noter que si $P'$ est l'un quelconque des points de $\Spec A'$ au-dessus de $P$, les degr\'es en $1/x$ de  $$\phi_j(P'),\;\resp\, \phi_j,\;\resp\,\phi_i(P')-\phi_j(P'),\;\resp\,\phi_i-\phi_j\; (i\neq j)$$   apparaissent  comme pentes de $$NP(M_{(P)}), \;\resp\,NP(M),\;\resp\,NP({{\sE}nd}\,M_{(P)}),\;\resp\,NP({{\sE}nd}\,M),$$  compte tenu de 3.1.2.    
  
   Prouvons $iv)\Rightarrow ii)$. Pour cela, on peut remplacer $A$ par $A_P$, puis par son compl\'et\'e (compte tenu de la remarque 3.2.2); supposons donc que $A$ soit un anneau de valuation discr\`ete complet, et $P$ le point ferm\'e de $Z=\Spec A$.  
  
  \noindent  On proc\`ede par r\'ecurrence sur le rang $\mu$, en remarquant que la condition $iv)$ est stable par facteur direct. Le cas $\mu=1$ est imm\'ediat. Par ailleurs, le cas $\rho = 0$ \'etant trivial, on peut supposer $\rho >0$. 
  
  Supposons d'abord qu'il n'y ait qu'un seul $\phi_j$. Il d\'ecoule alors de 3.3.1 que $\phi_j$ est \`a coefficients dans $A' \cap K= A$. En outre, la condition $iv)$ implique que $\phi_{j,-\rho}$ est une unit\'e de $A$, ce qui \'etablit $ii)$ dans ce cas.
  
  Supposons dor\'enavant qu'il y ait plusieurs $\phi_j$. Fixons $j$, et posons $$\rho'= \max_{i\neq j}\, - v_x(\phi_{i } -\phi_{j }) \in  \, ] 0, \rho  ] .$$ Il d\'ecoule de 3.3.1 que la composante $\phi^{\leq -\rho'}_j $ de $\phi_j$ de degr\'e $\leq -\rho'$ est \`a coefficients dans $A'$. En outre le module diff\'erentiel $$M' = M\otimes L_{-\phi^{\leq -\rho'}_j}$$ a au moins deux pentes, et son rang de Poincar\'e-Katz est $\rho'$. C'est donc aussi le rang de Poincar\'e-Katz de  ${{\sE}nd}\,M'$. Par ailleurs, la condition $iv)$ implique $NP({{\sE}nd}\,M'_{(P)})= NP({{\sE}nd}\,M')  $. On tire de l\`a que les coefficients $ \phi_{i,-\rho'}-\phi_{j,-\rho'}$ sont ou bien nuls, ou bien inversibles dans $A'$, et les deux cas interviennent.
  
  Le lemme de d\'e\-com\-po\-si\-tion 2.2.1 (point 1) s'applique alors \`a la d\'erivation $\delta =x^{\rho' +1} {\partial\over \partial x }  $  de $A'[[x^{1/e}]]$ (qui est un anneau complet d'id\'eal maximal $\frak m$ engendr\'e par l'id\'eal maximal de $A'$ et $x$),  et au $A'[[x^{1/e}]]$-r\'eseau ${\Bbb M'}$ engendr\'e par une base $\bf n$ comme dans 3.3.2. D'o\`u une d\'e\-com\-po\-si\-tion du module diff\'erentiel $ M'_{A'[[x^{1/e}]]}$, donc aussi de $ M_{A'[[x^{1/e}]]}$, ce qui permet de dimi\-nuer $\mu$ tout en respectant la condition $(\ast\ast)$ sur les polyg\^ones de Newton. On en d\'eduit, par r\'ecurrence, que la d\'ecom\-po\-si\-tion de Turrittin-Levelt de $M_{K'((x^{1/e}))}$ descend sur $A'((x^{1/e}))$. Il est alors clair que la condition $iv)$ implique $ii)$.

\smallskip  Pour terminer la preuve du point 1), il suffit d'\'etablir $iii)\Rightarrow i)$. La condition $iii)$ \'etant locale, elle implique (compte tenu de ce qui pr\'ec\`ede) $iv)$ et $ii)$ en tout point $P'$ d'un voisinage affine $U$ de Zariski de $P$ dans $Z$, ce qui, compte tenu de la remarque 3.2.2, implique $i)$. 
 
  \medskip D\'emontrons le point 2). On peut supposer qu'il n'y a qu'une seule pente $\rho >0$. On sait que l'on peut prendre pour $K'$ l'extension de $K$ engendr\'ee par les coefficients des $\phi_j$ (2.3.4), qui par hypoth\`ese sont entiers sur $A$. Il s'agit donc de d\'emontrer que les $A[\phi_{j,-k}]$ sont des extensions \'etales de $A$. Le th\'eor\`eme de puret\'e de Zariski-Nagata (\cf [{$ \bf SGA\,1$}, X.3.1) permet encore de remplacer $A$ par le compl\'et\'e d'un localis\'e en un premier de hauteur $1$ quelconque. On peut donc supposer de nouveau que $A$ est un anneau de valuation discr\`ete complet. 
  
    Commen\c cons par montrer que les $A[\phi_{j,-\rho}]$ sont \'etales sur $A$.  
 Pour tout \'el\'ement $\iota$ du sous-groupe d'inertie $I$ de $Gal(K'/K)$, et toute valeur propre $\phi_{j,-\rho}$ de $H(0)$ (notation de 3.2.2), $\iota(\phi_{j,-\rho})$ en est une autre $\phi_{i,-\rho}$, et $\phi_{i,-\rho}-\phi_{j,-\rho}$ est dans l'id\'eal maximal de $A'$. De ce que les segments de pente maximale $\rho$ co\"{\i}ncident pour $NP({{\sE}nd}\,M_{(P)})$ et pour $ NP({{\sE}nd}\,M) $, on d\'eduit que $\phi_{i,-\rho}=\phi_{j,-\rho}$. Donc les $\phi_{j,-\rho}$ sont fixes sous $I$, \ie sont dans une extension non ramifi\'ee de $A$ (contenue dans $A'$).  
 
En tordant $ M$ par $ L_{-\phi_{j,-\rho}x^{-\rho}}$, o\`u $ j$ r\'ealise le minimum des nombres rationnels
  $\min_{i\neq j}\, v_{A'}(\phi_{i,-\rho} -\phi_{j,-\rho}) $, et en it\'erant, on trouve que tous les $\phi_{j,-k}$ sont dans une extension non ramifi\'ee de $A$. \sq

 \bigskip\noindent  {\bf 3.6. Les diviseurs $D_{Z,\sigma}(M)$.} Rappelons que $Gal(F'/F)$ permute les composantes $M_{\bar\phi_j, F'}$ de la d\'ecomposition (3.2.1.1) en permutant les $\phi_j$ (\cf 2.3.1, 2.3.4). Soit $\sigma$ l'une des pentes non nulles de $M_F$ (\cf 2.4). Notons $J_{(\sigma)}$ l'ensemble des $j$ pour lesquels $\phi_j\in F'$ appara\^{\i}t dans la composante de $M_{F'}$ de pente $\sigma$. \'Ecrivons le terme de $\bar \phi_j$ de plus bas degr\'e en $x$ sous la forme $$\phi_{j,-\sigma}.x^{-\sigma} ,$$ avec $\phi_{j,-\sigma}\in K'$.  Pour $\sigma\in J_{(\sigma)}$, notons $\mu_j$ la dimension de $R_{j,F'}$ sur $F'$,  et posons $$\mu_{(\sigma)} = \sum_{j\in J_{(\sigma)}} \, \mu_j,$$ de sorte que $\sigma. \mu_{(\sigma)}$ est l'irr\'egularit\'e de la composante de $M_{F}$ de pente $\sigma$. 
 
   Il est alors clair que l'expression 
 $$\varphi_\sigma (x) = \prod_{j\in J_{(\sigma)}} \, (x^\sigma  -\phi_{j, -\sigma})^{\mu_j} \;\;\leqno{(3.6.1.1)}$$  est un polyn\^ome dans $K[x]$, de degr\'e \'egal \`a  $\sigma.\mu_{(\sigma)}$. 
 
 \smallskip Si tous les points de $Z=\Spec A$ sont semi-stables, on a m\^eme $\varphi_\sigma (x)\in A[x]$, et ce polyn\^ome d\'efinit un diviseur de Weil positif 
 $$D_{Z,\sigma}(M)= (\varphi_\sigma (x)) \in Div (\Spec A[x]).$$ Chacune de ses composantes est finie sur $Z$, et \'etale au-dessus des points stables.

  \bigskip \bigskip \bigskip    {\bf 4. \'Evitement des points tournants par \'ecla\-te\-ment.}   
 
 \bigskip Dans ce paragraphe, on s'int\'eresse \`a la structure for\-melle des modules \`a con\-ne\-xion in\-t\'e\-gra\-ble \`a plusieurs variables $x_1, x_2, \ldots, x_d$, \`a p\^oles le long du diviseur lisse  $ \, x_1 =0$. Le corollaire 4.3.3 jouera un r\^ole crucial au \S 7.

\medskip\noindent  {\bf 4.1. Connexions m\'eromorphes formelles sur $\frak X$.}   On suppose maintenant que $ A$ est une alg\`ebre noeth\'erienne int\`egre et for\-mellement lisse de dimension $d -1$ sur un corps alg\'ebriquement clos $k$.  En pratique $Z=\Spec  A$ sera une vari\'et\'e connexe lisse  sur  $k$, ou parfois le compl\'et\'e d'une telle vari\'et\'e en un point ferm\'e $P$. 

On consid\`ere un sch\'ema formel affine $x$-adique $${\frak X}\cong \Spf A[[x]],$$ dont on identifie la fibre sp\'eciale  ${\frak X}_{red}$ \`a $Z$ (sous-sch\'ema ferm\'e d\'efini par $x=0$).  On appelle faisceau des {\it fonctions m\'ero\-morphes for\-melles} sur $\frak X$ {\it \`a p\^oles le long de $Z$} le faisceau
$$ \sO_{\frak X}(\ast Z):= \sO_Z((x))= \sO_{\frak X}[{1\over x}].$$
 On appelle faisceau des {\it formes diff\'erentielles m\'ero\-morphes for\-melles} sur $\frak X$ {\it \`a p\^oles le long de $Z$} le $\sO_{\frak X}(\ast Z)$-module (localement libre de rang $d$) $$ \Omega^1_{\frak X}(\ast Z):=    \Omega^1_{\frak X}[{1\over x}]\;\cong\; \Omega^1_Z((x))\oplus  \sO_{\frak X}(\ast Z).dx.$$ Les formes diff\'erentielles \`a p\^oles logarithmiques le long de $Z$ en forment le sous-$\sO_{\frak X}$-module $$ \Omega^1_{\frak X}(\log Z) \cong \Omega^1_Z[[x]]\oplus  \sO_{\frak X} . {dx\over x}.$$ 
  \medskip On se donne par ailleurs un $A((x))$-module projectif de type fini ${\frak M}$ (vu comme $\sO_{\frak X}(\ast Z)$-module) muni d'une con\-ne\-xion {\it in\-t\'e\-gra\-ble} relative \`a $k$:
 $$\nabla:\; {\frak M}\to {\frak M}\otimes_{\sO_{\frak X}(\ast Z)} \Omega^1_{\frak X}(\ast Z) .$$
   Par abus de langage, on dira que ${\frak M}$ est une {\it con\-ne\-xion m\'ero\-morphe  for\-melle} sur $\frak X$ {\it \`a p\^oles le long de $ Z$}. 
    
     \medskip\noindent {\bf Remarque 4.1.1.} La projectivit\'e de ${\frak M}$ est automatique: c'est un fait g\'en\'eral, ind\'ependant de l'int\'egrabilit\'e, qui vient de ce que l'anneau diff\'erentiel $A((x))$ est simple (eu \'egard  \`a ${\partial}=x{d \over  dx}$ et aux $k$-d\'erivations de $A$), \cf [${\bf A1}$,  2.5.2.1].  
     
\medskip\noindent {\bf Remarques 4.1.2.} $1)$ Les con\-ne\-xions m\'ero\-morphes for\-melles de rang un \`a p\^oles le long de $ Z$  sont celles de la forme $ {\frak L}_\omega =  ( \sO_{\frak X}(\ast { Z}), \, \nabla(1)= \omega) ,\;\; \omega \in   \Omega^1_{\frak X}(\ast { Z}), \;d\omega =0. $   Avec les notations du paragraphe pr\'ec\'edent, on pose  $\langle \partial, \omega\rangle = \phi \in A((x)).$ 

\smallskip Une telle con\-ne\-xion ${\frak L}_\omega$ est r\'eguli\`ere si et seulement si $\omega \in   \Omega^1_{\frak X}(\log {  Z})$. Ajouter \`a $\omega$ une diff\'erentielle ferm\'ee dans $\Omega^1_{\frak X}(\log {  Z})$ revient donc \`a tordre ${\frak L}_\omega$ par une con\-ne\-xion r\'eguli\`ere de rang un, ce qui permet de se ramener au cas o\`u $\phi \in {1\over x}A[{1\over x}]$, de degr\'e en $1\over x$ \'egal au rang de Poincar\'e-Katz $\rho$ de ${\frak L}_\omega$ relativement \`a $\partial$ le long de $Z$ (\cf 2.5.3).  
   
 \smallskip\noindent   $2)$ Si $A$ est local complet de carac\-t\'e\-ris\-tique $0$, le lemme de Poincar\'e formel permet d'\'ecrire  
  $$ \omega= \epsilon {dx \over  x } +d(x^{-\rho}.\alpha),\leqno{(4.1.2.1)}$$ o\`u
  $ \epsilon  \in k,\; \alpha \in  \sO_{\frak X}  $. On peut alors \'ecrire la con\-ne\-xion ${\frak L}_\omega$ sous la forme $$ \,x^\epsilon e^{x^{-\rho}.\alpha}.\sO_{\frak X}(\ast Z) \leqno{(4.1.2.2)}$$  
    et on a $$\partial(\alpha) - \rho \alpha = (\phi+ \epsilon) x^\rho.\leqno{(4.1.2.3)}$$  
    Dans le cas o\`u $\phi \in {1\over x}A[{1\over x}]$, on a $\epsilon =0$, et on en d\'eduit qu'on peut choisir de mani\`ere unique $\alpha\in A[x]$,  nul ou de degr\'e $\leq \rho$ en $x$. Un point $P\in \Spec A$ est stable si et seulement si $\alpha(0)$ ne s'annule pas en $P$.
  
 \smallskip\noindent   $3)$ Si $A$ est local complet de carac\-t\'e\-ris\-tique $0$, il est bien connu que {\it les con\-ne\-xions m\'ero\-morphes for\-melles r\'egu\-li\`eres (\`a p\^oles le long de $ Z$)  sont extensions it\'er\'ees de con\-ne\-xions m\'ero\-morphes  for\-melles r\'egu\-li\`eres de rang un.}
 
\noindent Esquissons une d\'emonstration dans le cas $d=2$, c'est-\`a-dire $A\cong k[[x_2]]$. Soit $\bf n$ une base de $M$ dans laquelle la matrice $G $ de $ \partial$ est \`a coefficients dans $k[[x_2,x]]$ (il en existe d'apr\`es 3.3.2).   
 L'int\'egrabilit\'e entra\^{\i}ne que les valeurs propres de $G$ sont dans $k$. Un argument classique de changement de base (en commen\c cant par des cisaillements pour que ces valeurs propres soient distinctes modulo $\Z$) permet ensuite de supposer $G $ \`a coefficients dans $k[[x_2]]$. Par int\'egrabilit\'e, la matrice de $\nabla(x_2{ \partial \over {\partial x_2}} )$ est aussi \`a coefficients dans $k[[x_2]]$. On conclut par la th\'eorie \`a une variable $x_2$ (voir aussi [$\bf S2$ III.2.1.1]).

\bigskip\noindent  {\bf 4.2. Cas o\`u les points de $Z$ sont semi-stables ou stables.} Si tous les points de $Z$ sont semi-stables pour $\frak M$, on a une d\'ecomposition {\rm (3.2.1.2)} de l'image inverse ${\frak M}'$ de $\frak M$ sur $\Spf A'[[x^{1/e}]]$ compatible \`a l'action de $\partial = x{\partial\over \partial x}$, ainsi qu'une d\'ecomposition de $\frak M$ suivant les pentes (3.2.5). En outre, si tous les points de $Z$ sont stables pour $\frak M$, alors d'apr\`es 3.4.1, $A'$ est \'etale sur $A$, donc formellement lisse sur $k$.  

  \proclaim Proposition 4.2.1. Si tous les points de $Z$ sont semi-stables, la d\'ecomposition de $\frak M$ suivant les pentes est une d\'ecomposition de connexions m\'eromorphes formelles. En outre, si tous les points de $Z$ sont stables, alors la d\'ecomposition {\rm (3.2.1.2)} de ${\frak M}'$est une d\'ecomposition de connexions m\'eromorphes formelles. \par

\dem Le point est l'int\'egrabilit\'e, qui r\'esulte de 2.5.1,1) et 2.5.2. \sq 

  \smallskip   La remarque 2.5.3 s'applique aussi dans ce contexte.   
  
  \medskip Consid\'erons le dual $\sT_{{\frak X}, Z}$ de $\Omega^1_{\frak X}(\log Z)$, et le morphisme $$\iota: \sO_Z\to ( \sT_{{\frak X}, Z})_{\mid Z}$$ dual de l'application r\'esidu. Le choix d'une d\'erivation continue $\partial\in  \Gamma(\sT_{{\frak X}, Z})$ dont la restriction \`a $Z$ est $\iota(1)$ identifie le fibr\'e normal $N_{Z}\frak X$ \`a  $\Spec A[x]$, $\partial$ s'identifiant alors \`a $x{\partial \over \partial x}$, \cf [$\bf AB$, I.1.4.5]. 
  
  Supposons que les points de $Z$ soient semi-stables pour $\frak M$, et soit $\sigma$ une pente non nulle de $\frak M$. Alors le diviseur de Weil positif $D_{Z,\sigma}({\frak M})$ (\cf 3.6) s'interpr\`ete comme {\it diviseur sur  $N_{Z}\frak X\,$}. 

 \proclaim Lemme 4.2.2. Cette interpr\'etation est canonique, \ie ne d\'epend pas du choix de $\partial$. \par 
\dem D'apr\`es 3.2.5, on peut supposer que $\frak M$ a une seule pente $\sigma$. En outre, par descente, on se ram\`ene au cas $Z=Z',\, e=1$. Toute autre d\'erivation continue de $\sO(\frak X )$ ayant les m\^emes propri\'et\'es s'\'ecrit $u.\partial + vx\partial'$ o\`u $u$ est une unit\'e principale $x$-adique (\ie vaut $1$ modulo l'id\'eal $\sI_Z$ de $Z$), $v$ est entier $x$-adique, et $\partial'$ commute \`a $\partial$. D'apr\`es 2.5.1 2) (appliqu\'e \`a la composante ${\frak M}_{(\sigma)}$ de $\frak M$ de pente $\sigma$),  $vx\partial'$ ne contribue pas au terme de plus bas degr\'e en $x$ dans $\phi_j$.
Ainsi, l'inverse $ x^\sigma\phi_{j,-\sigma}^{-1}$ de ce terme est un \'el\'ement bien d\'efini de $\Gamma(\sI^\sigma_Z/\sI^{\sigma+1}_Z) \subset \sO(N_{Z}\frak X\,),$ et $D_{Z,\sigma}({\frak M})$ un diviseur bien d\'efini sur $N_{Z}\frak X\,$.
  \sq

\medskip\noindent {\bf Remarque  4.2.3.} Dans le cas particulier o\`u $Z$ est une courbe sur $k$ de carac\-t\'e\-ris\-tique nulle, ce $1$-cycle $D_{Z,\sigma}({\frak M})\subset N_{Z}\,\frak X$ co\"{\i}ncide avec celui que Sabbah ([$\bf S1$], [$\bf S2$, I.3.1])  a construit  en projetant sur $N_{Z}\,\frak X$ l'intersection, dans le cotangent $T^\ast N_{Z}\,\frak X$, du $2$-cycle micro-carac\-t\'e\-ris\-tique $CCh'_{Z,\sigma}({\frak M})$ (au sens d'Y. Laurent) et du $2$-cycle image inverse du point $1$ par le morphisme d'Euler $T^\ast N_Z\,{\frak X}\to {\bf A}^1$.  Nous n'aurons pas besoin de ce fait\footnote{$^{(9)}$}{\rm \sm En notant $(x_1,x_2,\xi_1,\xi_2)$ les co\-or\-don\-n\'ees sur $T^\ast N_{Z}\,\frak X$ (avec $x_1= x $, et $x_2$ une co\-or\-don\-n\'ee sur $Z$), le morphisme d'Euler est $x_1\xi_1$, et on peut montrer, en utilisant la remarque 4.1.2, 3), que dans le cas o\`u $\sigma$ est entier, $CCh'_{Z,\sigma}({\frak M})$ est combinaison des surfaces 
 $$  x_1^{\sigma+1} \xi_1  =  \phi_{ j, -\sigma} \,, \;\;\;\; \sigma . x_1^{\sigma}  \xi_2  +  { { \partial \phi_{ j, -\sigma} }\over{ \partial x_2 }}   = 0, $$ affect\'ees de la multiplicit\'e $\mu_j$.}.

\bigskip\noindent  {\bf 4.3. Stabilisation par \'eclatement.}  Un germe de courbe ${\frak C}$ sur ${\frak X}$ est un point ferm\'e de $\Spec A[[x]]\setminus Z$, ou, ce qui revient au m\^eme ([$\bf BoL$ 3.4]), un sous-sch\'ema formel $x$-adique ferm\'e de $\frak X$ du type $\Spf B$, o\`u $B$ est local int\`egre de dimension $1$. 

  On se donne un germe ${\frak C}$ de courbe lisse  (\ie $\cong \Spf B$ avec $B$ comme ci-dessus et r\'egu\-lier) sur ${\frak X}$ coupant $Z$ transversalement en un point ferm\'e $P$.   Un tel germe cor\-res\-pond \`a une section du morphisme structural ${\frak X}\to \Spf k[[x]]$, et quitte \`a changer l'isomorphisme ${\frak X}\cong \Spf A[[x]] $, on se ram\`ene au cas o\`u l'id\'eal de ${\frak C}$ est engendr\'e dans $\sO_{\frak X}$ par des g\'en\'erateurs de l'id\'eal de $P$ dans $Z$. 
 
 Consid\'erons la suite d'\'ecla\-te\-ments formels (\cf [{$\bf BoL$}]) $${\frak X}_{n+1}\to {\frak X}_{n}\to \cdots \to {\frak X}_{0}={\frak X}$$ du point $P_n$ intersection du transform\'e strict ${\frak C}_n$ de  ${\frak C}$ et de la composante $E_n$ du diviseur exceptionnel rencontrant ${\frak C}_n$, en partant de $ {\frak C}_{0}= {\frak C},  P_{0} =P$.
 
\noindent  Soit $Z_n'$ la r\'eunion des composantes de l'image inverse de $Z$ distinctes de $E_n$, et posons 
 $${\frak X}^0_n=  {\frak X}_n\setminus Z_n',\;\;\; E_n^0 = E_n\setminus (Z_n'\cap E_n).$$ 
  Alors ${\frak X}^0_n$ est un $k[[x]]$-sch\'ema formel affine, et $({\frak X}^0_n)_{red}= E_n^0$ (vari\'et\'e affine connexe lisse sur $k$).   
  
     Pour tout $n\geq 1$, on dispose alors de l'image inverse ${{\frak M}}_n$ de ${{\frak M}}$ sur ${\frak X}^0_n$, qui est une con\-ne\-xion m\'ero\-morphe  for\-melle sur ${\frak X}^0_n$  \`a p\^oles le long de $E_n^0$.

 \proclaim Th\'eor\`eme 4.3.1. Pour $n\geq   \rho({{\frak M}}) +1$, le point $P_n\in E_n^0$ est un point stable pour ${{\frak M}}_n$. \par 
 
 \dem  Soit $P_n'$ un point stable de $E_n^0$, et soit ${\frak C}'_n$ le germe de courbe sur ${\frak X}^0_n$ coupant $E_n^0$ transversalement en $P'_n$. L'image de $P_n'$ dans ${\frak X}_m,\, m<n,$ n'est autre que $P_m$. 
 Comme ${\frak C}$ coupe $x=0$ transversalement en $P$, on voit que pour $0<m<n$, tant ${\frak C}_m$ que ${\frak C}'_m$ coupent $E_m^0$ transversalement en $P_m$.
 
D'apr\`es 3.4.1, il s'agit de montrer, sous l'hypoth\`ese $n\geq   \rho({{\frak M}}) +1$, que $$(\ast\ast)_n\;\;\; NP_{P_n}({{\frak M}}_{n \mid {\frak C}_n })=  NP_{P'_n}({{\frak M}}_{n \mid {\frak C}_n'}),\;\; NP_{P_n}( {{\sE}nd}\,{{\frak M}}_{n \mid {\frak C}_n })=  NP_{P'_n}( {{\sE}nd}\,{{\frak M}}_{n \mid {\frak C}_n'}).$$ 
  
 Nous allons montrer plus: 
 
 \proclaim Lemme 4.3.2. Si $\, \rho({{\frak M}})\leq m<n,\,$ alors l'isomorphisme compos\'e de ${k}$-alg\`ebres
  $$\iota_m \;:\;\; {\cal O}_{{\frak C}_m,P_m}= {k}[[x ]] \hookrightarrow {\cal O}_{{\frak X}_m, P_m}\to  {\cal O}_{{\frak C}'_m, P_m} $$ induit un isomorphisme {\it horizontal}
 $${{\frak M}}_{m \mid {\frak C}_m}\cong  {{\frak M}}_{m \mid {\frak C}'_m} $$
(d'o\`u aussi un tel isomorphisme pour tout $m\geq 0$, et en par\-ti\-cu\-lier $(\ast\ast)_n$ pour $m=n$). \par

\dem Soit $x_2,\ldots, x_d$ un syst\`eme de co\-or\-don\-n\'ees \'etales en $P\in Z$, et prenons 
$$x'_2={ {x_2}\over {x^m}},  \ldots, x'_d= {{x_d}\over {x^m}}$$ pour syst\`eme de co\-or\-don\-n\'ees \'etales en $P_m\in   E_m^0$. Alors, comme on l'a vu ci-dessus, on peut supposer que $\frak C$ est d\'efini, sur $\frak X$,  par l'annulation  des $x_i$,  de sorte que ${\frak C}_m$ est d\'efinie par l'annulation des $x'_i$. On a $\hat A_P\cong k[[{\bf x}]] = {k}[[x_2,\ldots, x_d]]$. 
  
    Il suit de 2.5.1, 2) (avec $\partial'= \partial/\partial x_i$) et de 2.1.3 $ii)$ (avec $\partial= \partial/\partial x_i$) qu'il existe un ${k}[[x ,{\bf x}]]]$-r\'eseau ${\Bbb M}$ de $ \hat {{\frak M}}_P= {{\frak M}} \otimes_{\sO_{\frak X}{\ast Z}}{k}[[x ,{\bf x}]]][{1\over x }]$ et un entier $N\geq 0$, tels que pour tout multi-indice ${\bf j}= (j_2, \ldots, j_d) $ de longueur $\mid {\bf j}\mid = j_2+\ldots j_d$,
  $$\nabla({\partial^{\bf j} \over \partial  
{\bf x}^{\bf j}}) ({\Bbb M}) \;\subset\; x^{-[\rho( {{\frak M}}) \mid {\bf j}\mid]-N}{\Bbb M}\;\subset \; x^{-  m  \mid {\bf j}\mid -N} {\Bbb M} , $$ donc $$\nabla({\partial^{\bf j} \over \partial  
{\bf x'}^{\bf j}})(\pi^\ast {\Bbb M}) \subset  x^{-N}\pi^\ast {\Bbb M}. $$
Si $\car k= 0$,  cela permet de d\'efinir l'op\'erateur 
$$\Pi= \sum_{\bf j}  
 (-1)^{\mid {\bf j}\mid} {1\over {\bf j}!}  \nabla({\partial^{\bf j} \over \partial  
{\bf x'}^{\bf j}})\; \; : \;
\pi^\ast {\Bbb M} \longrightarrow x^{-N}\pi^\ast {\Bbb M}.\; $$
  On v\'erifie que pour tout $f\in k[[x,{\bf x'}]]$ et tout $e\in \pi^\ast {\Bbb M}$ , $\Pi(f e)= (f_{\mid x'_i=0}) \Pi(e)$, ce qui montre que $\Pi$ induit un homomorphisme ${k}((x))$-lin\'eaire 
     $$ {{\frak M}}_{\mid {\frak C} } \cong \pi^\ast {{\frak M}}_{\mid {\frak C}_m } = (\pi^\ast {\Bbb M}/({\bf x'} ))[{1\over x}] \,\to  \, (\pi^\ast {\Bbb M})[{1\over x}] = \pi^\ast { \hat {{\frak M}}_P}.$$
  Tout comme $\Pi$, il commute \`a $\partial \over \partial x$. En outre, $\Pi(e)
\equiv e \; {\rm mod}\;
 ({\bf x'})$, d'o\`u il suit que l'homomorphisme ${k}((x ))$-lin\'eaire compos\'e  
$$ {{\frak M}}_{\mid {\frak C }} \cong  (\pi^\ast { {{\frak M}}})_{\mid  C_m  }\,\to  \, \pi^\ast {\hat {{\frak M}}_P}  \,\to \, (\pi^\ast {\hat {{\frak M}}_P})_{\mid  C_m^\prime } \cong   {{\frak M}}_{\mid {\frak C}'}$$ n'est autre que l'isomorphisme induit par $\iota_m $, qui est donc horizontal. 

\medskip\noindent Si $\car k\neq 0$, le m\^eme argument s'applique \`a condition de remplacer $k[[{\bf x'}]]$ par l'enveloppe \`a puissances divis\'ees (le fait que l'homomorphisme de $k[[{\bf x'}]]$ dans l'enveloppe \`a puissances divis\'ees ne soit pas injectif ne pose pas probl\`eme). \sq

 \proclaim Corollaire 4.3.3. Supposons que $\,d \,(= \dim Z+  1) = 2$. Soit ${\frak X}_1\to \frak X$ l'\'ecla\-te\-ment formel des points tournants de $Z$ (s'il en est), et construisons une suite d'\'ecla\-te\-ments formels $${\frak X}_{n+1}\to {\frak X}_{n}$$ en \'eclatant les points tournants du diviseur exceptionnel $E_n$ dans ${\frak X}_n$.  Alors cette suite s'arr\^ete \`a un pas $ \leq  \rho({{\frak M}})+1 $. 
  \par
   \dem  Soit $P_n$ un point ferm\'e de non-croisement de $ E_n $, et soit ${\frak C}_n$ le germe de courbe   sur ${\frak X}_n$ coupant $E_n$ transversalement en $P_n$. Comme la projection $\pi_n:\, {\frak X}_n\to \frak X$ est une suite d'\'ecla\-te\-ments formels de points tournants (qui ne sont pas des points de croisement),  ${\frak C} = \pi_n({\frak C}_n)$  est un germe de courbe sur ${\frak X}$ coupant $Z$ transversalement en un point ferm\'e $P$, et ${\frak C}_n$ est le transform\'e strict de $\frak C$. Pour  $n\leq  \rho({{\frak M}})+1 $, le th\'eor\`eme 4.3.1 implique que $P_n$ est stable, \ie n'est pas un point tournant.\sq

 \bigskip \bigskip   

\bigskip {\bf 5. Points de croisement, et stabilisation par \'ecla\-te\-ment. }  

\medskip On s'int\'eresse maintenant \`a la structure for\-melle des modules \`a con\-ne\-xion in\-t\'e\-gra\-ble  \`a deux variables $y_1, y_2$, \`a p\^oles le long du diviseur \`a croisements normaux  $ \, y_1y_2=0$. 
Le lemme 5.3.1 jouera un r\^ole important au \S 6. Le th\'eor\`eme 5.4.1 sera crucial aux  \S\S$\,$   6 et 7.

\medskip\noindent  {\bf 5.1. Connexions m\'eromorphes formelles sur $\frak Y$.}   
 On consid\`ere le sch\'ema formel affine $(y_1,y_2)$-adique ${\frak Y}= \Spf  k[[y_1,y_2]] $. On note $Q$ le point ${\frak Y}_{red}$ (d\'efini par $y_1=y_2=0$). 

   \noindent  Un germe de courbe sur ${\frak Y}$ est un point ferm\'e de $\Spec A[[y_1,y_2]]\setminus Q$, ou, ce qui revient au m\^eme [$\bf BoL$ 3.4], un sous-sch\'ema formel $(y_1,y_2)$-adique ferm\'e de $\frak Y$ du type $\Spf B$, o\`u $B$ est local int\`egre de dimension $1$. 

On se donne des germes de courbes ${\frak Z}_1, \ldots, {\frak Z}_t$, d\'efinis par des irr\'eductibles $f_1,\ldots f_t\in k[[y_1,y_2]]$,  et on pose ${\frak Z}=\cup {\frak Z}_i,\; f=\prod f_i$.  Pour employer un langage g\'eom\'etrique comme en 4.1, on note  $$ \sO_{\frak Y}(\ast {\frak Z}):=  k[[y_1,y_2]][{1\over f}]  ,$$ et on appelle  {\it fonctions m\'ero\-morphes for\-melles} sur $\frak Y$ {\it \`a p\^oles le long de $\frak Z$} ses \'el\'ements; on note 
$$  \Omega^1_{\frak Y}(\ast {\frak Z}):=  {\Omega^1}^{cont}_{k[[y_1,y_2]][{{1}\over {f}}]} ,$$  et on appelle  {\it formes diff\'erentielles m\'ero\-morphes for\-melles} sur $\frak Y$ {\it \`a p\^oles le long de $\frak Z$} ses \'el\'ements.

  On se donne par ailleurs un $\sO_{\frak Y}(\ast {\frak Z})$-module projectif de type fini ${\frak N}$ muni d'une con\-ne\-xion {\it in\-t\'e\-gra\-ble} relative \`a $k$:
 $$\nabla:\; {{\frak N}}\to {{\frak N}}\otimes_{\sO_{\frak Y}(\ast {\frak Z}) } \Omega^1_{\frak Y}(\ast {\frak Z}) .$$ Par abus de langage, on dira que ${\frak N}$ est une {\it con\-ne\-xion m\'ero\-morphe  for\-melle} sur $\frak Y$ {\it \`a p\^oles le long de $\frak Z$}.

\smallskip \noindent On supposera toujours que le rang $\nu$ de ${{\frak N}}$ est strictement inf\'erieur \`a la carac\-t\'e\-ris\-tique de $k$ si celle-ci est non nulle.

\bigskip\noindent   
Le cas qui nous int\'eresse dans ce paragraphe est celui de $t=2,\, f_1=y_1, f_2= y_2$.  On consid\`ere dans ce cas le sous-$\sO_{\frak Y}$-module  
$$ \Omega^1_{\frak Y}(\log {\frak Z}) = \sO_{\frak Y}{dy_1\over y_1}\oplus \sO_{\frak Y}{dy_2\over y_2}$$ form\'e des formes diff\'erentielles logarithmiques le long de $\frak Z$.  

 \medskip\noindent {\bf Remarque 5.1.1.} La projectivit\'e de ${\frak N}$ est automatique: elle vient de ce que l'anneau diff\'erentiel $k[[y_1,y_2]][{{1}\over {y_1y_2}}]$ est simple (eu \'egard \`a $\partial_i = y_i {d\over dy_i}$, $i=1,2$),  \cf [${\bf A1}$,  2.5.2.1]. 

\noindent  Elle implique la libert\'e, car l'anneau $k[[y_1,y_2]][{{1}\over {y_1 y_2}}]$ est principal.

 \medskip\noindent {\bf Remarques 5.1.2.}   $1)$ Les con\-ne\-xions m\'ero\-morphes for\-melles de rang un \`a p\^oles le long de $\frak Z  $  sont celles de la forme  $\,{\frak L}_\omega =  ( \sO_{\frak Y}(\ast {\frak Z  }), \, \nabla(1)= \omega) ,\;\; \omega \in   \Omega^1_{\frak Y}(\ast {\frak Z  }), \;d\omega =0.  $  On \'ecrit $$\omega = \psi_1 {dy_1\over y_1}+  \psi_2 {dy_2\over y_2},\; \;\partial_1(\psi_2)= \partial_2(\psi_1). \leqno{(5.1.2.1)}$$
  Une telle con\-ne\-xion ${\frak L}_\omega$ est r\'eguli\`ere si et seulement si $\omega \in   \Omega^1_{\frak Y}(\log { \frak Z})$. 
  
\smallskip\noindent  $2)$  En carac\-t\'e\-ris\-tique $0$, on peut \'ecrire  
  $$ \omega= \epsilon_1 {dy_1 \over  y_1 }+ \epsilon_2 {dy_2 \over  y_2 } +d(y_1^{-\rho_1}.y_1^{-\rho_2}\alpha),\leqno{(5.1.2.2)}$$ o\`u
  $ \epsilon_1,\epsilon_2  \in k,\; \alpha \in  \sO_{\frak Y}, $  et $\rho_i$ est le rang de Poincar\'e-Katz de ${\frak L}_\omega$ relativement \`a $\partial_i $ le long de $y_i=0,\, i=1,2$.  On peut alors \'ecrire ${\frak L}_\omega$ sous la forme $$ \,y_1^{\epsilon_1}y_2^{\epsilon_2} e^{y_1^{-\rho_1}.y_1^{-\rho_2}.\alpha}.\sO_{\frak X}(\ast {\frak Z}) \leqno{(5.1.2.3)}$$ et on a $$\partial_i(\alpha) -\rho_i\alpha = (\psi_i+\epsilon_i)y_1^{\rho_1} y_2^{\rho_2}.\leqno{(5.1.2.4)}$$ 
   \smallskip\noindent   $3)$   En carac\-t\'e\-ris\-tique $0$, il est bien connu que les con\-ne\-xions m\'ero\-morphes for\-melles r\'egu\-li\`eres (\`a p\^oles le long de $\frak Z$)  sont extensions it\'er\'ees de con\-ne\-xions m\'ero\-morphes  for\-melles r\'egu\-li\`eres de rang un. Nous n'aurons pas besoin de ce r\'esultat (dont la preuve est un peu plus subtile que celle de 4.1.2, 3)).

 \bigskip\noindent  {\bf 5.2. Points de croisement semi-stables.}   
 On pose   $F_1 = k((y_2))((y_1))$ (\resp $F_2= k((y_1))((y_2))$), muni de la valuation $y_1$-adique (\resp $y_2$-adique). On note $\rho_1 = \rho_1({{\frak N}})$ (\resp $\rho_2 = \rho_2({{\frak N}})$) le rang de Poincar\'e-Katz du module diff\'erentiel ${{\frak N}}_{F_1}$ (\resp ${{\frak N}}_{F_2}$). 

\medskip Quitte \`a remplacer $y_1$ et $y_2$ par $y_1^{e}$ et $y_2^{e}$ avec $e\mid \nu!$, on dispose alors des d\'e\-com\-po\-si\-tions de Turrittin-Levelt, pour $i=1,2$ (avec les notations de 2.3.1): 
$${{\frak N}}_{F_i} = \oplus_j \; {{\frak N}}_{\bar\phi_{i,j} ,F_i}\,  \leqno (5.2.1.1)_i $$ $$\;\; {{\frak N}}_{\bar\phi_{i,j} ,F_i}\cong {{\frak L}}_{\phi_{i,j} ,F_i}\otimes {{\frak R}}_{i,j,F_i} .$$
Les $\phi_{1,j}\in {1\over y_1}k((y_2))[{1\over y_1}]$ sont de degr\'e $\geq -\rho_1  $ en $y_1$, l'un des $\phi_{1,j}$ au moins \'etant exactement de degr\'e $-\rho_1 $ (\resp $\phi_{2,j}\in {1\over y_2}k((y_1))[{1\over y_2}]$ de degr\'e $\geq -\rho_2  $ en $y_2$, l'un des $\phi_{2,j}$ au moins \'etant exactement de degr\'e $-\rho_2 $).

 \bigskip Ce paragraphe est d\'evolu \`a l'\'etude de la (double) question suivante.

  \medskip\noindent  {\bf Question 5.2.1.} {\it $i)$ A-t-on, quitte \`a remplacer $y_1$ et $y_2$ par $y_1^{e}$ et $y_2^{e}$ avec $e\mid \nu!$, une d\'ecom\-po\-si\-tion de ${{\frak N}} $ lui-m\^eme 
   $${{{{\frak N}}}}  =  \oplus_{h}  \; {{\frak N}}_{  \bar\omega_{h} }\,   \leqno (5.2.1.2)$$
avec\footnote{$^{(10)}$}{ \rm \sm il n'y a nulle part risque de confusion entre l'indice ${h}$ et le corps de base.}  $$\;\;  {{\frak N}}_{ \bar\omega_{h} } \cong  \;  {\frak L}_{  \omega_{h} } \otimes  {\frak R}_{ {h}},\;  $$   
o\`u 
\medskip\noindent - les $\omega_{h} \in  \Omega^1_{\frak Y}(\ast {\frak Z})$ d\'esignent des formes diff\'erentielles m\'ero\-morphes for\-melles {\rm ferm\'ees} ($d\omega_{h} =0$), dont les classes $\bar \omega_{h} $ modulo $\Omega^1_{\frak Y}(\log {\frak Z})$ sont {\rm deux \`a deux distinctes},
   \medskip\noindent -  $ {\frak L}_{\omega_{h} }  $ est le $\sO_{\frak Y}(\ast {\frak Z})$-module \`a con\-ne\-xion in\-t\'e\-gra\-ble de rang $1$ attach\'e \`a $\omega_{h} $,
 \medskip\noindent -  $ {\frak R}_{h} $ est un $\sO_{\frak Y}(\ast {\frak Z})$-module \`a con\-ne\-xion in\-t\'e\-gra\-ble {\rm r\'egu\-lier} (relativement \`a $\partial_1 = y_1 {d\over dy_1}$ et $\partial_2 = y_2 {d\over dy_2}$)? 
 \medskip  
  $ii)$ Si oui, quels sont, pour $i=1,2$, les liens entre les $\psi_{i,{h}} := \langle \partial_i, \bar \omega_{h} \rangle$ et les $\bar \phi_{i,j}$, \resp entre la d\'e\-com\-po\-si\-tion $(5.2.1.2)$ et les d\'e\-com\-po\-si\-tions  $(5.2.1.1)_i$ ? }
  
   \bigskip 
  
    \medskip La r\'eponse \`a 5.2.1.$i)$ est positive en rang $\nu=1$. En revanche, sans hypoth\`ese suppl\'ementaire, la question a une r\'eponse n\'egative:

  \medskip\noindent {\bf Contre-exemple 5.2.2.} Le contre-exemple 3.1.3 sert aussi bien ici; en fait il provient d'une con\-ne\-xion in\-t\'e\-gra\-ble, donn\'ee,  
  dans une base $ (m_1,m_2)$,  par
    $$\nabla(\partial_i)(m_1)= (-1)^{i-1} {y_2\over y_1}m_1,\; \nabla(\partial_i)(m_2)= (-1)^{i }  m_1,\; (i=1,2) .$$
 Une d\'e\-com\-po\-si\-tion $(5.2.1.2)$ ne peut avoir lieu, m\^eme apr\`es ramification, et m\^eme si l'on n'impose pas aux $\bar \omega_{h} $ d'\^etre deux \`a deux distinctes. En effet, on aurait compatibilit\'e \`a la d\'e\-com\-po\-si\-tion de Turritin-Levelt de ${\frak N}_{F_1}$ (\cf ci-dessous 5.3.1): 
 $${\frak N}_{F_1}\cong   {\frak L}_{{{y_2}\over {y_1}}, F_1}  \oplus   {\frak L}_{{0}, F_1} ,$$ o\`u
 $e^{y_2/y_1} m_1$ forme une base horizontale du facteur ${\frak L}_{{{y_2}\over {y_1}}, F_1} $, tandis que,
 en notant 
 $${ \epsilon}(x)=\sum_0^\infty (-1)^nn!x^{n+1} $$ la s\'erie d'Euler [${\bf E}$],
 $ \epsilon(y_1/y_2)m_1+ m_2$ forme une base horizontale du facteur ${\frak L}_{{0}, F_1} $. 
Or, en carac\-t\'e\-ris\-tique $0$, $ \epsilon(y_1/y_2)\notin k[[y_1,y_2]][{{1}\over {y_1 y_2}}]$. 
 
  \medskip
\proclaim D\'efininition 5.2.3.  Nous dirons que le point de croisement $Q$ est {\rm semi-stable} pour ${\frak N}$ si la question 5.2.1.$i)$ a une r\'eponse positive.\par

\medskip\noindent {\bf Remarque 5.2.4.} Si $\frak N$ n'a pas de p\^ole le long de $y_2=0$, cette d\'efinition est compatible \`a 3.2.4.
   
   \proclaim Lemme 5.2.5.   Si $Q$ est un point de croisement semi-stable, la d\'e\-com\-po\-si\-tion $(5.2.1.2)$ est unique.
   \par
  \noindent  Comme en 2.3.1.2), cela d\'ecoule de ce que si $\omega$ et $\omega'$ sont distincts modulo $\Omega^1_{\frak Y}(\log {\frak Z})$, alors 
   $ \, Hom_\nabla({\frak L}_\omega\otimes {\frak R}, {\frak L}_{\omega'}\otimes {\frak R}')= 0 \, $  car ${\frak L}_{\omega'-\omega}$ est irr\'egu\-lier. $\square$

  \bigskip\noindent  {\bf 5.3. } Supposons que $Q$ soit un point de croisement semi-stable, et examinons la question 5.2.1.$ii)$. 
    Observons qu'on a une d\'e\-com\-po\-si\-tion canonique, stable sous $\partial_i$: 
$${{k[[y_1,y_2]][{{1}\over {y_1 y_2}}]}\over{ k[[y_1,y_2]]}}\;\cong \; {1\over y_1}k[[y_2]][{1\over y_1}]\;\oplus \;{{{{1}\over {y_1 y_2}}  k[ {1\over y_1},{1\over y_2}]} }\;\oplus\; {1\over y_2}k[[y_1]][{1\over y_2}], \leqno (5.3.1.1)$$
ce qui donne, en regroupant les deux premiers termes:
$${{k[[y_1,y_2]][{{1}\over {y_1 y_2}}]}\over{ k[[y_1,y_2]]}}\;\cong \; {1\over y_1}k((y_2))[{1\over y_1}] \;\oplus \;{1\over y_2}k[[y_1]][{1\over y_2}]. \leqno (5.3.1.2)$$
Le morphisme naturel 
$$\varpi_1:\;\; {{k[[y_1,y_2]][{{1}\over {y_1 y_2}}]}\over{ k[[y_1,y_2]]}}\;\to \; {{k((y_2))((y_1))}\over{ k((y_2))[[y_1]]}}\;\cong \; {1\over y_1}k((y_2))[{1\over y_1}] \leqno (5.3.1.3)$$ cor\-res\-pond \`a la projection sur le premier facteur dans la d\'e\-com\-po\-si\-tion $(5.3.1.2)$, qui est compatible aux $\partial_i$. Sym\'etriquement pour  $$\varpi_2:\;\;  {{k[[y_1,y_2]][{{1}\over {y_1 y_2}}]}\over{ k[[y_1,y_2]]}}\;\to \; {{k((y_1))((y_2))}\over{ k((y_1))[[y_2]]}}\;\cong \; {1\over y_2}k((y_1))[{1\over y_2}]. \leqno (5.3.1.4)$$ 
 
 Compte tenu de l'unicit\'e de la d\'ecomposition de Turrittin-Levelt  (2.3.1, 2)), on obtient:

 \proclaim Lemme 5.3.1.  Pour $i=1,2$, la d\'ecomposition   $$ {\frak N}   =  \oplus_{h}  \; {{\frak N}}_{  \bar\omega_{h} }\,,\;\; {{\frak N}}_{ \bar\omega_{h} } \cong  \;  {\frak L}_{  \omega_{h} } \otimes  {\frak R}_{h}, \leqno (5.3.1.7)$$
  induit par tensorisation avec $F_i$ une d\'ecomposition qui {\it raffine} celle de Turrittin-Levelt   $$ {\frak N}_{\bar\phi_{i,j},F_i} = \bigoplus_{{h},\, \varpi_i(\psi_{i,{h}}) =\bar\phi_{i,j}}\, {{\frak N}}_{  \bar\omega_{h} }\otimes F_i . \;\; \square\;  \leqno (5.3.1.5)_i$$\par

   \medskip\noindent {\bf Remarque 5.3.2.}  Si ${\frak N}$ provient d'une si\-tua\-tion alg\'ebrique ou analytique,  il en est de m\^eme des $\bar\phi_{i,j}$, et donc aussi de m\^eme des $\bar\omega_j$. Ceci donne une preuve simplifi\'ee de [$\bf S2$, I.2.4.4, 2.4.5].

\bigskip\noindent  {\bf 5.4. Semi-stabilisation des points de croisement par \'ecla\-te\-ment.} Dans l'exemple 5.2.2, $Q$ n'est pas semi-stable, mais on constate que le probl\`eme dispara\^{\i}t apr\`es \'ecla\-te\-ment de $Q$. C. Sabbah a d\'emontr\'e qu'il s'agit l\`a d'un ph\'enom\`ene g\'en\'eral. 

\smallskip Consid\'erons une suite d'\'ecla\-te\-ments formels   
$$\pi:\,{\frak Y}'\to \frak Y$$ d'abord de $Q$ puis de points de croisement des  diviseurs exceptionnels successifs. On sait qu'une telle suite d'\'ecla\-te\-ments, dite {\it torique}, cor\-res\-pond \`a un \'eventail r\'egu\-lier du premier quadrant de $\R^2$.  Les cartes toriques sont isomorphes \`a $\A^2$ (convenablement compl\'et\'e). La trace de l'image inverse de $\frak Z$ dans une telle carte est la r\'eunion des axes de co\-or\-don\-n\'ees. Dans celle associ\'ee au c\^one d'ar\^etes passant par $(a,b)\in \N^2$ et par $(c,d)\in \N^2$, avec $ad-bc=1$, les co\-or\-don\-n\'ees adapt\'ees $(y'_1,y'_2)$ sont donn\'ees par $$y_1  = (y'_1)^a(y'_2)^b,\; y_2 = (y'_1)^c(y'_2)^d.\leqno{(5.4.1.1)}$$ Posons $\,\omega =  \phi_1{dy_1\over y_1}+ \phi_2{dy_2\over y_2}$. Alors on a $\,\pi^\ast\omega=  \phi'_1{dy'_1\over y'_1}+ \phi'_2{dy'_2\over y'_2},$ o\`u $$\phi'_1= a.\phi_1 + c.\phi_2,\; \; \phi'_2= b.\phi_1 + d.\phi_2 .\leqno{(5.4.1.2)}$$  Dualement,    
 $$y_1{\partial \over {\partial y_1}}= \pi_\ast(d\, y'_1{\partial \over  \partial y'_1} - c\, y'_2 {\partial \over  \partial y'_2}),\; y_2{{\partial \over  \partial y_2}}  = \pi_\ast(-b\, y'_1{\partial \over  \partial y'_1} + a\, y'_2{\partial \over  \partial y'_2}).\leqno{(5.4.1.3)}$$  
 
 Par image inverse ${\frak N}$ fournit un module \`a con\-ne\-xion in\-t\'e\-gra\-ble sur ${\frak Y}'$, et par compl\'etion en $(y'_1, y'_2)$, on obtient un module du type consid\'er\'e dans ce paragraphe, au voisinage formel du point de croisement $Q':\,y'_1=y'_2=0$.

  \proclaim Th\'eor\`eme 5.4.1. {\rm (Sabbah).}  Supposons $\car k=0$. Il existe une suite finie $\pi$ d'\'eclate\-ments formels toriques telle que tout point $Q'$ de croisement de $\pi^{-1}(Q)$ soit semi-stable pour l'image inverse de ${\frak N}$. Ceci reste valide si l'on continue \`a effectuer des \'ecla\-te\-ments toriques.\par 

 La preuve de [${\bf S}$, III 4.3.1] utilise la m\'ethode des orbites nilpotentes de [${\bf BaV}$]; en voici une qui \'evite les complications inh\'erentes \`a la m\'ethode de [${\bf BaV}$].

\medskip\noindent  \dem
  On raisonne par r\'ecurrence sur le rang $\nu$ du $k[[y_1,y_2]][{{1}\over {y_1 y_2}}]$-module libre ${\frak N} $ (\cf 5.1.1). Il n'y a rien \`a d\'emontrer pour $\nu \leq 1$. 
 
\proclaim Lemme 5.4.2.  ($\car k=0$ ou bien $ >  \nu$.) Apr\`es ramification mod\'er\'ee de m\^eme degr\'e autour de $y_1=0$ et de $y_2=0$, il existe une suite $\pi'$ d'\'ecla\-te\-ments toriques, un ensemble fini $\Lambda \subset {k}  $,  et,  dans chaque carte torique, une base ${\bf n}(\lambda)$ de l'image inverse de $$ {\frak N}^\Lambda := {\frak N}\otimes_{ {k}[[y_1,y_2]][{{1}\over {y_1 y_2}}]}  {k}[\lambda]_\Lambda[[y_1,y_2]][{{1}\over {y_1 y_2}}]$$  dans laquelle la matrice de 
$$\partial(\lambda) := y_1{\partial \over {\partial y_1}}+ \lambda y_2 {\partial \over {\partial y_2}} $$
n'a pas de p\^ole ou bien s'\'ecrit, dans les co\-or\-don\-n\'ees adapt\'ees $(y'_1,y'_2)$, sous la forme $$(y'_1)^{-r_1}(y'_2)^{-r_2}G(\lambda, y'_1, y'_2)$$ avec $r_1, r_2\geq 0$ non tous deux nuls, $\,G(\lambda, y'_1, y'_2) \in M_\nu({{{k}}}[\lambda]_{\Lambda}[[y'_1,y'_2]])$, et o\`u pour tout $\lambda_0\notin \Lambda$, $G(\lambda_0, 0, 0)\in M_\nu({{{k}}})$ n'est pas nilpotente. 
\par 
 
 \dem Soit $m\in {\frak N}$ un vecteur cyclique pour $ {\frak N}\otimes {\rm Frac} ({k}[[y_1,y_2]]) $ relativement \`a $\partial(0)= y_1{\partial \over {\partial y_1}}$. Dans une base convenable de $$\;\bigwedge^\nu  {\frak N}^\emptyset \cong {k}[\lambda][[y_1,y_2]][{{1}\over {y_1 y_2}}],\;$$   
 on peut donc \'ecrire d'une part  
 $$ m \wedge \partial(0) m \wedge \ldots \wedge \partial(0)^{\nu-1} m =  g(y_1,y_2) \in {k}[[y_1,y_2]],$$ o\`u $g$ n'est pas divisible par $y_1$ ni par $y_2$, et d'autre part  $$ m \wedge \partial(\lambda) m \wedge \ldots \wedge \partial(\lambda)^{\nu-1} m =  y_1^{-s_1}y_2^{-s_2}.g(\lambda,y_1,y_2) \in  {k}[\lambda][[y_1,y_2]], $$ 
o\`u $g(\lambda,y_1,y_2)$ n'est pas divisible par $y_1$ ni par $y_2$  et $  g(0,y_1,y_2)=y_1^{s_1}y_2^{s_2}g( y_1,y_2)\neq 0$.
  Alors pour tout $\lambda_0$ hors d'une partie finie $\Lambda'\subset {k}$,    $${\rm ord}_{y_1 }\,g(\lambda_0,y_1,0)\leq {\rm ord}_{y_1 }\,g(\lambda,y_1,0),\;{\rm ord}_{ y_2 }\,g(\lambda_0,0,y_2)\leq {\rm ord}_{ y_2 }\,g(\lambda,0,y_2).$$
    Ces bornes uniformes permettent de ``chasser" tous les diviseurs $g(\lambda_0,y_1,y_2)=0$ par \'ecla\-te\-ments toriques: il existe une suite $\pi''$  d'\'ecla\-te\-ments toriques telle que dans chaque carte torique, et pour tout $\lambda_0\notin\Lambda'$, le transform\'e strict de $g(\lambda_0,y_1,y_2)=0$ ne rencontre aucun croisement de $(\pi'')^{-1}(Q)$. Ainsi
 $$   {\bf m}(\lambda) = (m,  \partial(\lambda)(m),\ldots, \partial(\lambda)^{\nu-1}(m)) $$ 
fournit une base cyclique de l'image inverse de ${\frak N}^{\Lambda'} $ relativement \`a $ \partial(\lambda)$, et en induit une apr\`es sp\'ecialisation en tout $\lambda_0\notin \Lambda'$.
    
 Il en est encore de m\^eme si l'on effectue au d\'epart une ramification mod\'er\'ee de m\^eme degr\'e $e$ autour de $y_1=0$ et $y_2=0$ (une telle ramification commute aux \'ecla\-te\-ments toriques, donc \`a $\pi''$, et divise $\partial(\lambda)$ par $e$). Cela permet de supposer que les rangs de Poincar\'e-Katz $\rho_1,\rho_2$ le long des axes adapt\'es ${y'}_1=0, {y'}_2=0$ de chaque carte torique sont entiers. Modifions alors $\bf m(\lambda)$ en une base 
 $$  {\bf m}(\lambda)\pmatrix{1&0&&0\cr 0& (y'_1)^{\rho_1}(y'_2)^{\rho_2} &&0\cr &&\ldots &\cr 0&0&&  (y'_1)^{ \rho_1(\nu-1)}(y'_2)^{\rho_2(\nu-1)} } $$ et \'ecrivons la matrice de $\partial(\lambda)$ dans cette base sous la forme $${y'}_1^{-\rho_1}{y'}_2^{-\rho_2}H(\lambda, {y'}_1, {y'}_2)$$ avec $H(\lambda, {y'}_1, {y'}_2) \in M_\nu({{{k}}}[\lambda]_{\Lambda'}[[{y'}_1,{y'}_2]][{{1}\over {y'_1 y'_2}}])$. 
Le point 2) du th\'eor\`eme de Turrittin-Katz 2.1.2, appliqu\'e au corps $y'_1$-adiquement complet  ${{{k}}}(\lambda)(({y'}_2))(({y'}_1))$, montre que $H(\lambda, {y'}_1, {y'}_2)$ est \`a coefficients dans ${{{k}}}(\lambda)(({y'}_2))[[{y'}_1]] $ et que si $\rho_1>0$, $H(\lambda, 0, {y'}_2)$ n'est pas nilpotente.   Idem en \'echangeant $y'_2$ et $y'_1$. 
   En outre si l'un des $\rho_i$ est nul, disons $\rho_1=0$, on peut encore se ramener au cas o\`u $H(\lambda, 0, {y'}_2)$ n'est pas nilpotente en multipliant la base par $y'_1$ (ce qui a pour effet d'ajouter  $(y'_2)^{\rho_2}I_\nu$ \`a $H$). 
   
\noindent  On trouve donc une base ${\bf n}(\lambda)$ dans laquelle la matrice de $\partial(\lambda)$ 
   
 - n'a pas de p\^ole si  $\rho_1=\rho_2=0$,
   
 - s'\'ecrit sinon $${y'}_1^{-\rho_1}{y'}_2^{-\rho_2}H(\lambda, {y'}_1, {y'}_2),\;{\rm{avec}}\;\;H(\lambda, {y'}_1, {y'}_2) \in M_\nu({{{k}}}[\lambda]_{\Lambda'}[[{y'}_1,{y'}_2]]),$$ $H(\lambda, 0, {y'}_2)$ et $H(\lambda, {y'}_1, 0)$ \'etant non nilpotentes.
 Autrement dit 
$$h(\lambda_0, {y'}_1, {y'}_2, t) := \det (t . I_\nu - H(\lambda_0,{y'}_1,{y'}_2))-t^\nu \in {{{k}}}[[{y'}_1,{y'}_2]][t]$$ n'est pas divisible par ${y'}_1$ ni par ${y'}_2$, et il en est de m\^eme si on sp\'ecialise $t$ en $ t_0\in {k}$ convenable, fix\'e. 
  Alors pour tout $\lambda_0$ hors d'une partie finie $\Lambda \supset \Lambda'$,      $${\rm ord}_{y'_1 }\,h(\lambda_0,y'_1,0, t_0)\leq {\rm ord}_{y'_1 }\,h(\lambda,y'_1,0, t_0),$$ $${\rm ord}_{ y'_2 }\,h(\lambda_0,0,y'_2, t_0)\leq {\rm ord}_{ y'_2 }\,h(\lambda,0,y'_2, t_0).$$
    Gr\^ace \`a ces bornes uniformes, il existe une suite d'\'ecla\-te\-ments toriques telle que dans chaque carte torique et pour tout $\lambda_0  \notin \Lambda$, le transform\'e strict de $h(\lambda_0, {y'}_1, {y'}_2, t_0)=0$ ne rencontre aucun croisement du diviseur au-dessus de $Q$. Chacune de ces cartes admet pour co\-or\-don\-n\'ees adapt\'ees $y_1, y_2$ avec $y'_1= (y''_1)^a(y''_2)^b, \, y'_2=  (y''_1)^c(y''_2)^d$ avec $a,b,c,d \geq 0, \, ad-bc=1$. On voit donc que $\bf n(\lambda)$ fournit une base o\`u la matrice de $\partial(\lambda) $ s'\'ecrit $$(y''_1)^{-r_1}(y''_2)^{-r_2}G(\lambda, y''_1, y''_2)$$ avec $G(\lambda, y''_1, y''_2)= H(\lambda, {y'}_1, {y'}_2) \in M_\nu({{{k}}}[\lambda]_{\Lambda}[[y''_1,y''_2]])$, et $$\det (t  . I_\nu - G(\lambda_0,y''_1,y''_2))-t ^\nu  = h(\lambda_0, {y'}_1, {y'}_2, t )\in {k}[[y''_1,y''_2]][t]$$ ne s'annule pas en $y''_1=y''_2=0$. Donc on obtient que pour tout $\lambda_0\notin \Lambda$, $G(\lambda_0, 0, 0)\in M_\nu({{{k}}})$ n'est pas nilpotente.
 \sq
 
 \medskip A partir de ce lemme, la preuve du th\'eor\`eme 
suit celle de [${\bf S}$, III.4.3.1]:  
  en tordant par une con\-ne\-xion de rang un, on peut supposer $\bigwedge^\nu{\frak N}$ r\'egu\-lier. On compl\`ete la suite d'\'ecla\-te\-ments toriques en subdivisant l'\'eventail cor\-res\-pondant de mani\`ere \`a ce que tout c\^one poss\`ede une ar\^ete dont le vecteur primitif $(a,c)$ est tel que l'image de $a$ dans $k$ soit non nulle et que l'image de $c/a$ ne soit pas dans $\Lambda$ ({\it c'est uniquement l\`a que l'on utilise l'hypoth\`ese  $\car k=0$}), et pour laquelle le diviseur associ\'e est une composante du diviseur polaire de $\partial(c/a)$ si ce dernier est non vide.  
 
 Soit alors une carte torique, de co\-or\-don\-n\'ees adapt\'ees $y'_1, y'_2$. Supposons que $y'_1$ cor\-res\-ponde \`a une ar\^ete dont la pente $\lambda_0 = c/a$ n'est pas dans $ \Lambda$. Dans la base ${\bf n}(c/a)$ de l'image inverse de ${\frak N}$, la matrice de $y'_1{\partial\over \partial y'_1} = a .\partial(c/a)  $  est donc de la forme $$(y'_1)^{-r_1}(y'_2)^{-r_2}G(y'_1, y'_2)$$ avec $G(y'_1, y'_2) \in M_\nu({{{k}}}[[y'_1,y'_2]])$, et  $G(0, 0)$ n'est pas nilpotente d\`es lors que $r_1$ et $r_2$ ne sont pas tous deux nuls. Comme $\bigwedge^\nu{\frak N}$ est r\'egu\-lier,  $G(0, 0)$ est en outre de trace nulle, donc elle a deux valeurs propres distinctes. Si $r_1$ et $r_2$ ne sont pas tous deux nuls, le lemme de d\'e\-com\-po\-si\-tion 2.2.1 (points 1) et 3))  s'applique aux d\'erivations $\delta = (y'_1)^{ r_1+1}(y'_2)^{ r_2 }{\partial\over \partial y'_1} $ et $\delta' = (y'_1)^{ s_1}(y'_2)^{ s_2 +1}{\partial\over \partial y'_2} $ de $k[[y'_1,y'_2]]$ (pour $s_1,s_2\geq 0$ convenables) et au $k[[y'_1,y'_2]]$-r\'eseau engendr\'e par ${\bf n}(c/a)$ (on pourrait, alternativement, invoquer 2.5.1, 1) et 2.2.2); d'o\`u une d\'e\-com\-po\-si\-tion de la con\-ne\-xion  for\-melle  (in\-t\'e\-gra\-ble) ${\frak N}$, ce qui permet de dimi\-nuer $\nu$. 
 
 Si au contraire $r_1=r_2=0$, on remarque que la con\-ne\-xion est r\'egu\-li\`ere le long de $y'_1=0$ et $y'_2=0$: en effet, elle l'est non seulement pour la d\'erivation $y'_1{\partial\over \partial y'_1} $ mais aussi pour toute d\'erivation du type $\partial(\lambda_1) $ avec $\lambda_1$ hors de $\Lambda$ et distinct de l'image de $c/a $ (on peut aussi invoquer 2.5.2).   
 
 \medskip Enfin, pour montrer qu'on peut choisir la ramification en $y_1$ et $y_2$ de degr\'e $e$ divisant $\nu!$, on utilise le m\^eme argument galoisien que dans 2.3.1, 3).
 \sq

 \bigskip\noindent {\bf 5.5. Points de croisement stables.} Reprenons les notations de 5.3.  Disons que le point de croisement $Q$ est {\it stable} pour ${\frak N}$ s'il est semi-stable, et si les $\omega_j$ peuvent \^etre choisis de telle sorte que pour $i=1,2$, les \'el\'ements $\phi_{i, j, r_{ij}}$ et leurs diff\'erences $\phi_{i, j, r_{ij}}- \phi_{i, j', r_{ij'}} $ ne s'annulent pas en $Q$ (lorsqu'ils sont non nuls). 
 
 Si $\car k=0$, il revient au m\^eme, en vertu de (5.1.2.3), de demander que pour $\omega= \omega_j$ et $\omega = \omega_j-\omega_{j'}$, l'\'el\'ement  $\alpha \in k[[y_1,y_2]]$  associ\'e  (\cf 5.1.2) ne s'annule  pas en $Q$ s'il est non nul (on peut toujours supposer, modulo $\Omega^1_{\frak Y}(\log {\frak Z})^{d=0}$, que les r\'esidus de $\omega$ sont nuls). 
 
    \smallskip\noindent {\bf Remarque 5.5.1.}  Si $Q$ est un point de croisement stable, les in\'egalit\'es de 5.3.2 sont des \'egalit\'es.

  \medskip Cette notion de point de croisement stable  \'equivaut  \`a la notion de ``bonne structure for\-melle" de [$\bf S2$] (modulo l'argument galoisien indiqu\'e \`a la fin de la preuve de 5.4.1). Dans la suite, elle ne nous servira pas, et nous ne la mentionnons que pour faire le lien pr\'ecis avec la probl\'ematique de [$\bf S2$] d'une part (\cf 8.1 ci-dessous), et par sym\'etrie avec la d\'efinition 3.4.2 d'autre part, sym\'etrie pr\'ecis\'ee par le r\'esultat suivant.  
   
   \proclaim Proposition 5.5.2. $i)$ Soit ${\frak Y}'\to \frak Y=\Spf k[[y_1,y_2]]$ l'\'eclatement du point $Q:\, y_1=y_2=0$, et soit $\frak X = \Spec k[x_2, {1\over x_2}][[x]] $ la carte for\-melle affine de  ${\frak Y}'$ donn\'ee par les co\-or\-don\-n\'ees  $\,x=y_1,\, x_2=y_2/y_1.$  Soit $\frak N$ une con\-ne\-xion m\'eromorphe for\-melle sur $\frak Y$ \`a p\^oles le long de $y_1y_2=0$, et soit ${\frak N}'$ son image inverse sur $\frak X$ (qui est une con\-ne\-xion m\'eromorphe for\-melle  \`a p\^oles le long de $Z= {\frak X}_{red} = \Spec k[x_2, {1\over x_2}]$). 
   \smallskip Si 
   $Q$ est (semi-)stable pour $\frak N$, alors tout point $P$ de $Z$ est (semi-)stable pour ${\frak N}'$.
   \medskip\noindent $ii)$ Soit $Z= \Spec A$ une courbe affine lisse connexe munie d'une co\-or\-don\-n\'ee \'etale $x_2$ et soit  ${\frak X}'\to \frak X = \Spf A[[x]]$ l'\'eclatement du point $P:\, x=x_2=0$. Soit ${\frak Y} $ le compl\'et\'e de ${\frak X}'$ en le point d'intersection $Q$ du diviseur exceptionnel et du transform\'e strict de $Z = {\frak X}_{red}$. \'Ecrivons ${\frak Y}  = \Spf k[[y_1,y_2]]$ o\`u 
   $\,y_1=x/x_2,\, y_2=x_2.$  Soit $\frak M$ une con\-ne\-xion m\'eromorphe for\-melle sur $\frak X$ \`a p\^oles le long de $Z$, et soit ${\frak M}'$ son image inverse sur $\frak Y$ (qui est une con\-ne\-xion m\'eromorphe for\-melle \`a p\^oles le long de $y_1y_2=0$). 
  \smallskip  Si $P$ est (semi-)stable pour $\frak M$, alors $Q$ est (semi-)stable pour ${\frak M}'$. \par 

 La d\'emonstration, directe, est laiss\'ee au lecteur. Prendre garde que dans $i)$ et $ii)$, les r\'eciproques ne sont pas vraies (\cf contre-exemple 5.2.2).

 \medskip  Par ailleurs, il n'est pas difficile de compl\'eter 5.4.1 par:
  
     \proclaim Scholie 5.5.3. Supposons que  $Q$ soit un point de croisement semi-stable. Alors il existe une suite finie $\pi$ d'\'eclatements formels toriques telle que tout point de croisement $Q'$ de $\pi^{-1}(Q)$ soit stable pour l'image inverse de $\frak N$.\par 

Voir [$\bf S2$] III.1.3).

 \bigskip \bigskip   
 \vfill\eject
  {\bf 6. Semi-continuit\'e du rang de Poincar\'e-Katz.}

\bigskip\noindent  {\bf 6.1. \'Enonc\'e.} Soit $k$ un corps alg\'ebriquement clos de carac\-t\'e\-ris\-tique nulle. Soient ${\frak Y}$ le sch\'ema  formel $(y_1,y_2)$-adique $ \Spf  k[[y_1,y_2]]$, et ${\frak Z }= {\frak Z}_i\cup\ldots \cup  {\frak Z}_t$ une r\'eunion de germes de courbes sur $\frak Y$. Si $\frak C$ est un autre germe de courbe, la {\it multiplicit\'e d'intersection} $ ({\frak C},{\frak Z}_i)_Q$  est la dimension sur $k$ du quotient de $k[[y_1,y_2]]$ par la somme de l'id\'eal de $\frak C$ et de l'id\'eal de ${\frak Z}_i$. On note $\tilde{\frak C}\to \frak C$ la normalisation, $O$ le point $\tilde{\frak C}_{red}$.

On se donne une {\it con\-ne\-xion m\'ero\-morphe  for\-melle} ${\frak N}$ sur $\frak Y$ {\it \`a p\^oles le long de $\frak Z$}, c'est-\`a-dire un $\sO_{\frak Y}(\ast {\frak Z})$-module  projectif de type fini muni d'une con\-ne\-xion {\it in\-t\'e\-gra\-ble} relative \`a $k$.

Le r\'esultat suivant compare le rang de Poincar\'e-Katz de l'image inverse ${\frak N}_{\tilde{\frak C}}$ de ${\frak N}$ sur $\tilde\frak C$ aux rangs de Poincar\'e-Katz de ${\frak N}$ le long des ${\frak Z}_i$.

\proclaim Th\'eor\`eme 6.1.1. On a  $$ \rho_O({\frak N}_{\tilde {\frak C} }) \leq \sum_i\, ({\frak C},{\frak Z}_i)_Q\,. \,\rho_{{\frak Z}_i}({\frak N}  ).  \leqno (6.1.1.1) $$ 
 \par   
 
 Le cas par\-ti\-cu\-lier de (6.1.1.1) o\`u le second membre est nul donne:
 
   \proclaim Corollaire 6.1.2. {\rm (Deligne [${\bf De1}$]).} Si ${\frak N}$ est r\'eguli\`ere le long des ${\frak Z}_i$, alors ${\frak N}_{\tilde {\frak C} }$ est r\'eguli\`ere en $O$. \par   
   
\noindent La preuve que nous en donnons ici est apparemment la premi\`ere preuve purement alg\'ebrique de ce r\'esultat; pour plus de d\'etails sur ce sujet, voir [$\bf{B}$] et [$\bf A3$].
 
   \proclaim Corollaire 6.1.3.  Soient $X$ et $Y$ des vari\'et\'es alg\'ebriques lisses sur $k$ (\resp des vari\'et\'es analytiques complexes). Soit $f: Y\to X$ un morphisme lisse, de dimension relative $1$, \`a fibres connexes,  et soit $Z$ une hypersurface de $Y$ finie sur $X$ (via $f$). Soit ${\sN}$ un module \`a  con\-ne\-xion m\'ero\-morphe {\rm in\-t\'e\-gra\-ble} sur $Y$ \`a p\^oles le long de $Z$. 
  Alors la somme des rangs de Poincar\'e-Katz $$\; \displaystyle\sum_{z\in Z, \,f(z)=x} \, {  \rho}_z({\sN}_{(x)})\;$$ 
est une fonction semi-continue inf\'erieurement de $x\in X$ .
\par
 \dem (de 6.1.1 $\Rightarrow$ 6.1.3) Soit $\partial$ une d\'erivation le long des fibres de $f$. En lisant le rang de Poincar\'e-Katz sur l'\'equation diff\'erentielle donn\'ee par un vecteur cyclique local (relatif \`a ${\sN}$ et $\partial$), on voit facilement que la fonction dont il est question est constructible pour la topologie naturelle sur $X$. Pour \'etablir sa semi-continuit\'e, on peut donc se limiter aux (germes de) courbes lisses trac\'ees sur $X$ passant par un point fix\'e arbitraire $y_0$. Il suffit d'appliquer 6.1.1 en prenant pour  $\frak C$ le germe formel de $f^{-1}(x_0)$ au voisinage de chacun des points $Q$ de $Z\cap f^{-1}(x_0)$, et pour $\frak Z_i$ les germes des branches de $Z$ passant par $Q$. \sq  
 
\medskip Le corollaire vaudrait tout aussi bien pour des vari\'et\'es analytiques r\'eelles, voire des espaces $k$-analytiques lisses de Berkovich sur un corps $k$ complet non archim\'edien.

  \bigskip\noindent  {\bf 6.2. Quelques propri\'et\'es des diviseurs $D_{Z,\sigma}({\frak M})$.}  
   Soit $R$ un anneau local noeth\'erien complet, d'id\'eal maximal $\frak m$ et de corps r\'esiduel $k$. Soit $ {\frak Y}'  $ un $R$-sch\'ema formel $\frak m$-adique quasi-compact plat topologiquement de pr\'esentation finie.
   Soient ${\frak Z}'_i$ des germes de courbes lisses sur ${{ \frak Y}'}$, et posons ${\frak Z}'= (\cup {\frak Z}'_i) \cup  {\frak Y}'_{red}$. 
 
 On suppose que les composantes irr\'eductibles de ${\frak Y}'_{red}$ sont des courbes lisses et que les composantes de ${\frak Z}'$ se coupent transversalement.
 
 Soit $\bar Z$ l'une des composantes irr\'eductibles de ${\frak Y}'_{red}$.  On suppose que $\bar Z$ est r\'eguli\`erement immerg\'e dans ${\frak Y}'$, de sorte qu'on peut parler du fibr\'e normal  $ N_{\bar Z }\,{{ \frak Y}'}  $. On note $P(N_{\bar Z }\,{{ \frak Y}'} )$ son compl\'et\'e projectif 
 (adjonction d'une section \`a l'infini au-dessus de $\bar Z$).

 Soit $Z= \Spec A$ un ouvert affine dense de $\bar Z$ ne rencontrant pas les croisements de ${\frak Z}'$. On suppose que ${{ \frak Y}'}$ admet un sous-sch\'ema formel ouvert ${\frak X} \cong \Spf A[[x]]$ (la topologie $\frak m$-adique sur $\sO_{{ \frak Y}'}$ induisant la topologie $x$-adique sur $\sO_{\frak X}$) avec $Z =  {\frak X}_{red}$.

  \medskip\noindent {\bf Exemple 6.2.1.} Repla\c cons-nous provisoirement dans la si\-tua\-tion de 6.1.1, et consid\'erons une suite $\pi: {{\frak Y}'}\to \frak Y  $ d'\'eclatements ponctuels formels telle que les transform\'es stricts ${\frak Z}'_i$ des ${\frak Z}_i$ coupent ${{\frak Y}'}_{red}$ transversalement. On  prend pour  $\bar Z$ l'une des composantes irr\'eductibles du diviseur exceptionnel $\pi^{-1}(Q)$, et pour $Z$ un ouvert dense de $\bar Z$ ne rencontrant pas les croisements. C'est l'exemple qui nous servira en 6.3 et en 7.3.

 \medskip  Soit alors $\frak N'$ une con\-ne\-xion m\'ero\-morphe for\-melle sur ${\frak Y}'$ \`a p\^oles le long de ${\frak Z}'$. Sa restriction $\frak M = {\frak N }'_{\mid \frak X}$  au sous-sch\'ema formel affine $\frak X$ est du type consid\'er\'e au \S 4. 
 
  \smallskip On suppose que  {\it tous les points de $Z$ sont semi-stables pour $\frak M$}. On dispose alors, pour chaque pente non nulle $\sigma$ de $\frak M$ (le long de $Z$), du diviseur positif $\, D_{Z,\sigma}({\frak M})\,$ sur le fibr\'e normal $\,N_{Z}\,\frak X \cong \Spec A[x]\,$ (\cf 3.6, 4.2.2). Rappelons que $\, D_{Z,\sigma}({\frak M})= (\varphi_\sigma(x)) \,$ est donn\'e par l'\'equation 
 $$ \varphi_\sigma(x):=  \prod_{j\in J_{(\sigma)}} \, (x^\sigma  -\phi_{j, -\sigma})^{\mu_j}=0, \;\;\leqno{(6.2.1.1)}$$ 
 $\phi_{j, -\sigma}$ \'etant le coefficient, dans une extension finie de $A$, du terme de plus bas degr\'e ($= -\sigma$) en $x$ dans $\phi_j$.  Ce diviseur $D_{Z,\sigma}({\frak M})$ est fini sur $Z$, de degr\'e sur $ Z$ \'egal \`a l'irr\'egularit\'e $\sigma.\mu_{(\sigma)}$ de $\frak M_{(\sigma)}$ (o\`u $\mu_{(\sigma)}= \sum_{j\in J_{(\sigma)}}\, \mu_j$ est la dimension de la composante $\frak M_{(\sigma)}$ de $\frak M$ de pente $\sigma$).

 \smallskip En prenant l'adh\'erence de Zariski dans $P(N_{\bar Z }\,{{ \frak Y}'} )$, on obtient un diviseur de Weil de $P(N_{\bar Z }\,{{ \frak Y}'} )$, not\'e $\bar D_{\bar Z,\sigma}({\frak N }')$, dont le lemme suivant pr\'ecise l'intersection avec la section \`a l'infini $(\infty)$.

  \proclaim Lemme 6.2.3.   $i)$ $\bar D_{\bar Z,\sigma}({\frak N }')$ ne rencontre $(\infty)$ qu'au dessus de $\bar Z\setminus Z$. 
 \medskip \noindent $ii)$  Si $\sigma = \rho_{\bar Z}({\frak N}')$, $\bar D_{\bar Z,\sigma}({\frak N }')$ ne rencontre $(\infty)$ qu'au dessus des points de croisement de $\bar Z$ et ${\frak Z}'$. 
 \medskip \noindent $iii)$   Au-dessus d'un point de croisement semi-stable $Q'$ de $ \bar Z $ et d'une autre composante $\bar T$ de ${{ \frak Y}'}_{red}$ (\resp de $ \bar Z $ et de ${{\frak Z}'_i}$), la multiplicit\'e d'inter\-sec\-tion $$\,(\bar D_{\bar Z,\sigma}({\frak N }') \,,  \,(\infty) )\,$$ est major\'ee par $\,\mu_{(\sigma)}.\rho_{{\bar T}}({\frak N}')\,$ (\resp $\,\mu_{(\sigma)}.\rho_{{{\frak Z}'_i}}({\frak N}')$).  
 \medskip \noindent $iv)$ Au-dessus d'un tel point de croisement, la multiplicit\'e d'intersection $$((\sum_{\sigma >0}\, \bar D_{\bar Z,\sigma}({\frak N }')) \,, \,(\infty) ) \,$$  est major\'ee par ${\rm ir}_{{\bar T}}({\frak N}') \,$ (\resp $ \,{\rm ir}_{{{\frak Z}'_i}}({\frak N}')$). \par   

 \dem  $i)$ est clair sur l'\'equation (6.2.1.1) de $D_{Z,\sigma}({\frak M})$.
 
\smallskip\noindent $ii)$ Soit $P$ un point, \'eventuellement non semi-stable, de $\bar Z\setminus Z$ qui ne soit pas au croisement avec $\frak Z$.   Prenons une co\-or\-don\-n\'ee locale $x_2$ sur $\bar Z$ qui s'annule \`a l'ordre $1$ en $P$ (de sorte que le compl\'et\'e de ${{ \frak Y}'}$ en $P$ s'identifie \`a $\Spf k[[x ,x_2]]$).  On lit sur l'\'equation (6.2.1.1) (quitte \`a ramifier pour se ramener au cas plus net o\`u $\sigma$ est entier) que  $$  (\bar D_{\bar Z,\sigma}({\frak N }'), (\infty) )\,   =     \,\max ( 0, -  \sum_{j\in J_{(\sigma)}} \, \mu_j. \ord_{x_2} \, \phi_{ j, -\sigma}).  \leqno{(6.2.3.1)}$$ Or, d'apr\`es 3.3.1,  $ \ord_{x_2} \, \phi_{ j, -\sigma} \geq 0$ lorsque $\sigma = \rho_{\bar Z}({\frak N}')$. 
 
\smallskip\noindent $iii)$ et $iv)$. Identifions le compl\'et\'e de ${{ \frak Y}'}$ en $Q'$ \`a $\Spf k[[y_1 ,y_2]]$ o\`u $y_1= x$ et o\`u $y_2=0 $ d\'efinit $\bar  T$ (\resp ${{\frak Z}'_i}$). Reprenons les notations de 5.3. Notons $\nu_h$ le rang du facteur ${\frak N}'_{\bar\omega_h}$, et $H_{(\sigma)}$ l'ensemble des indices $h$ pour lesquels $\max(0, -\ord_{y_1}\,\varpi_1(\psi_{1, h}))=\sigma$, de sorte que $\mu_{(\sigma)} = \sum_{h\in H_{(\sigma)}}\, \nu_h$. Pour $h\in H_{(\sigma)}$, notons $\psi_{1,h, -\sigma}$ le terme de plus bas degr\'e en $y_1$ de $\psi_{1, h}$, qui est aussi celui de $\varpi_1(\psi_{1, h})$ puisque $\sigma>0$. 

\noindent En vertu de la d\'ecomposition (5.3.1.5)$_2$, on a   
$$ \max(0, - \ord_{y_2} \, \psi_{ 2,{h}}) =  \max(0, - \ord_{y_2} \, \varpi_2(\psi_{ 2,{h}})) \leq   \rho_{{\bar T}}({\frak N}'),\;\;  resp \; \leq \, \rho_{{{\frak Z}'_i}}({\frak N}') ,\leqno{(6.2.3.2)}$$ 
$$\sum_{h}\, \nu_{h}.\max(0, - \ord_{y_2} \, \psi_{ 2,{h}}) \leq   {\rm ir}_{{\bar T}}({\frak N}') ,\;\;  resp \;\leq \, {\rm ir}_{{{\frak Z}'_i}}({\frak N}').\leqno{(6.2.3.3)}$$
  En vertu de la d\'ecomposition (5.3.1.5)$_1$, on peut r\'e\'ecrire  $\varphi_\sigma(x)$ sous la forme
$$ \prod_{{h}\in H_{(\sigma)}}\,  \, (y_1^\sigma  -\psi_{1,{h},-\sigma})^{\nu_{h}}=0 , \;\;\leqno{(6.2.3.4)}$$ 
et (6.2.3.1) sous la forme $$ (\bar D_{\bar Z,\sigma}({\frak N }'), (\infty) )\,=   \,\max ( 0, -  \sum_{{h}\in H_{(\sigma)}} \, \nu_j. \ord_{y_2} \, \psi_{ 1,{h}, -\sigma}).  \leqno{(6.2.3.5)}$$
Comparant (6.2.3.2), (6.2.3.3) et (6.2.3.5), on voit qu'il suffit de d\'emontrer que 
$$\max(0,- \ord_{y_2} \, \psi_{ 1,{h}, -\sigma}) \leq  \max(0, - \ord_{y_2} \, \psi_{ 2,{h}}).\leqno{(6.2.3.6)}$$
Or $ \ord_{y_2} \, \psi_{ 1,{h}, -\sigma}  \geq  \ord_{y_2} \, \psi_{ 1,{h} } ,$ et par 
  int\'egrabilit\'e (5.1.2.1),  $\max(0, -\ord_{y_2}\,  \psi_{1,{h}} )=$ $$=  \max(0, -\ord_{y_2}\,  \partial_{ 2} \psi_{1,{h}} ) = \max(0, -\ord_{y_2}\,  \partial_{ 1} \psi_{2,{h}} ) \leq  \max(0, - \ord_{y_2}\,   \psi_{2,{h}} ) .\;\;\square
$$     
     
   \proclaim Proposition 6.2.4. Supposons $\bar Z$ propre, et que tous les points de croisement de $\bar Z$ soient semi-stables pour ${\frak N}'$. Alors on a l'in\'egalit\'e 
   $$  \displaystyle{ (-deg(N_{\bar Z }\,{  \frak Y}' ))\,.\, \rho_{\bar Z}({\frak N}')
    \leq \sum_{{\bar T}\cap {\bar Z}\neq \emptyset} \, \rho_{ \bar T }({\frak N}') +  \sum_{{\frak Z}'_i\cap{ \bar Z}\neq \emptyset} \, \rho_{{\frak Z}'_i}({\frak N}'). } \leqno{(6.2.4.1)}$$   En outre, si tous les points de $\bar Z$ sont semi-stables, on a l'in\'egalit\'e 
     $$ \displaystyle{(-deg(N_{\bar Z }\,{{ \frak Y}'}))\,.\, {\rm ir}_{\bar Z}({\frak N}')  \leq \sum_{{\bar T}\cap {\bar Z}\neq \emptyset} \,  {\rm ir}_{ \bar T }({\frak N}') +  \sum_{{\frak Z}'_i\cap{ \bar Z}\neq \emptyset} \,  {\rm ir}_{{\frak Z}'_i}({\frak N}'). } \leqno{(6.2.4.2)}$$
     \par 
   
    \dem Compte tenu de 6.2.3, cela d\'ecoule de la formule 
    $$ -  deg(N_{\bar Z }\,{{ \frak Y}'}).\delta\, = \, \bar C.(\infty) - \bar C.(0)$$ valable pour toute courbe ferm\'ee $\bar C\subset P(N_{\bar Z }\,{{ \frak Y}'})$ dont la projection sur $\bar Z$ est finie de degr\'e $\delta$  (formule qui se d\'emontre en prenant l'image inverse par $\tilde  C\to \bar Z$, o\`u $\tilde  C$ d\'esigne la normalis\'ee de $ \bar C$, et en utilisant la section naturelle de $\,N_{\bar Z }\,{{ \frak Y}'} \times_{\bar Z} \tilde  C\to \tilde  C$). On applique cette formule aux composantes $\bar C$ de $\bar D_{\bar Z,\sigma}({\frak N }')$.  \sq

 \bigskip\noindent  {\bf 6.3. Preuve de 6.1.1. } Comme dans l'exemple 6.2.1, effectuons une suite d'\'eclatements ponctuels formels 
$\,\pi : {\frak Y}'\to {\frak Y}\,$ (en commen\c cant par \'eclater $Q$), de sorte que 
 
\smallskip  \noindent $i)$ $\pi^{-1}({\frak Z}\cup {\frak C})$ soit \`a croisements normaux, 

\smallskip \noindent $ii)$ chaque composante irr\'eductible du diviseur exceptionnel ne coupe qu'un au plus des   $ {\frak Z}_i'$ ($=$ transform\'e  strict  de ${\frak Z}_i$),

\smallskip \noindent $iii)$ tout point de croisement du diviseur exceptionnel soit semi-stable pour ${\frak N}' = \pi^\ast({\frak N})$  (c'est loisible en vertu de 5.4.1). 

 \medskip Le diviseur exceptionnel est arborescent. On note ${\Bbb T}$ l'arbre dual: ses sommets $v$ sont en bijection avec les composantes $E_v$ de $\pi^{-1}Q$.  Vu $ii)$, il est loisible d'indexer  $ {\frak Z}'_i$  par le sommet $v$ cor\-res\-pondant  \`a la composante $E_v$ qu'il coupe (pour un sous-ensemble fini $V$ de sommets). 
On note $v_0$ le sommet cor\-res\-pondant  \`a la composante $E_{v_0}$ que coupe ${\frak C}'$.

 On note enfin $A$  la matrice (sym\'etrique) d'intersection des $E_v,\; v\in {\Bbb T}$, affect\'ee du signe $-$; ses coefficients diagonaux $A_{vv}$ sont les entiers naturels $- deg(N_{E_v}{\frak Y}') $, et ses coefficients non diagonaux $A_{vw}$ valent $-1$ si $v$ et $w$ sont voisins (c'est \`a dire si $E_v$ et $E_w$ se coupent) et $0$ sinon. 

En appliquant (6.2.4.1) \`a $\bar S = E_v$, on obtient 
$$ \sum_{w\in {\Bbb T}}\;  A_{vw}\,. \,{\rho}_{E_{w}}({\frak N}' ) \leq   {\rho}_{{\frak Z}'_v}({\frak N}' ) $$
si $v\in V$, et 
$$ \sum_{w\in {\Bbb T}}\;  A_{vw}\,. \,{\rho}_{E_{w}}({\frak N}' ) \leq 0$$ sinon.  D'o\`u
$$ {\rho}_{E_u}({\frak N}' ) \leq \sum_{v\in V} \,  ( A^{-1})_{uv}\,. \,{\rho}_{{\frak Z}'_v}({\frak N}' ).\leqno{(6.3.1.1)}$$
Si on se donne, pour chaque $v$, une paire de germes de courbes ${\frak C}'_v, {\frak C}''_v$ coupant $E_v$ transversalement hors des croisements, il est classique que $$( A^{-1})_{vw}=(\pi({\frak C}'_v).\pi({\frak C}''_w))_{Q }$$ (la d\'emonstration, par r\'ecurrence sur le nombre d'\'ecla\-te\-ments, est rappel\'ee dans [${\bf S2}$, I.3.2.8]). On peut prendre en par\-ti\-cu\-lier ${\frak C}'_{v_0}= {\frak C}'$ et ${\frak C}''_v= {\frak Z}'_v$ si $v\in V$, d'o\`u  
$$ {\rho}_{E_{v_0}}({\frak N}' ) \leq \sum_{v\in V}\, (\pi({\frak C} ).\pi({\frak Z}_v))_{Q }\,. \, {\rho}_{{\frak Z}'_v}({\frak N}' ).\leqno{(6.3.1.2)}$$
Par ailleurs, par 3.1.1, on a 
$$ {\rho}_{E_{v_0}}({\frak N}' ) \geq   {\rho}_{Q'}({\frak N}'_{\mid {\frak C}' }) =  {\rho}_O({\frak N}_{\tilde{\frak C} }) .\leqno{(6.3.1.3)}$$
Ceci prouve (6.1.1.1).\sq
  
\medskip\noindent {\bf Remarque 6.3.1.} On obtient une variante - moins \'economique - de cette preuve en rempla\c cant l'usage de 3.1.1 (\`a la fin) par l'usage de 4.3.1, qui permet de choisir $\pi$ de telle sorte que $Q'$ soit un point stable de $E_{v_0}$ (auquel cas (6.3.1.3) est une \'egalit\'e).

\bigskip \bigskip     
   
\bigskip {\bf 7. Semi-continuit\'e de l'irr\'egularit\'e.}

\bigskip\noindent  {\bf 7.1. \'Enonc\'e.} Soit $(Y, Q)$ un germe de surface analytique complexe. Soient $C, Z_1,\ldots, Z_t$ des germes de courbes analytiques dans $(Y, Q)$ (passant par $Q$), et soit $Z=\cup Z_i$.  Soit $\sN$  une {\it con\-ne\-xion m\'ero\-morphe int\'egrable} sur $ Y$ {\it \`a p\^oles le long de $  Z$} (c'est-\`a-dire un $\sO_{ Y}(\ast {  Z})$-module projectif de type fini muni d'une con\-ne\-xion {\it in\-t\'e\-gra\-ble}). On note $(\tilde C, O)\to (C,Q)$ la normalisation.

  Le r\'esultat suivant compare l'irr\'egularit\'e de l'image inverse $\sN_{\tilde C}$ de $\sN$ sur $\tilde C$ aux irr\'egularit\'es de $\sN$ le long des ${ Z}_i$.
    
\proclaim Th\'eor\`eme 7.1.1. On a  $$  {\rm ir}_O(\sN_{\tilde {  C} }) \leq \sum_i\, ({ C},{  Z}_i)_Q\,. \, {\rm ir}_{{  Z}_i}(\sN ). \leqno (7.1.1.1) $$\par   

 Le th\'eor\`eme vaut aussi bien dans le contexte formel plus g\'en\'eral de 6.1.1. Nous avons pr\'ef\'er\'e cette fois le cadre analytique ou la preuve s'\'ecrit de mani\`ere plus lisible. 
 
  \proclaim Corollaire 7.1.2.  {\rm (Conjecture de Malgrange).}  Soit $f: Y\to X$ un morphisme lisse de vari\'et\'es analytiques complexes, de dimension relative $1$, \`a fibres connexes, et soit $Z$ une hypersurface de $Y$ finie sur $X$ (via $f$). Soit $\sN$ un module \`a  con\-ne\-xion m\'ero\-morphe {\rm in\-t\'e\-gra\-ble} sur $Y$ \`a p\^oles le long de $Z$. 
 \smallskip Alors  la
somme des irr\'egularit\'es $$\; \displaystyle\sum_{z\in Z, \,f(z)=x} \, {\rm ir}_z(\sN_{(x)})\;$$ est une fonction semi-continue inf\'erieurement de $x\in X$ .
\par
 La d\'emonstration de 7.1.1 $\Rightarrow$ 7.1.2 est compl\`etement analogue \`a celle de 6.1.1 $\Rightarrow$ 6.1.3.
 
 \medskip  Le corollaire vaudrait tout aussi bien pour des $k$-vari\'et\'es alg\'ebriques lisses, des vari\'et\'es analytiques r\'eelles, voire des espaces $k$-analytiques lisses de Berkovich sur un corps $k$ complet non archim\'edien.

  \bigskip\noindent  {\bf 7.2. Un cas par\-ti\-cu\-lier.} C'est celui, que nous allons traiter directement, o\`u { $t= 1$ ou $2$, ${ Z}_i$ est d\'efini par l'\'equation $y_i=0, (i\leq t)$, et o\`u $ C$ est un germe de courbe lisse coupant ${ Z}_i$ transversalement}.   
  
     Dans cette si\-tua\-tion, prouvons que  
$$  {\rm ir}_Q(\sN_{\mid { C} }) \leq \sum \,{\rm ir}_{{ Z}_i}(\sN  ). \leqno {(7.2.1.1)} $$   
Consid\'erons pour cela la famille $C_\lambda$ de (germes analytiques de) courbes translat\'ees de $C$ le long de l'axe $y_1$. Pour $\lambda\neq 0$ proche de $0$ et assez g\'en\'eral, $C_\lambda$ coupe $Z $ transversalement en $t$ points stables $Q_i, (i\leq t)$. On a donc $$\sum \, {\rm ir}_{Q_i}(\sN_{\mid { C_\lambda} })= \sum \,{\rm ir}_{{ Z}_i}(\sN  ).$$
D'apr\`es un r\'esultat de semi-continuit\'e de Deligne [$\bf De2$]\footnote{$^{(11)}$}{\sm aussi d\'emontr\'e, ult\'erieurement et par voie transcendante, par Mebkhout [$\bf Me1$]. La preuve de Deligne est de nature alg\'ebrico-formelle, et s'appuie sur la caract\'erisation de l'irr\'egularit\'e donn\'ee dans [$\bf GL$].}, on a 
$$ {\rm ir}_Q(\sN_{\mid { C} })+\nu  \leq \sum \, ({\rm ir}_{Q_i}(\sN_{\mid { C_\lambda} })+ \nu) ,  \leqno {(7.2.1.2)}$$ (o\`u $\nu$ d\'esigne le rang de $\sN$), d'o\`u 
$$  {\rm ir}_Q(\sN_{\mid { C} }) \leq \sum \, {\rm ir}_{{ Z}_i}(\sN  ) +(t-1)\nu .  \leqno {(7.2.1.3)_t}  $$ D'o\`u le r\'esultat si $t=1$ (observons incidemment que ce r\'esultat est moins \'el\'ementaire que son analogue 3.1.1 pour le rang de Poincar\'e-Katz). 
\par\noindent Si $t=2$, 
on d\'eduit de (7.2.1.3)$_2$ que, si $E$ d\'esigne le diviseur exceptionnel de l'\'eclat\'e de $Q$, la con\-ne\-xion image inverse $\sN'$ par l'\'eclatement v\'erifie
 $${\rm ir}_{E }{\sN'  }\leq {\rm ir}_{Z_1}{\sN}+ {\rm ir}_{Z_2}{\sN }+  \nu.  $$
En prenant l'image inverse de $\sN$ par  $y_1\mapsto y_1^n,\, y_2\mapsto y_2^n$ (changement de variables qui commute \`a l'\'ecla\-te\-ment), les irr\'egularit\'es sont multipli\'ees par $ n$, d'o\`u l'on d\'eduit
$${\rm ir}_{E }{\sN'  }\leq {\rm ir}_{Z_1}{\sN}+ {\rm ir}_{Z_2}{\sN }+ \nu/n, $$
et, en passant \`a la limite,
 $${\rm ir}_{E }{\sN'  }\leq {\rm ir}_{Z_1}{\sN}+ {\rm ir}_{Z_2}{\sN }. \leqno {(7.2.1.4)} $$
En appliquant (7.2.1.1) (pour $t=1$) au transform\'e strict de $C$ et \`
 a $E$, on en d\'eduit (7.2.1.1) (pour $t=2$).\sq 

\medskip\noindent {\bf Remarque 7.2.1.} Il est tentant d'essayer de prouver 7.1.1 par la m\'ethode de 6.5, en rempla\c cant le recours \`a (6.2.4.1) par le recours \`a (6.2.4.2). Pour appliquer directement (6.2.4.2), qui requiert que tous les points du diviseur exceptionnel, y compris les croisements, soient semi-stables, on pourrait remplacer $\Bbb T$ par l'arbre infini cor\-res\-pondant \`a l'it\'eration altern\'ee de l'\'evitement des points tournants par \'eclatement (4.3.3) et de la semi-stabilisation des croisements par \'eclatement (5.4.1). Ce faisant, on se heurte \`a l'obstacle suivant: tirer l'analogue de (6.3.1.1) cor\-res\-pondant \`a cet arbre infini reviendrait \`a intervertir les sommations (index\'ees par les sommets de l'arbre) dans une s\'erie double non commutativement convergente. 

\noindent Nous allons contourner cet obstacle en tronquant la s\'erie double \`a un ordre fini et en contr\^olant les ``effets de bord" dus \`a l'interversion des sommations.

 \bigskip\noindent  {\bf 7.3. Preuve de 7.1.1.} Effectuons une suite d'\'eclatements ponctuels  
$\, \pi_1 : {  Y}_1\to {  Y}\,$  (en commen\c cant par \'eclater $Q$), de sorte que 
 
\smallskip  \noindent $i)$ $\pi_1^{-1}({  Z}\cup {  C})$ soit \`a croisements normaux, 

\smallskip \noindent $ii)$ chaque composante irr\'eductible du diviseur exceptionnel ne coupe qu'un au plus des  transform\'es  stricts  de ${  Z}_i$,

\smallskip \noindent $iii)$ la composante du diviseur exceptionnel que coupe le transform\'e strict $C_1$ de $C$ ne contienne aucun point tournant  pour ${\sN}_1= \pi_1^\ast({\sN})$ (c'est loisible en vertu de 4.3.3),  

\smallskip \noindent $iv)$ tout point de croisement du diviseur exceptionnel soit semi-stable pour ${\sN}_1 $  (c'est loisible en vertu de 5.4.1). 

\smallskip Notons ${\Bbb T}_1$ l'arbre associ\'e, comme en 6.5, au diviseur exceptionnel. 
  Indexons les composantes $E_{1,v}$ du diviseur exceptionnel par les sommets $v$ de ${\Bbb T}$.  Notons $v_0$ le sommet cor\-res\-pondant  \`a la composante $E_{1,v_0}$ que coupe  $C_1$.    
 
\smallskip   Notons ${V}$ l'ensemble des sommets $v\in  {\Bbb T}_1$ pour lesquels la composante  $E_{1,v}$ coupe l'un - et un seul par $ii)$,  not\'e $Z_{1,v}$  - des transform\'es stricts des ${  Z}_i$.  

\smallskip   Notons ${V'}$ l'ensemble des sommets $v\in  {\Bbb T}_1$ pour lesquels $E_{1,v}$ a des points non semi-stables. D'apr\`es $i)$ et $iii)$, ni $V$ ni ${V'}$ ne contiennent $v_0$. 

\smallskip Consid\'erons une suite d'\'eclatements ponctuels
$$  Y_n\to \cdots \to Y_2\to Y_1 $$ o\`u l'on commence par \'eclater successivement les points non semi-stables des transform\'es stricts des $E_v,\, v\in V'$ (ce qui cor\-res\-pond \`a ajouter \`a chaque fois une ar\^ete \`a l'arbre partant de l'un des $v\in V'$),  puis, successivement, les points de croisement pour les semi-stabiliser (ce qui cor\-res\-pond \`a subdiviser les ar\^etes qu'on vient d'ajouter). Pour $m\leq n$, notons $\pi_m: Y_m\to Y$ la suite d'\'eclatements compos\'ee et $\sN_m$ l'image inverse de $\sN$ sur  $Y_m$.

\smallskip Notons ${\Bbb T}_m$ l'arbre associ\'e au diviseur exceptionnel  $\pi_m^{-1}(Q)$,  $E_{m,v} $ la composante de $\pi_m^{-1}(Q)$ attach\'ee \`a $v\in {\Bbb T}_m$,  $Z_{m,v}$ le transform\'e strict de $Z_{1,v}$ (pour $v\in V$), $C_{m }$ le transform\'e strict de $C_{1 }$, et $Q_m = C_m \cap E_{m, v_0}$. 
 Observons que ${\Bbb T}_1 $  est un sous-arbre de ${\Bbb T}_m$, que les sommets de ${\Bbb T}_m$ s'identifient \`a des sommets de ${\Bbb T}_{m+1}$, mais que ${\Bbb T}_m$ n'est pas n\'ecessairement un sous-arbre de ${\Bbb T}_{m+1}$. 
\smallskip Notons $A_m$ la matrice (sym\'etrique) d'intersection des $E_{m,v},\,  v\in {\Bbb T}_m$, affect\'ee du signe $-$. Rappelons que ses coefficients non nuls non diagonaux valent $-1$, et que son inverse est \`a coefficients positifs. 

Notons  $B_m$ la matrice (antisym\'etrique) ayant pour coefficients les entiers
$$(B_m)_{v w}=(A_m)^{-1}_{v_0v}(A_m)_{v w}.{\rm ir}_{E_{m,w}}({\sN_m}) - (A_m)^{-1}_{v_0w}(A_m)_{wv}.{\rm ir}_{E_{m,v}}({\sN_m}) ,\; v,w\in {\Bbb T}_m .$$

\proclaim Lemme 7.3.1. $i)$ Si $w\in {\Bbb T}_{m+1}\setminus {\Bbb T}_{m }$ cor\-res\-pond \`a l'\'eclatement d'un point non semi-stable de $E_{m, v},\, v\in {V'}$, alors $$(A_{m+1})^{-1}_{v_0w}= (A_m)^{-1}_{v_0v},$$
$${\rm ir}_{E_{m+1, w}}({\sN_{m+1}})\leq {\rm ir}_{E_{m,v}}({\sN_m}),$$
$$(B_{m+1})_{vw}\geq 0.$$
 $ii)$ Si $w\in {\Bbb T}_{m+1}\setminus {\Bbb T}_{m }$ cor\-res\-pond \`a l'\'ecla\-te\-ment du point d'intersection $E_{m,v}\cap E_{m,v'},\; v,v'\in {\Bbb T}_{m }$, alors 
$$(A_{m+1})^{-1}_{v_0w}= (A_{m })^{-1}_{v_0v}+ (A_{m})^{-1}_{v_0v'},$$
$${\rm ir}_{E_{m+1, w}}({\sN_{m+1}})\leq {\rm ir}_{E_{m,v}}({\sN_m})+ {\rm ir}_{E_{m,v'}}({\sN_m}),$$
$$(B_{m+1})_{vw}\geq (B_{m })_{vv'}.$$\par

 \dem Les assertions sur $A_.^{-1}$, dont les coefficients sont interpr\'et\'es comme multiplicit\'es d'intersection, sont laiss\'ees au lecteur (voir par exemple [${\bf S2}$, I.3.2.9, I.3.2.10]).
Les assertions sur l'irr\'egularit\'e d\'ecoulent de (7.2.1.3)$_{t=1}$ et (7.2.1.4). Ces assertions impliquent celles sur $B_. $. \sq 
 
\smallskip  En par\-ti\-cu\-lier, on a $(B_{n})_{vw}\geq 0$ pour tout  $v\in V',\,w\in  {\Bbb T}_n\setminus {\Bbb T}_1$. 
\par\noindent Par antisym\'etrie, on a   
$$\sum_{v\in {\Bbb T}_1, w\in {\Bbb T}_n }\, (B_n)_{v,w} = \sum_{(v,w), v\in {V'}, w\in {\Bbb T}_n\setminus {\Bbb T}_1}\, (B_n)_{v,w}, \leqno {(7.3.1.1)}$$
qui est donc une quantit\'e positive. 

\smallskip L'in\'egalit\'e (6.2.4.2) appliqu\'ee \`a $\bar S = E_{n,v}$ (et \`a la suite d'\'eclatements formels $\hat\pi_n$ associ\'ee \`a $\pi_n$, \cf 6.2.1) donne 
$$\sum_{ w\in {\Bbb T}_n }\,  (A_n)_{vw}.{\rm ir}_{E_{n,w}}({\sN_n})\leq  {\rm ir}_{Z_{n,v}}({\sN_n})  $$ ou $\leq 0$, selon que $v\in V$ ou non. D\`es lors
$$\sum_{v\in {\Bbb T}_1, w\in {\Bbb T}_n }\, (A_n)^{-1}_{v_0v}(A_n)_{vw}.{\rm ir}_{E_{n,w}}({\sN_n})\leq \sum_{v\in V}  (A_n)^{-1}_{v_0v}.{\rm ir}_{Z_{n,v}}({\sN_n}). \leqno {(7.3.1.2)}$$ Par ailleurs 
$$\sum_{v\in {\Bbb T}_1, w\in {\Bbb T}_n}\,   (A_n)^{-1}_{v_0w}(A_n)_{wv}.{\rm ir}_{E_{n,v}}({\sN_n}) = {\rm ir}_{E_{n,v_0}}({\sN_n}).\leqno {(7.3.1.3)}$$   
 La diff\'erence $$ \sum_{v\in {\Bbb T}_1, w\in {\Bbb T}_n }\, (A_n)^{-1}_{v_0v}(A_n)_{vw}.{\rm ir}_{E_{n,w}}({\sN_n})  - 
  \sum_{v\in {\Bbb T}_1, w\in {\Bbb T}_n}\,   (A_n)^{-1}_{v_0w}(A_n)_{wv}.{\rm ir}_{E_{n,v}}({\sN_n}) ,$$ qui n'est autre que $ \sum_{v\in {\Bbb T}_1, w\in {\Bbb T}_n }\, (B_n)_{v,w}$, est donc positive, ce qui donne
  $${\rm ir}_{E_{n,v_0}}({\sN_n})\leq \sum_{v\in V}  (A_n)^{-1}_{v_0v}.{\rm ir}_{Z_{n,v}}({\sN_n})= \sum_{i} \,  ({ C},{  Z}_i)_Q\,. \, {\rm ir}_{{  Z}_i}(\sN )  .\leqno {(7.3.1.4)}$$  
 Enfin, puisque $E_{n,v_0}$ n'a pas de point tournant, on a $\,{\rm ir}_{E_{n,v_0}}({\sN_n})=  {\rm ir}_{Q_n}{\sN_{ \mid C_n} }= {\rm ir}_{O}{\sN_{\tilde C} }$. Combin\'e \`a (7.3.1.4), ceci prouve (7.1.1.1).
\sq

 \bigskip \bigskip

\bigskip {\bf 8. Autour de la conjecture de Sabbah.  }  

\medskip\noindent  {\bf 8.1. \'Enonc\'e.} Soit $k$ un corps alg\'ebriquement clos de carac\-t\'e\-ris\-tique nulle. Pla\c cons-nous de nouveau dans la si\-tua\-tion de 4.1: $\frak M$ est une con\-ne\-xion m\'ero\-mor\-phe formelle sur $\frak X \cong \Spf A[[{x}]]$ \`a p\^oles le long de ${x}=0$.  

Supposons d'abord que $Z= \Spec A$ soit une courbe connexe lisse sur $k$. 

D'apr\`es 4.3.3, en \'eclatant successivement les points tournants qui apparaissent, on aboutit au bout d'un nombre fini d'\'etapes \`a la si\-tua\-tion o\`u tous les points de non-croisement du diviseur polaire sont stables pour l'image inverse de la con\-ne\-xion. Mais les points de croisement pourraient n'\^etre pas stables, ni m\^eme semi-stables.

Parall\`element, d'apr\`es 5.4.1 et 5.5.3, en \'eclatant successivement les points de croisement qui apparaissent, on aboutit au bout d'un nombre fini d'\'etapes \`a la si\-tua\-tion o\`u tous les points de croisement du diviseur polaire sont stables pour l'image inverse de la con\-ne\-xion. Mais il pourrait y avoir des points tournants, et m\^eme des points non semi-stables. 

La conjecture de Sabbah pr\'edit qu'on peut mener \`a bien les deux processus \`a la fois:

\proclaim Conjecture 8.1.1. {\rm (Sabbah)}. Il existe une suite finie $\pi: {\frak X}'\to {\frak X}$ d'\'eclatements ponctuels formels telle que tous les points de $  \pi^{-1}(Z)$ soient stables pour $\pi^\ast \frak M$.\par

\medskip\noindent  {\bf 8.2. Connexions rigides absolues.} \`A pr\'esent, supposons seulement que $Z$ soit une vari\'et\'e affine connexe lisse sur $k$. Au sch\'ema formel $\frak X $, on peut associer sa {\it fibre g\'en\'erique} au sens de Raynaud $X = {\frak X}_\frak k$, qui est une vari\'et\'e analytique rigide (affino\"{\i}de) sur $\frak k= k(({x}))$. Ses points sont les germes de courbes sur $\frak X$ (\cf [$\bf BoL$]). On obtient une vari\'et\'e analytique rigide canoniquement isomorphe en partant d'un \'eclat\'e formel ponctuel ${\frak X}'$ de $\frak X$, en identifiant un germe de courbe sur $\frak X$ et son transform\'e strict sur ${\frak X}'$; un tel \'eclat\'e est dit {\it mod\`ele formel} de $X$. On a un morphisme de sites annel\'es $ \, sp: \, X \to {\frak X}'\, $ (dit de sp\'ecialisation)
 et le tube $sp^{-1}(S)$ de
 tout sous-sch\'ema localement ferm\'e $S$ de $Z' = {\frak X}'_{red}$ est un ouvert admissible de $X$ (\cf [$\bf Ber$, 0.2]). 
 
 \medskip \noindent {\bf Exemple 8.2.1.} Si $Z\subset {\bf A}^1 $, $X$ est le compl\'ementaire d'un nombre fini de disques unit\'e ouverts dans le disque unit\'e ferm\'e. Si ${\frak X}'$ est un \'eclatement torique, et $Q$ est le point de croisement (\resp $T$ est la r\'eunion des deux composantes de $Z'$) correspondant au c\^one d'ar\^etes $(a,c), (b,d)$ de l'\'eventail associ\'e, le tube de $Q$ (\resp $\bar T$) est la couronne ouverte (\resp ferm\'ee) de rayons $\vert {x}\vert^{d/b}$ et  $\vert {x}\vert^{c/a}$.  
  
 \medskip \noindent Toute connexion m\'eromorphe formelle $\frak M$ \`a p\^oles le long de $Z$ fournit un module localement libre $M$ sur $X$ muni d'une connexion ``absolue"  $ \,M \to M\otimes \Omega^1_{X/k} \;$ (o\`u $\Omega^1_{X/k}\,= \Omega^1_{X/\frak k}\oplus \sO_X.dx$). 
 Une telle connexion est dite {\it r\'eguli\`ere} si $M$ l'est (il revient au m\^eme de dire que l'image inverse de $M$ sur l'un des mod\`eles formels l'est). Elle est dite {\it \'el\'ementaire} si elle est produit tensoriel d'une connexion de rang un et d'une connexion r\'eguli\`ere.   

\proclaim Th\'eor\`eme 8.2.2. Il existe un recouvrement ouvert $(U_\alpha)$ de $X$, et des morphismes \'etales finis $V_\alpha\to U_\alpha$ tels que l'image inverse de $M$ sur $V_\alpha$ soit somme directe de connexions \'el\'ementaires. \par
  
 \dem \'Etant donn\'e un point de $X$, vu comme germe de courbe sur un mod\`ele formel $\frak X'$, on se ram\`ene, par \'eclatement et en passant \`a une extension finie de $\frak k$, \`a supposer que le germe est lisse et coupe transversalement $Z' = {\frak X'}_{red}$. D'apr\`es 4.3.1, on peut m\^eme supposer qu'il coupe $Z'$ dans un ouvert lisse connexe $S\subset Z'$ form\'e de points semi-stables. On peut prendre pour $U_\alpha$ le tube de $S$. On couvre ainsi $X$.  
 \sq
 
 \proclaim Conjecture 8.2.3.  Il existe un tel recouvrement  $(U_\alpha)$  {\rm fini}. \par
 
 \medskip\noindent {\bf Remarques 8.2.4.} $1)$ Dans le cas o\`u $X$ est un point, 8.2.2 \'equivaut au th\'eor\`eme de Turrittin-Levelt sur $\frak k$.
 
 \smallskip \noindent  $2)$ Dans la situation de 8.1 (o\`u $X$ est une courbe affino\"{\i}de), la conjecture 8.2.3 d\'ecoule de celle de Sabbah: en effet, dans un mod\`ele formel ${\frak X}'$ dans lequel tous les points de $\pi^{-1}(Z)$ sont semi-stables, indexons par $\alpha$ les croisements $Q'$ de $Z'= \pi^{-1}(Z)$ ainsi que les composantes $S$ de la partie lisse de $\pi^{-1}(Z)$; on peut alors prendre pour $U_\alpha$ les tubes de ces sous-sch\'emas de $Z'$.  
 
 \smallskip \noindent $3)$ Dans 8.2.2, on ne doit pas s'attendre \`a trouver un recouvrement admissible de $X$, en raison de ce que les d\'ecompositions formelles aux croisements $Q'$ ne correspondent pas en g\'en\'eral aux d\'ecompositions de Turrittin-Levelt sur les strates  $S$ voisines, mais \`a des raffinements d'icelles, \cf 5.3.
 
 \medskip\noindent  {\bf 8.3. Probl\`eme de compacit\'e.} Il est tentant d'essayer de prouver 8.2.2 par un argument de quasi-compacit\'e, par exemple en consid\'erant l'espace de Zariski-Riemann $X^{ZR} = \limp {\frak X}'$. Il s'av\`ere plus commode de passer au plus grand quotient s\'epar\'e de $X^{ZR}$ (\cf [$\bf F$]), qui est l'espace $\frak k$-analytique de Berkovich $X^{an}$ associ\'e \`a $X$. On a des morphismes de sites annel\'es $ \,  X \to X^{ZR} \to X^{an}$,  le compos\'e est injectif et d'image dense, et $X^{an}$ est compact. 
 
 On peut encore attacher \`a $\frak M$ un module localement libre \`a connexion absolue sur $X^{an}$, et  essayer de traiter 8.2.3 en prouvant un analogue de 8.2.2 pour $X^{an}$.
  
 Dans la situation 8.2.1, les points de $X^{an}$ sont tr\`es faciles \`a d\'ecrire (\cf [$\bf Ber$, 3.6]); on constate que ceux qui posent probl\`eme sont les points \`a valeurs, non pas dans le compl\'et\'e $\hat{\bar{\frak k}}$ de la cl\^oture alg\'ebrique de $\,\frak k$, mais dans son extension maximalement compl\`ete $\displaystyle{\,\{ \sum_{r\in \Omega }\, a_r x^r \}\,}$ (o\`u $\,\Omega\,$ parcourt les sous-ensembles bien ordonn\'es de $\,\Q$); les arguments du style 4.3.2 s'av\`erent en effet inop\'erants lorsque le premier point d'accumulation des $\,\Omega\subset \Q_{>0}\,$ qui interviennent est inf\'erieur au rang de Poincar\'e-Katz de $\,\frak M$.

  \bigskip\bigskip{\bf Appendice: sp\'ecialisation du polyg\^one de Newton. }
 
 \medskip On se place dans le cadre et les hypoth\`eses de 3.2: $A$ est noeth\'erien int\'egrale\-ment clos, $M$ est un module diff\'erentiel de rang $\mu$ sur $A((x))$, et $\mu$ est strictement inf\'erieur aux caract\'eristiques r\'esiduelles de $A$ si celles-ci sont non nulles. 
 
 Rappelons que nos polyg\^ones de Newton sont plac\'es, par convention, de mani\`ere \`a ce que le sommet le plus \`a droite soit le point de coordonn\'ees $(\mu, 0)$.
 
  \proclaim Th\'eor\`eme A.1.  Pour tout point $P$ de $Z= \Spec A$, 
  $NP(M_{(P)})\subset NP(M)$.\par
 
  \proclaim Corollaire A.2.  Soient $X$ et $Y$ des vari\'et\'es alg\'ebriques lisses sur $k$ (\resp des vari\'et\'es analytiques complexes). Soit $f: Y\to X$ un morphisme lisse, de dimension relative $1$, \`a fibres connexes,  et soit $Z$ une hypersurface de $Y$  finie et {\rm \'etale} sur $X$ (via $f$). Soit ${\sN}$ un module \`a  con\-ne\-xion m\'ero\-morphe relative sur $Y$ \`a p\^oles le long de $Z$ (\cf 1.1). 
  Alors le bord du polyg\^one de Newton  $NP_z({\sN}_{(f(z))})\;$  
 ne peut que cro\^{\i}tre (au sens large) par sp\'ecialisation sur $Z$.
\par

 Pour d\'emontrer A.1, on peut supposer et on supposera que $Z$ est un trait dont $P$ est le point ferm\'e. Alors $M$ est libre de rang $\mu$. On note $K$ le corps de fractions de $A$, $k$ le corps r\'esiduel.   
 Notons $\mu_{(\sigma)}$ (\resp ${\bar\mu}_{(\sigma)}$) le rang de la partie de pente $ \sigma $ de $M$ (\resp de $M_k = M_{(P)}$). En outre, quitte \`a effecture une ramification, on peut supposer que toutes les pentes (tant de $M$ que de $M_k$) sont enti\`eres.

 En suivant [$\bf GL$], choisissons un  $A[[x]]$-r\'eseau $\Bbb M$ de $M$ et consid\'erons la suite double de sous-$A[[x]]$-modules de $M$ (index\'ee par $(m,n)\in \N^2$):  
 $${\Bbb M}_{m,0} = {\Bbb M},\;\; {\Bbb M}_{m,n+1}= {\Bbb M}_{m, n }+ x^{m +1 }  {{d}\over{dx}}({\Bbb M}_{m, n }).$$  
 Notons que les quotients ${\Bbb M}_{m,n}/\Bbb M$ sont des $A$-modules de type fini.

  \proclaim Lemme A.3. On a 
    $$\lambda_{m } := \lim_n {1\over n}\dim_K  \,({\Bbb M}_{{m } ,n })_K/{\Bbb M}_K = \sum_{\sigma>m} (\sigma -m  )\mu_{(\sigma)} , \leqno (A.3.1) $$
   $${\bar\lambda}_{m }:= \lim_n {1\over n}\dim_k \,({\rm Im}\, ({\Bbb M}_{{m },n})_k \to M_k/{\Bbb M}_k) = \sum_{\sigma>m}(\sigma -m  ){\bar\mu}_{(\sigma)}  \leqno (A.3.2)  $$  et 
   $$\lambda_{m }\,\geq \, {\bar\lambda}_{m }\, \geq \, 0.  \leqno (A.3.3) $$      \par
  La valeur des limites se d\'eduit ais\'ement de 2.1.3.ii et 2.3.1 (voir aussi  [$\bf GL$], qui explicite l'\'egalit\'e $\lambda_0= {\rm ir}\, M$). L'in\'egalit\'e entre  limites est essentiellement le lemme de [$\bf De2$]\footnote{$^{(12)}$}{\sm cela r\'esulte formellement de ce que le double $A[[x]]$-dual $ {\Bbb M}_{m,n}^{\ast\ast}$ est libre sur $A[[x]]$ et que le quotient $ {\Bbb M}_{m ,n}^{\ast\ast}/{\Bbb M}$ est plat sur $A$; en effet, on a $\dim_k \,({\rm Im}\, ({\Bbb M}_{m ,n})_k \to M_k/{\Bbb M}_k)= \dim_k ({\Bbb M}_{m ,n}^{\ast\ast}/{\Bbb M})_k - \dim_k ({\Bbb M}_{m ,n}^{\ast\ast}/{\Bbb M}_{m ,n})_k$.}. \sq
   
\smallskip Pour $x\in [0,\mu]$, notons $f(x)$ (\resp ${\bar f}(x)$) la fonction convexe, affine par morceaux, \`a valeurs n\'egatives, dont le graphe borde $NP(M)$ (\resp $NP(M_k)$). Pour $\xi\in [0, \infty[$, notons $\displaystyle f^\ast(\xi)= \sup_x\, (x\xi - f(x))$ (\resp ${\bar f^\ast}(\xi)$) la transform\'ee de Legendre de $f$ (\resp $\bar f$). Le couple $(f, f^\ast)$ v\'erifie l'in\'egalit\'e de Young 
 $\, x\xi  \leq f(x)+ f^\ast(\xi) $. 
 
 Compte tenu de (A.3.1) et (A.3.2), on a 
 $$ f^\ast(m )= m\mu+ \lambda_m,\;  {\bar f^\ast}(m )= m\mu+ {\bar \lambda}_m  \leqno (A.3.4)$$ pour tout $m\geq 0$, et l'in\'egalit\'e de Young pour $(f,f^\ast)$ donne 
 $$\lambda_m\geq  - f(x) - m(\mu-x)  \leqno (A.3.5)$$ pour tout $x\in [0,\mu]$.

 Par ailleurs, pour $x = x_m <\mu$ \' egal \`a l'abscisse d'un sommet de $NP(M_k)$, et pour $\xi = m $ \'egal \`a la pente de $NP(M_k)$ \`a proximit\'e droite de $x$, l'in\'egalit\'e de Young pour $({\bar f}, {\bar f^\ast})$ est une \'egalit\'e et s'\'ecrit 
$$ {\bar\lambda}_m= -{\bar f}(x_m)- m(\mu-x_m).\leqno (A.3.6)$$ 
On d\'eduit enfin de (A.3.3)(A.3.5)(A.3.6) l'in\'egalit\'e 
$${\bar f}(x_m)\geq f(x_m), \leqno (A.3.7)$$ qui prouve le th\'eor\`eme.\sq

\smallskip \noindent {\bf Remarque A.4.} On voit facilement aussi que la fonction dont le graphe borde $NP(M)$ n'est autre que la transform\'ee de Legendre, sur $[0,\mu]$, de la fonction affine par morceaux $f^\ast$ qui interpole (A.3.4) sur $[0,\infty[$ ($\lambda_m$ \'etant d\'efini dans (A.3.1)).
  
\bigskip

  \centerline{\it Bibliographie}
 \bigskip
  \item{[$\bf{A1}$]} Andr\'e Y., 
Diff\'erentielles non commutatives et th\'eorie de Galois
diff\'erentielle ou aux diff\'erences, Ann. Scient. E.N.S., {\bf 5} (2001), 1-55.
\medskip \item{[$\bf{A2}$]}  Andr\'e Y., Filtrations de type Hasse-Arf et monodromie $p$-adique,
  Invent. Math.  {\bf 148} (2002),  285-317.
\medskip \item{[$\bf{A3}$]} Andr\'e Y., An algebraic proof of Deligne's regularity criterion, \`a para\^{\i}tre dans les actes du colloque du R.I.M.S. Kyoto 2006. 
  \medskip
 \item{[$\bf{AB}$]} Andr\'e Y., Baldassarri F., {\it De Rham Cohomology of
Differential Modules on Algebraic Varieties},
Progress in Mathematics, Vol. {\bf 189}, Birkh\"auser (2001). Deuxi\`eme \'edition enti\`erement refondue en pr\'eparation.
\medskip
\item{[$\bf{BaV}$]} Babbitt D., Varadarajan V., {Deformation of nilpotent matrices over rings and reduction of analytic families of meromorphic differential equations}, Mem. Amer. Math. Soc. {\bf 55}, n. 325 (1985). 
 \medskip
\item{[$\bf{B}$]} Baldassarri F., Towards an algebraic proof of Deligne's regularity criterion. An informal survey of open problems,  Milan J. Math.  {\bf 73}  (2005), 237--258.
 \medskip
\item{[$\bf{Be}$]} Berkovich V., \'Etale cohomology for non-Archimedean analytic spaces, Publ. math. Inst. Hautes ƒtudes Sci. {\bf 78} (1993), 5--161.
 \medskip
\item{[$\bf{Ber}$]} Berthelot P., Cohomologie rigide et cohomologie rigide \`a supports propres, pr\'epublication (1996).
\medskip
 \item{[$\bf{BoL}$]} Bosch S., L\"utkebohmert W., Formal and rigid geometry I, Math. Ann. 295 (1993), 291-317.
  \medskip
 \item{[$\bf{CD}$]} Christol G., Dwork B., Modules diff\'erentiels sur des couronnes,
 Ann. Inst. Fourier {\bf 44} (3) (1994), 663--701.
\medskip
 \item{[$\bf{ De1}$]} Deligne P.,  {\it \'Equations diff\'erentielles \`a points singuliers r\'egu\-liers}, Springer Lecture Notes in Math. {\bf 163} (1970) + erratum.
\medskip
 \item{[$\bf{ De2}$]} Deligne P.,  Lettre \`a N. Katz  (1/12/1976), \`a para\^{\i}tre dans un volume de Documents math\'ematiques, S.M.F.
\medskip
\item{[$\bf{DGS}$]} Dwork B., Gerotto G., Sullivan F., {\it An Introduction to
{\rm G}-functions}, Annals of Mathematical Studies {\bf 133},
Princeton University (1994).
\medskip
\item{[$\bf{E}$]} Euler L., De seriebus divergentibus, in: Opera omnia I.14 Teubner (1925).
\medskip
 \item{[$\bf{F}$]}  Fujiwara K., Theory of tubular neighborhoods in \'etale topology, Duke Math. J. {\bf 80} (1995),  15--57. 
\medskip
 \item{[$\bf{GL}$]} G\'erard R., Levelt A., Invariants mesurant l'irr\'egularit\'e en un point singulier des syst\`emes d'\'equations diff\'erentielles lin\'eaires, Ann. Inst. Fourier (Grenoble) {\bf 23} (1973), no. 1, 157--195.
  \medskip  \item{[$\bf{Ka}$]} Katz N., {Nilpotent connections and the monodromy theorem.
Applications of a result of Turrittin},
Publ. Math. IHES {\bf 39} (1970), 175-232.
\medskip  \item{[$\bf{Ke}$]} Kedlaya S., Semi-stable reduction for overconvergent $F$-isocrystals, I, II, III, pr\'epublications (2006). 
  \item {[$\bf{LvdE}$]} Levelt A., van den Essen A., {Irregular singularities in several variables}, Mem. Amer. Math. Soc. {\bf 40}, n. 270 A.M.S. (1982).
 \medskip
\item {[$\bf{L}$]} Levelt A., Jordan decom\-po\-si\-tion for a class of singular
differential operators,
Ark. Math. {\bf 13} (1975), 1-27.
 \medskip
\item {[$\bf{M}$]} Malgrange B., Connexions m\'ero\-morphes, II: le r\'eseau canonique, Invent. Math. {\bf 124} (1996), 367-387.
  \medskip
\item {[$\bf{Me1}$]} Mebkhout Z., Sur le th\'eor\`eme de semi-continuit\'e de l'irrŽgularit\'e des \'equations diff\'erentielles, in: {\it Differential systems and singularities} (Luminy, 1983). Ast\'erisque {\bf 130} (1985), 365--419.
 \medskip
\item {[$\bf{Me2}$]} Mebkhout Z., Le th\'eor\`eme de positivit\'e de l'irr\'egularit\'e pour les $D_X$-modules, in: {\it The Grothendieck Festschrift}, Vol. III, 83--132, Progr. Math. {\bf 88} (1990) Birkh\"auser. 
  \medskip
 \item {[$\bf{S1}$]} Sabbah C., {\'Equations diff\'erentielles \`a points singuliers irr\'egu\-liers en dimension $2$ }, Ann. Inst. Fourier {\bf 43} 5 (1993) 1619-1688.
  \medskip
 \item {[$\bf{S2}$]} Sabbah C., {\it \'Equations diff\'erentielles \`a points singuliers irr\'egu\-liers et ph\'eno\-m\`ene de Stokes en dimension $2$},
Ast\'erisque {\bf 263} (2000).
  \medskip \item {[$\bf{Sc}$]}  Sch\"afke R., Formal fundamental solutions of irregular singular differential  equations depending upon parameters.  J. Dynam. Control Systems  {\bf 7}  (2001),  no. 4, 501--533. 
 \medskip
 \item {[$\bf{Se}$]} Serre J.-P., {\it Corps locaux}, Hermann, 1968. 
 \medskip
 \item {[$\bf{T}$]} Turrittin H., Convergent solutions of ordinary differential equations in the neighborhood of an irregular point, Acta Math. {\bf 93} (1955), 27-66.
  \medskip
 \item {[$\bf{W}$]}ÊWasow W., {\it Linear turning point theory}, Applied Mathematical Sciences, {\bf 54}, Springer-Verlag (1985).
   \medskip
 \item {[$\bf{SGA1}$]} Rev\^etements \'etales et groupe fondamental, S\'eminaire de g\'eom\'etrie alg\'ebrique du Bois Marie 1960--61. Dirig\'e par A. Grothendieck.  Documents Math\'ema\-tiques  {\bf 3}. Soci\'et\'e Math\'ematique de France, Paris (2003).

    \bigskip\bigskip  
    
    \sm{D\'epartement de Math\'ematiques,   \'Ecole Normale Sup\'erieure, }
   
   \sm{45 rue d'Ulm, F-75005 Paris. Courriel: yves.andre@ens.fr}

 \end